\documentclass[11pt ,reqno]{amsart}
\usepackage{amssymb,amscd,verbatim, amsthm}
\usepackage[dvips]{color}
\usepackage[mathscr]{eucal}

\newcommand\uset{\underset}
\newcommand\ub{\underbrace}
\newcommand\ovl{\overline}

\newcommand \und[2]{\underset{#1}{#2}}
\newcommand \ovr[2]{\overset{#1}{#2}}
\newcommand \unb[2]{\underset{#1}{{\underbrace{#2}}}}

\newcommand \fk[1]{{{\mathfrak #1}}}
\newcommand \C[1]{{\mathcal #1}}
\newcommand \ov[1]{{\overline {#1}}}
\newcommand \ch[1]{{\check{#1}}}
\newcommand \bb[1]{{\mathbb #1}}

\newcommand\lr{\overset{LR}{<}}
\newcommand\elr{\overset{LR}{\thickapprox}}

\newcommand \wti[1]{{\widetilde {#1}}}
\newcommand \wht[1]{{\widehat {#1}}}

\newcommand \bA{{\bb A}}
\newcommand \bC{{\bb C}}
\newcommand \bF{{\bb F}}
\newcommand \bH{{\bb H}}
\newcommand \bN{{\bb N}}
\newcommand \bR{{\bb R}}
\newcommand \bZ{{\bb Z}}

\newcommand\cha{{\check \alpha}}

\newcommand\cfm{{\check{\fk m}}}

\newcommand\one{1\!\!1}

\newcommand\CI{{\C I}}
\newcommand\CO{{\C O}}
\newcommand\CR{{\C R}}

\newcommand \vA{{\check A}}
\newcommand \va{\check{\mathfrak a}}
\newcommand \vO{{\check \CO}}

\newcommand\cCO{{\check{\CO} }}

\newcommand\ie{{\it i.e.~ }}

\newcommand\eg{{\it e.g.~ }}

\newcommand\ep{{\epsilon}}
\newcommand\la{{\lambda}}
\newcommand\La{{\Lambda}}
\newcommand\om{{\omega}}
\newcommand\al{{\alpha}}
\newcommand\sig{{\sigma}}

\newcommand \vG{{\check G}}

\newcommand \vg{\check{\fk g}}
\newcommand \vm{\check{\fk m}}
\newcommand \ve{{\check e}}
\newcommand \vch{{\check h}}
\newcommand \vf{{\check f}}

\newcommand \vh{{\check h}}

\newcommand \fz{\mathfrak z} 

\newcommand\Hom{\operatorname{Hom}}
\newcommand\Ad{\operatorname{Ad}}
\newcommand\Img{\operatorname{Im}}
\newcommand\ad{\operatorname{ad}}
\newcommand\rk{\operatorname{rk}}
\newcommand\Ind{\operatorname{Ind}}
\newcommand\ind{\operatorname{ind}}

\newcommand\fg{\mathfrak g}

\newtheorem*{corollary}{Corollary}
\newtheorem*{definition}{Definition}
\newtheorem*{lemma}{Lemma}
\newtheorem*{proposition}{Proposition}
\newtheorem*{remark}{Remark}

\newtheorem*{theorem}{Theorem}

\numberwithin{equation}{subsection}

\begin{document}
\today

\bigskip
\title{UNITARY SPHERICAL SPECTRUM FOR SPLIT CLASSICAL GROUPS}

\author{Dan Barbasch}
       \address[D. Barbasch]{Dept. of Mathematics\\
                Cornell University\\Ithaca, NY 14850}
        \email{barbasch@math.cornell.edu}

\maketitle

\bigskip
\section{Introduction}\label{sec:1}

This paper gives a complete classification of the spherical unitary 
dual of the split groups
$Sp(n)$ and $So(n)$ over the real and $p-$adic field. 
{
While there is a large overlap with  \cite{B1},\ \cite{B2},\ \cite{B3} and
\cite{BM3},  the formulation of the answer is new, and so are a
significant number of the techniques employed. Completely new is the
proof of necessary conditions for 
unitarity in the real case. As is explained in these references, in
general, in the $p-$adic case, the classification of the spherical
unitary dual is equivalent to the classification of the unitary
irreducible finite dimensional representations of an affine graded
Hecke algebra.   In concrete terms the problem then becomes to
determine whether a certain hermitian form for each Weyl group
representation is positive semidefinite.    
For the real case, following a suggestion of D. Vogan, I find a set of
$K-$types which I call \textit{relevant} which detect the
nonunitarity. They have the property that they are in 1-1
correspondence with certain irreducible Weyl group 
representations (also called relevant) so that the hermitian forms
are \textit{the same} in the real and $p-$adic case. The fact that
these relevant $W-$types detect unitarity in the $p-$adic case is also
new. Thus the same proof applies in both cases.}
Since the answer is independent of the field, this establishes a form
of the Lefschetz principle.
 
\bigskip
Let $G$ be a split symplectic or orthogonal group over a local field
$\bF$ which is either $\bR$ or a $p-$adic field. 
Fix a maximal compact subgroup $K.$ In the real case,
there is only one conjugacy class. In the $p-$adic case, let $K=G(\C
R)$ where $\bF\supset\C
R\supset\C P,$ with $\C R$ the ring of integers and $\C P$ the
maximal prime ideal. Fix also a  rational Borel
subgroup $B=AN.$ Then $G=KB$. A
representation $(\pi,V)$ (admissible) is called spherical if $V^K\ne
(0).$ 

\medskip
The classification of irreducible admissible spherical modules
is well known. For every irreducible spherical $\pi,$
there is a character $\chi\in\widehat{A}$ such that
$\chi |_{A\cap K}=triv,$ and $\pi$ is the unique spherical subquotient
of $Ind_{B}^G[\chi\otimes\one].$ We will call a character $\chi$ whose
restriction to  $A\cap K$ is trivial, \textit{unramified}. Write
$X(\chi)$ for the induced module (principal series) and $L(\chi)$ for
the irreducible spherical subquotient. Two such modules $L(\chi)$ and
$L(\chi')$ are equivalent if and only if there is an element in the
Weyl group $W$ such that $w\chi=\chi'.$ An $L(\chi)$ admits a
nondegenerate hermitian form if
and only if there is $w\in W$ such that $w\chi=-\ovl{\chi}.$ 

The character $\chi$ is called \textit{real} if it takes only positive
real values. For real groups, $\chi$ is real if and only if
$L(\chi)$ has real infinitesimal character (\cite{K}, chapter 16). As
is proved there, any unitary representation of a real reductive group
with nonreal infinitesimal character is unitarily induced from a
unitary representation with real infinitesimal character on a proper
Levi component. So for real groups it makes sense to consider  only
real infinitesimal character. In the $p-$adic case, $\chi$ is called
real if the infinitesimal character is real in the sense of \cite{BM2}.
The results in \cite{BM1} show that the problem of determining the
unitary irreducible representations  with Iwahori fixed vectors is
equivalent to the same problem for the Iwahori-Hecke algebra. In
\cite{BM2}, it is shown  that the problem of classifying the unitary dual 
for the Hecke algebra reduces to determining the unitary dual with
real infinitesimal character of some smaller Hecke algebra (not necessarily
one for a proper Levi subgroup). So for $p-$adic groups as well it is
sufficient to consider only real $\chi.$  

So we start by parametrizing real unramified characters of $A.$ 
Since $G$ is split, $A\cong (\bF^\times)^n$ where $n$ is the rank. 
Define
\begin{equation}
  \label{eq:1.1}
  \fk a^*= X^*(A)\otimes_{\bZ}\bR,
\end{equation}
where $X^*(A)$ is the lattice of characters of the algebraic torus $A.$ 
Each element $\nu\in\fk a^*$ defines an unramified character
$\chi_\nu$ of $A$, characterized  by the formula
\begin{equation}
  \label{eq:1.2}
  \chi_\nu(\tau(f))=|f|^{\langle \tau,\nu\rangle},\qquad f\in\bF^\times,
\end{equation}
where  $\tau$ is an element of the lattice of one
parameter subgroups $X_*(A).$ Since the torus is split, each element of
$X_*(A)$ can be regarded as a homomorphism of $\bF^\times$ into $A.$
The pairing in the exponent in (\ref{eq:1.2}) corresponds to the
natural identification of $\fk a^*$ with $\Hom[X_*(A),\bR].$
The map $\nu\longrightarrow \chi_\nu$ from $\fk a^*$ to real
unramified characters of $A$ is an isomorphism. We will often
identify the two sets writing simply $\chi\in\fk a^*.$   

\medskip
Let $\vG$ be the (complex) dual group, and let $\vA$ be the torus dual to
$A.$ Then $\fk a^*\otimes_\bR\bC$ is canonically isomorphic to $\va,$
the Lie algebra of $\vA.$ So we can regard $\chi$ as an element of $\va.$ 
We attach to each $\chi$ a nilpotent orbit $\vO(\chi)$ as follows.  By the
Jacobson-Morozov theorem, there is a 1-1 correspondence between
nilpotent orbits $\cCO$ and $\vG$-conjugacy classes of Lie triples 
$\{\ve,\vh,\vf\}$; the correspondence satisfies $\ve\in\cCO$. 
Choose the Lie triple such that
$\vh\in\va.$ Then there are many $\cCO$ such that $\chi$ can be written as 
$w\chi=\vh/2 +\nu$ with  $\nu\in\fk z(\ve,\vh,\vf),$ the centralizer
in $\vg$ of the triple. For example
this is always possible with $\cCO=(0).$ The results in \cite{BM1}
guarantee that for any $\chi$ there is a unique $\cCO(\chi)$ satisfying 
\begin{enumerate}
\item there exists $w\in W$ such that $w\chi=\frac12\vh +\nu$ with
  $\nu\in\fk z(\ve,\vh,\vf),$ 
\item if $\chi$ satisfies property (1) for any other $\cCO',$ then
  $\cCO'\subset\ovl{\cCO}(\chi).$
\end{enumerate}
 
Here is another characterization of the orbit $\cCO(\chi).$
Let
$$
\vg_1:=\{\ x\in \vg\ :\ [\chi,x]=x \ \},\qquad 
\vg_0:=\{ x\in \vg\ :\ [\chi,x]=0 \ \}.
$$ 
Then $\vG_0$, the Lie group corresponding to the Lie algebra $\vg_0$
has an open dense orbit in $\vg_1.$ Its $\vG$ saturation 
in $\vg$ is $\vO(\chi).$

The pair $(\cCO(\chi),\nu)$ has further nice properties. For example
assume that $\nu=0$ in (1) above. Then the representation $L(\chi)$ is
one of the parameters 
that the Arthur conjectures predict to play a role in the residual
spectrum. In particular, $L(\chi)$ should be unitary. In
the $p-$adic case one can verify the unitarity directly as follows.  
In \cite{BM1} it is shown how to calculate the
Iwahori-Matsumoto dual of $L(\chi)$ in the Kazhdan-Lusztig
classification of representations with Iwahori-fixed vector. It turns
out that in the case $\nu=0,$ it is a tempered module, and therefore unitary.
Since the results in \cite{BM1} show that the Iwahori-Matsumoto
involution preserves unitarity, $L(\chi)$ is unitary as well. 
In the real case, a direct proof of the unitarity of $L(\chi)$ (still
with $\nu=0$ as in (1) above) is given in \cite{B3}, and in section
\ref{sec:9} of this paper.  

\medskip
 In the classical Lie algebras, the centralizer $\fk z(\ve,\vh,\vf)$
 is a product of symplectic and orthogonal Lie algebras. We will often
 abbreviate it as $\fk z(\cCO).$ The orbit $\cCO$ is called
 \textit{distinguished} if $\fk z(\cCO)$ does not contain a nontrivial
 torus; equivalently,  the orbit does not meet any  proper Levi
 component.  Let $\vm_{BC}$ be the centralizer of
 a Cartan subalgebra in $\fk z(\cCO).$ This is the Levi component of a
 parabolic subalgebra. The subalgebra
 $\vm_{BC}$ is the Levi subalgebra attached to $\cCO$ by the
 Bala-Carter classification of nilpotent  orbits. The intersection of
 $\vO$ with $\vm_{BC}$ is the other datum  attached to $\cCO,$ a
 distinguished  orbit in $\vm_{BC}.$ We will usually  denote it
 $\vm_{BC}(\vO)$ if we need  to emphasize the dependence on the
 nilpotent orbit.   Let $M_{BC}\subset G$ be the Levi subgroup
 whose Lie algebra  $\fk m_{BC}$ has $\vm_{BC}$ as its dual.

The parameter $\chi$ gives rise to a spherical irreducible
representation $L_{M_{BC}}(\chi)$ on $M_{BC}$ as well as a $L(\chi).$ 
Then $L(\chi)$ is the unique spherical irreducible subquotient of
\begin{equation}
  \label{eq:1.3}
  I_{M_{BC}}(\chi):=Ind_{M_{BC}}^G [ L_{M_{BC}}(\chi)].
\end{equation}
To motivate why we consider $M_{BC}(\cCO),$ we need to recall some
facts about the Kazhdan-Lusztig classification of  representations
with Iwahori fixed vectors in the $p$-adic case. Denote by $\tau$ the 
Iwahori-Matsumoto involution. Then the space of Iwahori fixed vectors of
$\tau(L(\chi))$ is a $W-$representation (see \ref{sec:5.2}), and contains the 
$W-$representation $sgn.$ Irreducible representations with
Iwahori-fixed vectors are parametrized by Kazhdan-Lusztig data; these are $\vG$
conjugacy classes of $(\ve,\chi,\psi)$ where $\ve\in\vg$ is such that
$[\chi,\ve]=\ve,$ and $\psi$ is an irreducible representation of the
component group $A(\chi,\ve).$ To each such parameter there is
associated a standard module $X(\ve,\chi,\psi)$ which contains a
unique irreducible submodule $L(\ve,\chi,\psi).$ All other factors
have parameters $(\ve',\chi',\psi')$ such that 
$$
\vO(\ve)\subset\ovl{\vO(\ve')},\qquad
\vO(\ve)\ne\vO(\ve').
$$
As explained in section 4 and 8 in \cite{BM1}, $X(\ve',\chi',\psi')$
contains $sgn$ if and only if $\psi'=triv.$ Thus if we assume $\vO$ satisfies
(1) and (2) with respect to $\chi$, it follows that  
$X(\ve,\chi,triv)=L(\ve,\chi,triv).$ We would like it to equal
$I_{M_{BC}}$ but this is not true. In general (for an $M$ which
contains $M_{BC}$),   
$L(\ve,\chi,triv)=Ind_M^G[X_M(\ve,\chi,triv)]$ if and only if the component
$A_{M}(\ve,\chi)$ equals the component group $A(\ve,\chi).$ We will 
enlarge $M_{BC}(\cCO)$ to an $M_{KL}$ so that
$A_{M_{KL}}(\ve,\chi)=A(\ve,\chi).$ Note that if
$\vm\subset\vm'$, then $A_M(\ve,\chi)\subset
A_{M'}(\ve,\chi)$. Then  
\begin{equation}
  \label{eq:1.4}
Ind_{M_{KL}}^G[X_{M_{KL}}(\ve,\chi,triv)]=X(\ve,\chi,triv)
=L(\ve,\chi,triv),
\end{equation}
and denoting
\begin{equation}
  \label{eq:1.5imkl}
  I_{M_{KL}}(\chi):=Ind_{M_{KL}}^G[L_{M_{KL}}(\chi)],  
\end{equation}
the equality
\begin{equation}
  \label{eq:1.5}
L(\chi)=I_{M_{KL}}(\chi)
\end{equation}
follows by applying
$\tau$. We remark that $M_{KL}$ depends on $\chi$ as well as $\ve.$
It will be described explicitly in section \ref{sec:2}. A more general
discussion about how canonical $\vm_{KL}$ is, appears in \cite{BC2}.

\medskip
In the real case, we use the same Levi components as in the $p-$adic case. 
Then equality (\ref{eq:1.5}) does not hold for any proper Levi
component. A result essential for the
paper is that equality \textit{does}  hold at the level of 
multiplicities of the \textit{relevant} $K-$types (section
\ref{sec:5.2}, see particularly equation  (\ref{eq:relirred})). 

\medskip
We will use the data $(\cCO,\nu)$ to parametrize the unitary dual. Fix
an $\cCO.$ A  representation $L(\chi)$ will be called a 
\textit{complementary series attached to $\cCO,$} if it is unitary, and
$\vO(\chi)=\cCO.$ To describe it, we need
to give the set of $\nu$ such that $L(\chi)$ with $\chi=\vh/2
+\nu$ is unitary. Viewed as an element of $\fk z(\cCO),$ the element
$\nu$ gives rise to a spherical parameter $( (0),\nu)$ where $(0)$
denotes the trivial nilpotent orbit. The main result in section
\ref{sec:3.2} says that the $\nu$ giving rise to the 
complementary series for $\cCO$ coincide with  the ones giving rise to
the complementary series for $(0)$ on $\fk z(\cCO).$ This is suggestive
of Langlands functoriality. 

It is natural to conjecture that such a result will hold for all split
groups. Recent work of D. Ciubotaru for $F_4,$ and by 
D. Ciubotaru and myself for $E_6, E_7, E_8,$  show that
this is generally true, but there are exceptions.

\bigskip
I give a more detailed outline of the paper. Section \ref{sec:2}
reviews notation from earlier papers. Section \ref{sec:3}  gives a
statement of the main results.  A representation is called
\textit{spherical unipotent} if its parameter is of the form $\vh/2$
for the neutral element of a Lie triple associated to a nilpotent
orbit $\vO.$ 
{
The unitarity of the spherical unipotent representations is dealt with in
section \ref{sec:9}. For the $p-$adic case I simply cite
\cite{BM3}. The unitarity is a consequence of the fact that the
Iwahori-Matsumoto involution preserves unitarity, and spehrical
unipotent representations correspond to tempered representations.  
The real case (sketched in \cite{B2}) is more involved. It is redone in section
\ref{sec:9.5}. The proofs are simpler than the original ones because I
take advantage of the fact that wave front sets, asymptotic supports
and associated varieties ``coincide'' due to \cite{SV}. Section
\ref{sec:10.1} proves an irreducibility result in the real case which
is clear in the $p-$adic case from the work of Kazhdan-Lusztig.  This
is needed for determining the complementary series (definition
\ref{d:2.1}  in section \ref{2.1}).  
}
Sections \ref{sec:4} and \ref{sec:5} deal with the nonunitarity. The
decomposition $\chi=\vh/2 +\nu$ is introduced in section
\ref{sec:3}. It is more common to parametrize the $\chi$ by representatives in
$\va$ which are dominant with respect to some positive root system. We
use Bourbaki's standard realization of the positive system. It
is quite messy to determine the data $(\cCO,\nu)$ from a dominant
parameter, because of the nature of the nilpotent orbits and the 
Weyl group.  Sections \ref{sec:2.3}-\ref{sec:2.4} give a
combinatorial description of $(\cCO,\nu)$ starting from a dominant  $\chi.$ 
Siddhartha Sahi made significant contributions to the phrasing of these
algorithms. 

In the classical cases, the orbit
$\cCO$ is given in terms of partitions. To such a partition we
associate the Levi component  
$$
\check{\fk m_{BC}}:= gl(a_1)\times\dots\times gl(a_k)\times \check{\fk
  g}_0(n_0) 
$$ 
given by the Bala-Carter classification. (The $\vg_0$ in this formula
is \textit{not related} to the one  just after conditions (1) and (2)).
The intersection of $\vO$ with $\vm_{BC}$ is an orbit of the form
$$
(a_1)\times\dots\times (a_r)\times \cCO_0
$$
where $\cCO_0$ is a distinguished nilpotent orbit, and $(a_i)$ is the principal
nilpotent orbit on $gl(a_i).$ This is the distinguished orbit
associated to $\cCO$ by Bala-Carter.  
Then $\chi$ gives rise to irreducible spherical modules 
$L(\chi)$ and $I_{M_{BC}}(\chi)$ as in (\ref{eq:1.3}) 
and $I_{M_{KL}}(\chi)$ as in (\ref{eq:1.5}). The module $L(\chi)$
is the irreducible spherical subquotient of $I_M(\chi)$. As already
mentioned, $I_{M_{KL}}(\chi)=L(\chi)$  in the $p-$adic case, but not
the real case. In all cases,   
the multiplicities of the \textit{relevant $K-$types} in $L(\chi),$
$I_{M_{KL}}(\chi)$ coincide (section \ref{sec:2.6} and \ref{sec:5.2}). 
These are representations of the Weyl group in the $p-$adic case,
representations of the maximal compact subgroup in the real case. 
{Their definition is in chapter \ref{sec:4}; they
are a small finite set of representations which provide necessary
conditions for unitarity which are also sufficient. Their role is
outlined at the beginning of chapter \ref{sec:4}. In particular the issue is
addressed of how the relevant $K-$types allow us to deal with the $p-$adic
case only. A more general class of  $K-$types for split real groups (named
\textit{petite} $K$-types), on which the intertwining operator is
equal to the $p$-adic operator, is defined in \cite{B6}, and the
proofs are more conceptual. In \cite{oda} a more general notion of
\textit{single petaled} is studied independently, and for different
purposes.  The proofs in chapter \ref{sec:4} are the original ones,
and included for completeness. The interested reader can also consult 
\cite{BCP} for more general results of this kind.
}
The determination of the nonunitary parameters proceeds by induction
on the rank of $\vg$ and by the inclusion relations  of the closure
of the orbit $\cCO.$ Section \ref{sec:5} completes the induction step;
it shows that conditions (B) in section \ref{2.1} are necessary. 
The last part of the induction step is actually done in sections
\ref{sec:3.2} and \ref{sec:3.3}.

\bigskip
I would like to thank David Vogan for generously sharing his ideas
about the relation between $K-$types, Weyl group representations and
signatures. They were the catalyst for this paper. The referree as
well as Dan Ciubotaru, Peter Trapa and Siddhartha Sahi have made
numerous suggestions which have improved the exposition and content of
the paper.  

\bigskip
This research was supported by NSF grants DMS-9706758,DMS-0070561 and
DMS-03001712. 

\bigskip
\section{Description of the spherical parameters} \label{sec:2}

\subsection{Explicit Langlands parameters}\label{1.1} We consider 
spherical irreducible representations of the split connected classical 
groups of rank $n$ of type $B, C, D,$ precisely,  $G=So(2n+1),\ G=Sp(2n)$ and
$G=So(2n)$.  These groups will be denoted by $G(n)$ when there is no
danger of confusion ($n$ is the rank). Levi components will be written
as
\begin{equation}
  \label{eq:1.1levi}
  M=GL(k_1)\times\dots\times GL(k_r)\times G_0(n_0),
\end{equation}
where $G_0(n_0)$ is the factor of the same type as $G.$ The
{corresponding complex} Lie
algebras are denoted $\fk g (n)$ and $\fk m=
gl(k_1)\times\dots\times gl(k_r)\times\fk g_0(n_0) .$

\bigskip
As already explained in the introduction, we deal with \textit{real
  unramified characters} only, and this is sufficient for determining 
  the full spherical unitary dual.  

\medskip
In the case of
classical groups, such a character can be represented by a vector of
size the rank of the group. Two such vectors parametrize the same 
irreducible spherical module if they are conjugate via the Weyl group
which acts by permutations and sign changes for type $B,\ C$ and by
permutations and an even number of sign changes in type $D.$

For a given \textit{unramified} 
 $\chi,$ let $L(\chi)$ be the corresponding irreducible spherical module. We
 will occasionally refer to $\chi$ as the infinitesimal character.

\medskip
For any nilpotent orbit $\cCO\subset\vg$ we attach an unramified parameter
$\chi_{\cCO}\in \fk a^*$ as
follows. Recall from the introduction that $\fk a^*\otimes_\bR \bC$ is
canonically isomorphic to $\va$. Let $\{\ve,\vh,\vf\}$ 
be representatives for the Lie triple associated to a nilpotent
orbit $\cCO.$ Then $\chi_\cCO:=\check h/2.$  

Conversely, to each $\chi$ we will attach a nilpotent orbit
$\cCO\subset \check{\fk g}$ and Levi components $M_{BC}$, 
$M_{KL}:=GL(k_1)\times \dots \times
GL(k_r)\times G_0(n_0)$.  In addition we will specify an even
nilpotent orbit $\check \CO_0\subset \check{\fk g}_0(n_0)$ with
unramified character 
$\chi_0:=\chi_{\cCO_0}$ on $\fk g_0(n_0)$, and unramified characters
$\chi_i$ on the  $GL(k_i)$. These data have the property that
$L(\chi)$ is the spherical subquotient of   
\begin{equation}\label{1.1.0}
Ind_{M_{KL}}^G [\bigotimes_i L(\chi_i)\otimes L(\chi_0)].
\end{equation}

\subsection{}\label{1.2a} We introduce the following notation (a
variant of the one used by Zelevinski \cite{ZE}). 

\begin{definition}\label{d:1.1} A {\rm string} is a sequence
\begin{equation}\notag
(a,a+1,\dots ,b-1,b)
\end{equation}
of numbers increasing by 1 from $a$ to $b.$ A set of strings is
called {\rm nested} if for any two strings either the coordinates
{\rm do not} differ by integers, or if they do, then their
coordinates, say   
$(a_1,\dots ,b_1)$ and $(a_2,\dots, b_2),$  satisfy 
\begin{equation}\label{eq:d1}
a_1\le a_2\le b_2\le b_1\qquad \text{ or }\qquad a_2\le a_1\le b_1\le b_2,
\end{equation}
or 
\begin{equation}\label{eq:d2}
b_1+1<a_2 \qquad \text{ or } \qquad b_2+1<a_1.\qed
\end{equation}
A set of strings is called {\rm strongly nested} if the coordinates
of any two  strings either {\rm do not} differ by integers or else
satisfy (\ref{eq:d1}).
\end{definition}
\noindent Each string represents a 1-dimensional spherical representation 
of a $GL(n_i)$ with $n_i=b_i-a_i+1.$ The matchup is
\begin{equation}
  \label{eq:1.1.1a}
  (a,\dots ,b)\longleftrightarrow \big| \det
  \big|^{\frac{a+b}{2}},\quad \text{of}\quad GL(b-a+1),
\end{equation}
{where $\det$ is the determinant character of $GL(n).$}

\medskip
In the case of $G=GL(n),$ we record the following result. For
the $p-$adic case, it originates in the work of Zelevinski, and 
Bernstein-Zelevinski (\cite{ZE} and references therein).  
To each set of strings $(a_1,\dots ,b_1;\dots ;a_k,\dots ,b_k)$ we can
attach a Levi component $M_{BC}:=\prod_{1\le i\le k} GL(n_i),$ and an
induced module  
\begin{equation}
  \label{eq:1.1.1b}
  I(\chi):=Ind_{M_{BC}}^{GL(n)}[\bigotimes L(\chi_i)]
\end{equation}
where $\chi_i$ are obtained from the strings  as in (\ref{eq:1.1.1a}).

In general, if the set of strings is not nested, then the
corresponding induced module is not irreducible.
The coordinates of  $\chi$ in  $\fk a^*\simeq \bR^n$ determine a
set of nested strings as follows. 
{Extract the longest string starting with the smallest
  coordinate in {$\chi$}. Continue, in the same way, to extract sequences from the
remainder until there are no elements left. This set of strings is, up
to the order of the strings, the unique set of nested strings one can
form out of the entries of $\chi.$}

\begin{theorem}
  \label{t:2.1}
Suppose $\bF$ is $p$-adic. Let $(a_1,\dots ,b_1;\dots a_r,\dots ,b_r)$
be a set of nested strings, and $M:=GL(b_1-a_1+1)\times\dots \times
GL(b_r-a_r+1).$ Then
\begin{equation*}
  L(\chi)=Ind_{M}^{GL(n)}\Big[\ \big|\det\big|^{\frac{a_1+b_1}{2}}\cdot\dots
  \cdot \big|\det\big|^{\frac{a_r+b_r}{2}}\ \Big].
\end{equation*}
\end{theorem}

{In the language of section \ref{1.1}, $M_{BC}=M_{KL}=M,$ where $M$
    is the one defined in the theorem.}
The nilpotent orbit $\cCO$ corresponds to the
partition of $n$ with entries $(b_i-a_i+1);$ it is the orbit $\vO(\chi)$
satisfying (1) and (2) in the introduction, with respect to
$\chi=(a_1,\dots ,b_1;\dots a_r,\dots ,b_r)$.     

\medskip
For the real case (still $GL(n)$), the induced module in theorem
\ref{t:2.1}  fails to be irreducible. However 
\begin{equation}
  \label{eq:equalityrelevant}
  [\mu:L(\chi)]=
[\mu: Ind_{M}^{GL(n)}\Big[\ \big|\det\big|^{\frac{a_1+b_1}{2}}\cdot\dots
  \cdot \big|\det\big|^{\frac{a_r+b_r}{2}}\ \Big] 
\end{equation}
whenever $\mu$ is a \textbf{relevant} $K-$type. This is essential for
many of the arguments.  {Relevant
  $K$-types will be defined in chapter \ref{sec:4}, particularly
  section \ref{sec:4.7}.}

 We will generalize this procedure to the other classical
groups. As before, the induced modules that we construct fail to be
irreducible in the real case. {Instead we will establish that
  \begin{equation}
    \label{eq:relirred}
    \dim\Hom_K[\mu :L(\chi)]=\dim\Hom_K[\mu :I_{M_{KL}}(\chi)]
  \end{equation}
for \textit{relevant}   $K-$types, where  $I_{M_{KL}}(\chi)$ is defined
  by  (\ref{eq:1.5imkl}), and make essential use of this fact. } 

\medskip
\subsection{Nilpotent orbits}\label{sec:2.3} In this section we attach a
set of parameters to each nilpotent orbit $\cCO\subset \vg.$ Let
$\{\ch e,\ch h,\ch f\}$ be a Lie triple so that $\ch e\in\cCO,$ and
let $\fk z(\cCO)$ be its centralizer. In order for $\chi$ to be a
parameter attached to $\cCO$ we require that
\begin{equation}
  \label{eq:2.3.1}
  \chi=\vh/2 +\nu,\qquad \nu\in\fk z(\cCO),\ \text{ semisimple,}
\end{equation}
but also that if
\begin{equation}
  \label{eq:2.3.2}
  \chi = \ch h'/2 +\nu',\qquad \nu'\in\fk z(\cCO'),\ \text{ semisimple}
\end{equation}
for another nilpotent orbit $\cCO'\subset\vg,$ 
then $\cCO'\subset\ovl{\cCO}.$ In \cite{BM1}, it is shown that the
orbit of $\chi,$ uniquely determines  $\cCO$ and the conjugacy class
of $\nu\in\fk z(\cCO).$ We describe the pairs $(\cCO,\nu)$ explicitly
in the classical cases.

Nilpotent orbits are parametrized by partitions 
\begin{equation}
  \label{eq:2.3.3}
(\unb{r_1}{1,\dots ,1},\unb{r_2}{2,\dots ,2},\dots ,
\unb{r_j}{j,\dots   ,j},\dots ),
\end{equation}
satisfying the following constraints.
\begin{description}
\item[$\vg$ type $A_{n-1}$] $\ gl(n)$, partitions of $n.$
 
\item[$\vg$ type $B_n$] $\ so(2n+1)$,
partitions of $2n+1$ such that every even part occurs an even number of times.

\item[$\vg$ type $C_n$] $\ sp(2n)$,
partitions of $2n$ such that every odd part occurs an even number of
times.

\item[$\vg$ type $D_n$] $\ so(2n)$,
partitions of $2n$ such that every even part occurs an even number of
times. In the case when every part of the partition is  even, there
are two conjugacy classes of nilpotent orbits with the same Jordan
blocks, labelled (I) and (II). The two orbits are conjugate under the
action of  $O(2n).$ 

\end{description}

The Bala-Carter classification is particularly well suited for
describing the parameter spaces attached to the $\cCO\subset \vg.$  
An orbit is called \textit{distinguished} if it does not meet
any proper Levi component. In type A, the only distinguished orbit is
the principal nilpotent orbit, where the partition has only one
part. In the other cases, the distinguished orbits are the ones where
each part of the partition occurs at most once. In particular, these
are \textit {even nilpotent orbits}, \ie  $\ad\vh$ has even
eigenvalues only on $\vg$. Let $\cCO\subset \vg$ be an arbitrary nilpotent
orbit. We need to put it into as small as possible Levi component
$\vm.$ In type A, if the partition is $(a_1,\dots, a_k),$ the Levi
component is $\vm_{BC}=gl(a_1)\times\dots\times gl(a_k).$ In the other
classical types, the orbit $\cCO$ meets a proper Levi component 
if and only if one of
the $r_j>1.$ So separate as many pairs $(a,a)$ from the partition as
possible, and rewrite it as
\begin{equation}
  \label{eq:2.3.4}
  ((a_1,a_1),\dots ,(a_k,a_k);d_1,\dots ,d_l),
\end{equation}
with $d_i<d_{i+1}.$ The Levi component $\vm_{BC}$ attached to this
nilpotent by Bala-Carter is
\begin{equation}
  \label{eq:2.3.5}
\vm_{BC}= gl(a_1)\times\dots\times gl(a_k)\times \vg_0(n_0)\quad
n_0:=n-\sum a_i,
\end{equation}
The distinguished nilpotent orbit is the one with partition $(d_i)$ on
$\vg(n_0),$ principal nilpotent on each $gl(a_j).$ The $\chi$ of
the form $\vh/2 +\nu$ are the ones with $\nu$ an element of  the center
of $\vm_{BC}.$ The explicit form is 
\begin{equation}
  \label{eq:2.3.6}
  (\dots
  ;-\frac{a_i-1}{2}+\nu_i,\dots,\frac{a_i-1}{2}+\nu_i,\dots;\vh_0/2 ),
\end{equation}
where $\vh_0$ is the middle element of a triple corresponding to
$(d_i)$. We will write out $(d_i)$ and $\vh_0/2$ in sections
\ref{sec:2.3a}-\ref{sec:2.3d}. 

\medskip
We will consider more general cases where we write the partition of
$\cCO$ in the form (\ref{eq:2.3.4}) so that the $d_i$ are not
necessarily distinct, but $(d_i)$ forms an even nilpotent orbit in
$\vg_0(n_0).$ {This will be the situation for $\vm_{KL}.$}  

\medskip
The parameter $\chi$ determines an irreducible spherical module
$L(\chi)$ for $G$ as well as an $L_M(\chi)$ for $M=M_{BC}$ or $M_{KL}$ 
of the form
\begin{equation}
  \label{eq:2.3.7}
L_1(\chi_1)\otimes\dots \otimes L_k(\chi_k)\otimes  L_0(\chi_0),
\end{equation}
where the $L_i(\chi_i)$ {$i=1,\dots ,k$} are one dimensional. 
We will consider the relation between the induced module
\begin{equation}
  \label{eq:2.3.8}
  I_M(\chi):=Ind_M^G[L_M(\chi)],
\end{equation}
and $L(\chi)$.

\subsection{G of Type A}\label{sec:2.3a}\quad
We write the $\vh/2$ for a nilpotent $\cCO$ corresponding to
$(a_1,\dots ,a_k)$ with $a_i\le a_{i+1}$ as 
\begin{equation*}
  (\dots ;-\frac{a_i-1}{2},\dots ,\frac{a_i-1}{2};\dots ).
\end{equation*}
The parameters of the form $\chi=\vh/2+\nu$ are  then
\begin{equation}
  \label{2.3a.1}
  (\dots ;-\frac{a_i-1}{2}+\nu_i,\dots ,\frac{a_i-1}{2}+\nu_i;\dots ).
\end{equation}
Conversely, given a parameter as a concatenation of strings
\begin{equation}
  \label{eq:2.3a.2}
\chi=(\dots ;A_i,\dots ,B_i;\dots ),
\end{equation}
it is of the form $\vh/2+\nu$ where $\vh$ is the neutral element for
the nilpotent orbit with partition $(A_i+B_i+1)$ (the parts need not
be in any particular order) and
$\nu_i=\frac{A_i-B_i}{2}$. We recall the following well known result about
closures of nilpotent orbits.
\begin{lemma}
  \label{l:2.3a}
Assume $\cCO$ and $\cCO'$ correspond to the {(increasing)} partitions
$(a_1,\dots , a_k)$ and $(b_1,\dots ,b_k)$ respectively, where some of
the $a_i$ or $b_j$ may be zero in order to have the same number $k.$  
The following are equivalent
\begin{enumerate}
\item $\cCO'\subset\ovl{\cCO}.$
\item $\sum_{i\ge s} a_i\ge \sum_{i\ge s} b_i$ for all $k\ge s\ge 1.$ 
\end{enumerate}
\end{lemma}

\begin{proposition}
  \label{p:2.3a}
A parameter $\chi$ as in (\ref{2.3a.1}) is attached to $\cCO$ in the
sense of satisfying (\ref{eq:2.3.1}) and (\ref{eq:2.3.2}) 
 if and only if it is nested.
\end{proposition}
\begin{proof}
Assume the strings are not nested. There must be two strings
\begin{equation}
  \label{eq:2.3a.3}
  (A,\dots ,B),\qquad (C,\dots , D)
\end{equation}
such that $A-C\in\bZ$, and $A< C\le B< D,$ or $C=B+1$. Then by
conjugating $\chi$ by the Weyl group to a $\chi'$, we can rearrange
the coordinates of the two strings in (\ref{eq:2.3a.3}) so that the strings 
{\begin{equation}
  \label{eq:2.3a.4}
  (A,\dots ,D),\quad (C,\dots B),\qquad\text{ or }\quad (A,\dots
  ,B,B+1=C,\dots ,D).
\end{equation}
}
appear instead. Then by the lemma, $\chi'=\vh'/2+\nu'$ for a strictly larger
nilpotent $\cCO'$. 

Conversely, assume $\chi=\vh/2+\nu$, so it is written as strings, and
they are nested. The nilpotent orbit for which the neutral element is
$\vh/2$ has partition given by the lengths of the strings, say
$(a_1,\dots a_k)$ in increasing order. 
If $\chi$ is nested, then $a_k$ is the length of
the longest string of entries we can extract from the coordinates of
$\chi,$ $a_{k-1}$ the longest string we can extract from the remaining
coordinates and so on. Then (2) of lemma {\ref{l:2.3a}} precludes the
possibility that some conjugate $\chi'$ equals $\vh'/2+\nu'$ for a strictly
larger nilpotent orbit.    
\end{proof}
In type A, $\vm_{KL}=\vm_{BC}.$
\subsection{G of Type B}\label{sec:2.3b}
Rearrange the parts of the partition of $\cCO\subset
sp(2n,\bC),$ in the form (\ref{eq:2.3.4}), 
\begin{equation}
  \label{eq:2.3b.1}
((a_1,a_1),\dots
,(a_k,a_k);2x_0,\dots, 2x_{2m})  
\end{equation}
The $d_i$ have been relabeled as $2x_i$ and a $2x_0=0$ is added if
necessary, to
insure that there is an odd number. The $x_i$ are integers, because
all the odd parts of the partition of $\cCO$ occur an even number of
times, and were therefore extracted as $(a_i,a_i).$ The $\chi$ of the
form $\vh/2+\nu$ are 
\begin{equation}
  \label{eq:2.3b.2}
  \begin{aligned}
(\dots;-\frac{a_i-1}{2}+\nu_i,&\dots ,
\frac{a_i-1}{2}+\nu_i;\dots;\\
&\unb{n_{1/2}}{1/2,\dots ,1/2},\dots
,\unb{n_{x_{2m}-1/2}}{x_{2m}-1/2,\dots ,x_{2m}-1/2} ).  
  \end{aligned}
\end{equation}
where
\begin{equation}
  \label{eq:2.3b.3}
  n_{l-1/2}=\#\{x_i\ge l\}.
\end{equation}
Lemma \ref{l:2.3a} holds for this type verbatim. So the following
proposition holds.
\begin{proposition}
  \label{p:2.3b}
A parameter $\chi=\vh/2+\nu$ cannot be conjugated to one of the form
$\vh'/2 +\nu'$ for any larger nilpotent $\cCO'$ if and only if 
\begin{enumerate}
\item the set of strings satisfying
  $\frac{a_i-1}{2}+\nu_i-\frac{a_j-1}{2}-\nu_j\in\bZ$ are nested.
\item the strings satisfying $\frac{a_i-1}{2}+\nu_i\in 1/2 + \bZ$ satisfy the
  additional condition that either  $x_{2m}+1/2<
  -\frac{a_i-1}{2}+\nu_i$ or there is $j$ such that
{\begin{equation}
    \label{eq:2.3b.4}
    x_j+1/2<-\frac{a_i-1}{2}+\nu_i\le \frac{a_i-1}{2}+\nu_i\le x_{j+1}-1/2.
  \end{equation}
}
\end{enumerate}
\end{proposition}
The Levi component $\vm_{KL}$ is obtained from $\vm_{BC}$ as
follows. Consider the strings for which $a_i$ is even, and $\nu_i=0$. If
$a_i$ is not equal to any $2x_j,$ then remove one pair $(a_i,a_i)$, 
and add two $2x_j=a_i$ to the last part of (\ref{eq:2.3b.1}). 
For example, if the nilpotent orbit $\vO$ is
\begin{equation}
  \label{eq:2.3b.5}
  (2,2,2,3,3,4,4),
\end{equation}
then the parameters of the form $\vh/2+\nu$ are
\begin{equation}
\label{eq:2.3b.6}
  \begin{aligned}  
(&-1/2+\nu_1,1/2+\nu_1;-1+\nu_2,\nu_2,1+\nu_2;\\
&-3/2+\nu_3,-1/2+\nu_3,1/2+\nu_3,3/2+\nu_3;1/2)
\end{aligned}    
\end{equation}
The Levi component is $\vm_{BC}=gl(2)\times gl(3)\times gl(4)\times \vg(1)$. If
$\nu_3\ne 0,$ then $\vm_{BC}=\vm_{KL}.$ But if $\nu_3=0,$ then
$\vm_{KL}=gl(2)\times gl(3)\times \vg(5).$ The parameter is rewritten

  \begin{align}  \label{eq:2.3b.7}
 &\cCO\longleftrightarrow ((2,2)(3,3);2,4,4) \\
 &\chi\longleftrightarrow
 (-1/2+\nu_1,1/2+\nu_1;-1+\nu_2,\nu_2,1+\nu_2;1/2,1/2,1/2,3/2,3/2).\notag
  \end{align}

The explanation is as follows. For a partition (\ref{eq:2.3.3}),
\begin{equation}
  \label{eq:2.3b.8}
 \fk z(\cCO)=sp(r_1)\times so(r_2)\times sp(r_3)\times\dots
\end{equation} 
and the centralizer in $\vG$ is a product of $Sp(r_{2j+1})$ and
$O(r_{2j}),$ \ie $Sp$ for the odd parts, $O$ for the even parts. 
Thus  the component group {$A(\check h,\check e)$}, which by
\cite{BV2} also equals $A(\ve)$,  is a product of $\bZ_2$, one for each 
$r_{2j}\ne 0.$ Then $A(\chi,\ve)=A(\nu,\vh,\ve).$ In general
$A_{M_{BC}}(\chi,\ve)=A_{M_{BC}}(\ve)$  embeds canonically into
$A(\chi,\ve)$, but the two are  not necesarily equal. In this case
they are unless one of the $\nu_i=0$ for an even $a_i$ with the additional
property that there is no $2x_j=a_i.$  

We can rewrite each of the remaining  strings 
\begin{equation}
  \label{eq:2.3b.9}
  (-\frac{a_i-1}{2}+\nu_i,\dots ,\frac{a_i-1}{2}+\nu_i)
\end{equation}
 as 
\begin{align}
\chi_i:=&(f_i+\tau_i,f_i+1+\tau_i,\dots,F_i+\tau_i),  \label{eq:2.3b.10}\\
&\text{ satisfying}\notag\\  
&f_i\in\bZ+1/2,\quad 
0\le\tau_i\le 1/2,\quad F_i=f_i+a_i,\label{eq:2.3b.11}.   \\
 &{ |f_i+\tau_i|\ge |F_i+\tau_i| \text{ if }
\tau_i=1/2} \notag
\end{align}
This is done as follows. We can immediately get an expression like
(\ref{eq:2.3b.10}) with $0\le \tau_i<1,$ by defining $f_i$ to be the
largest element in $\bZ+1/2$ less than or equal to
$-\frac{a_i-1}{2}+\nu_i$. If $\tau_i\le 1/2$ we are done. Otherwise,
use the Weyl group to change the signs of all entries of the string,
and put them in increasing order. This replaces $f_i$ by $-F_i-1,$ and
$\tau_i$ by $1-\tau_i.$  The presentation of the strings subject to
(\ref{eq:2.3b.11}) is unique
except when $\tau_i=1/2.$ In this case the argument just given
provides the presentation $(f_i+1/2,\dots ,F_i+1/2),$ but also
provides the presentation 
\begin{equation}
(-F_i-1+1/2,\dots ,-f_i-1+1/2).
\label{eq:2.3b.12}
\end{equation}
We choose between (\ref{eq:2.3b.10}) and (\ref{eq:2.3b.12}) the one
whose leftmost term is larger in absolute value. That is, we require
$|f_i+\tau_i|\ge |F_i+\tau_i|$ whenever $\tau_i=1/2.$ 

\subsection{G of Type C}\label{sec:2.3c}\quad 
Rearrange the parts of the partition of \newline $\cCO\subset
so(2n+1,\bC),$ in the form (\ref{eq:2.3.4}), 
\begin{equation}
  \label{eq:2.3c.1}
((a_1,a_1),\dots ,(a_k,a_k);2x_0+1,\dots, 2x_{2m}+1);  
\end{equation}
The $d_i$ have been relabeled as $2x_i+1$. In this case it is automatic
that there is an odd number of nonzero $x_i.$
The $x_i$ are integers, because
all the even parts of the partition of $\cCO$ occur an even number of
times, and were threrefore extracted as $(a_i,a_i).$ The $\chi$ of the
form $\vh/2+\nu$ are 
\begin{equation}
  \label{eq:2.3c.2}
(\dots ;-\frac{a_i-1}{2}+\nu_i,\dots , \frac{a_i-1}{2}+\nu_i;\dots 
;\unb{n_{0}}{0,\dots ,0},\dots ,\unb{n_{x_{2m}}}{x_{2m},\dots ,x_{2m}}). 
\end{equation}
where
\begin{equation}
  \label{eq:2.3c.3}
  n_{l}=
  \begin{cases}
m\quad &\text{ if } l=0,\\
\#\{ x_i\ge l\}\quad &\text{ if } l\ne 0 .    
  \end{cases}
\end{equation}
Lemma \ref{l:2.3a} holds for this type verbatim. So the following
proposition holds.
\begin{proposition}
  \label{p:2.3c}
A parameter $\chi=\vh/2+\nu$ cannot be conjugated to one of the form
$\vh'/2 +\nu'$ for any larger nilpotent $\cCO'$ if and only if 
\begin{enumerate}
\item the set of strings satisfying
  $\frac{a_i-1}{2}+\nu_i-\frac{a_j-1}{2}-\nu_j\in \bZ$ are nested.
\item the strings satisfying $\frac{a_i-1}{2}+\nu_i\in\bZ$ satisfy the
  additional condition that either  $x_{2m}+1<
  -\frac{a_i-1}{2}+\nu_i$ or there is $j$ such that
{\begin{equation}
    \label{eq:2.3c.4}
    x_j+1<-\frac{a_i-1}{2}+\nu_i\le \frac{a_i-1}{2}+\nu_i\le x_{j+1}.
  \end{equation}
}
\end{enumerate}
\end{proposition}
The Levi component $\vm_{KL}$ is obtained from $\vm_{BC}$ as
follows. Consider the strings for which $a_i$ is odd and $\nu_i=0$. If
$a_i$ is not equal to any $2x_j+1,$ then remove one pair $(a_i,a_i)$, 
and add two $2x_j+1=a_i$ to the last part of (\ref{eq:2.3c.1}). 
For example, if the nilpotent orbit is
\begin{equation}
  \label{eq:2.3c.5}
  (1,1,1,3,3,4,4)=((1,1),(3,3),(4,4);1),
\end{equation}
then the parameters of the form $\vh/2+\nu$ are
\begin{equation}
\label{eq:2.3c.6}
  \begin{aligned}  
(&\nu_1;-1+\nu_2,\nu_2,1+\nu_2;\\
&-3/2+\nu_3,-1/2+\nu_3,1/2+\nu_3,3/2+\nu_3)
\end{aligned}    
\end{equation}
The Levi component is $\vm_{BC}=gl(1)\times gl(3)\times gl(4)$. If
$\nu_2\ne 0,$ then $\vm_{BC}=\vm_{KL}.$ But if $\nu_2=0,$ then
$\vm_{KL}=gl(1)\times gl(4)\times \vg(3).$ The parameter is rewritten
  \begin{align}  \label{eq:2.3c.7}
 &\cCO\longleftrightarrow ((1,1),(4,4);1,3,3) \\
 &\chi\longleftrightarrow
 (\nu_1;-3/2+\nu_3,-1/2+\nu_3,1/2+\nu_3,3/2+\nu_3;0,1,1).\notag
  \end{align}
The Levi component $\vm_{KL}$ is unchanged if $\nu_1=0.$

The explanation is as follows. For a partition (\ref{eq:2.3.3}),
\begin{equation}
  \label{eq:2.3c.8}
 \fk z(\cCO)=so(r_1)\times sp(r_2)\times so(r_3)\times\dots  
\end{equation} 
and the centralizer in $\vG$ is a product of $O(r_{2j+1})$ and
$Sp(r_{2j}),$ \ie $O$ for the odd parts, $Sp$ for the even parts. 
Thus  the component group is a product of $\bZ_2$, one for each
$r_{2j+1}\ne 0.$ Then $A(\chi,\ve)=A(\nu,\vh,\ve),$ and 
so $A_{M_{BC}}(\chi,\ve)=A(\chi,\ve)$
unless one of the $\nu_i=0$ for an odd $a_i$ with the additional
property that there is no $2x_j+1=a_i.$  

We can rewrite each of the remaining  strings 
\begin{equation}
  \label{eq:2.3c.9}
  (-\frac{a_i-1}{2}+\nu_i,\dots ,\frac{a_i-1}{2}+\nu_i)
\end{equation}
 as 
\begin{align}
\chi_i:=&(f_i+\tau_i,f_i+1+\tau_i,\dots,F_i+\tau_i),  \label{eq:2.3c.10}\\
&\text{ satisfying}\notag \\  
&f_i\in\bZ,\quad 
0\le\tau_i\le 1/2,\quad F_i=f_i+a_i \label{eq:2.3c.11}\\
&{|f_i+\tau_i|\ge |F_i+\tau_i| \text{ if } \tau_i=1/2}\notag.  
\end{align}
This is done as follows. We can immediately get an expression like
(\ref{eq:2.3c.10}) with $0\le \tau_i<1,$ by defining $f_i$ to be the
largest element in $\bZ$ less than or equal to
$-\frac{a_i-1}{2}+\nu_i$. If $\tau_i\le 1/2$ we are done. Otherwise,
use the Weyl group to change the signs of all entries of the string,
and put them in increasing order. This replaces $f_i$ by $-F_i-1,$ and
$\tau_i$ by $1-\tau_i.$  The presentation of the strings subject to
(\ref{eq:2.3c.11}) is unique
except when $\tau_i=1/2.$ In this case the argument just given also
provides the presentation 
\begin{equation}
(-F_i-1+1/2,\dots ,-f_i-1+1/2).
\label{eq:2.3c.12}
\end{equation}
We choose between (\ref{eq:2.3c.10}) and (\ref{eq:2.3c.12}) the one
whose leftmost term is larger in absolute value. That is, we require
$|f_i+\tau_i|\ge |F_i+\tau_i|$ whenever $\tau_i=1/2.$

\subsection{G of Type D}\label{sec:2.3d}\quad
Rearrange the parts of the partition of $\cCO\subset
so(2n,\bC),$ in the form (\ref{eq:2.3.4}), 
\begin{equation}
  \label{eq:2.3d.1}
( (a_1,a_1),\dots ,(a_k,a_k);2x_0+1,\dots, 2x_{2m-1}+1)  
\end{equation}
The $d_i$ have been relabeled as $2x_i+1$. In this case it is
automatic that there is an even number of nonzero $2x_i+1.$
The $x_i$ are integers, because
all the even parts of the partition of $\cCO$ occur an even number of
times, and were therefore extracted as $(a_i,a_i).$ The $\chi$ of the
form $\vh/2+\nu$ are 
\begin{equation}
  \label{eq:2.3d.2}
(\dots;-\frac{a_i-1}{2}+\nu_i,\dots , \frac{a_i-1}{2}+\nu_i;\dots;
\unb{n_{0}}{0,\dots ,0},\dots ,\unb{ n_{x_{2m-1} }}{x_{2m-1},\dots
  ,x_{2m-1}} ).  
\end{equation}
where
\begin{equation}
  \label{eq:2.3d.3}
  n_{l}=
  \begin{cases}
m\quad &\text{ if } l=0,\\
\#\{ x_i\ge l\}\quad &\text{ if } l\ne 0 .    
  \end{cases}
\end{equation}
Lemma \ref{l:2.3a} holds for this type verbatim. So the following
proposition holds.
\begin{proposition}
  \label{p:2.3d}
A parameter $\chi=\vh/2+\nu$ cannot be conjugated to one of the form
$\vh'/2 +\nu'$ for any larger nilpotent $\cCO'$ if and only if 
\begin{enumerate}
\item the set of strings satisfying
  $\frac{a_i-1}{2}+\nu_i-\frac{a_j-1}{2}-\nu_j\in \bZ$ are nested.
\item the strings satisfying $\frac{a_i-1}{2}+\nu_i\in\bZ$ satisfy the
  additional condition that either  $x_{2m-1}+1<
  -\frac{a_i-1}{2}+\nu_i$ or there is $j$ such that
{\begin{equation}
    \label{eq:2.3d.4}
    x_j+1<-\frac{a_i-1}{2}+\nu_i\le \frac{a_i-1}{2}+\nu_i\le x_{j+1}  .
  \end{equation}
}
\end{enumerate}
\end{proposition}
The Levi component $\vm_{KL}$ is obtained from $\vm_{BC}$ as
follows. Consider the strings for which $a_i$ is odd and $\nu_i=0$. If
$a_i$ is not equal to any $2x_j+1,$ then remove one pair $(a_i,a_i)$, 
and add two $2x_j+1=a_i$ to the last part of (\ref{eq:2.3d.1}). 
For example, if the nilpotent orbit is
\begin{equation}
  \label{eq:2.3d.5}
  (1,1,3,3,4,4),
\end{equation}
then the parameters of the form $\vh/2+\nu$ are
\begin{equation}
\label{eq:2.3d.6}
  \begin{aligned}  
(&\nu_1;-1+\nu_2,\nu_2,1+\nu_2;\\
&-3/2+\nu_3,-1/2+\nu_3,1/2+\nu_3,3/2+\nu_3)
\end{aligned}    
\end{equation}
The Levi component is $\vm_{BC}=gl(1)\times gl(3)\times gl(4)$. If
$\nu_2\ne 0$ and  $\nu_1\ne 0$, then $\vm_{BC}=\vm_{KL}.$ If $\nu_2=0$ and
$\nu_1\ne 0,$ then
$\vm_{KL}=\vg(3)\times gl(1)\times gl(4).$ If $\nu_2\ne 0$ and
$\nu_1=0,$ then $\vm_{KL}=gl(3)\times gl(4)\times \vg(1).$ If
$\nu_1=\nu_2=0,$ then $\vm_{KL}=gl(4)\times \vg(4).$ The parameter is rewritten
  \begin{align}  \label{eq:2.3d.7}
 &\cCO\longleftrightarrow ((1,1),(4,4);3,3) \\
 &\chi\longleftrightarrow
 (\nu_1;-3/2+\nu_3,-1/2+\nu_3;1/2+\nu_3,3/2+\nu_3;0,1,1).\notag
  \end{align}

The explanation is as follows. For a partition (\ref{eq:2.3.3}),
\begin{equation}
  \label{eq:2.3d.8}
 \fk z(\cCO)=so(r_1)\times sp(r_2)\times so(r_3)\times\dots  
\end{equation} 
and the centralizer in $\vG$ is a product of $O(r_{2j+1})$ and
$Sp(r_{2j}),$ \ie $O$ for the odd parts, $Sp$ for the even parts. 
Thus  the component group is a product of $\bZ_2$, one for each
$r_{2j+1}\ne 0.$ Then $A(\chi,\ve)=A(\nu,\vh,\ve),$ and 
so $A_{M_{BC}}(\chi,\ve)=A(\chi,\ve)$
unless one of the $\nu_i=0$ for an odd $a_i$ with the additional
property that there is no $2x_j+1=a_i.$  

We can rewrite each of the remaining  strings 
\begin{equation}
  \label{eq:2.3d.9}
  (-\frac{a_i-1}{2}+\nu_i,\dots ,\frac{a_i-1}{2}+\nu_i)
\end{equation}
 as 
\begin{align}
\chi_i:=&(f_i+\tau_i,f_i+1+\tau_i,\dots,F_i+\tau_i),  \label{eq:2.3d.10}\\
&\text{ satisfying}\qquad  f_i\in\bZ,\quad 
0\le\tau_i\le 1/2,\quad F_i=f_i+a_i \label{eq:2.3d.11}\\
&{|f_i+\tau_i|\ge |F_i+\tau_i| \text{ if } \tau_i=1/2}\notag.  
\end{align}
This is done as in types B and C, but see the remarks which have to do
with the fact that $-Id$ is not in the Weyl group. 
We can immediately get an expression like
(\ref{eq:2.3d.10}) with $0\le \tau_i<1,$ by defining $f_i$ to be the
largest element in $\bZ$ less than or equal to
$-\frac{a_i-1}{2}+\nu_i$. If $\tau_i\le 1/2$ we are done. Otherwise,
use the Weyl group to change the signs of all entries of the string,
and put them in increasing order. This replaces $f_i$ by $-F_i-1,$ and
$\tau_i$ by $1-\tau_i.$  The presentation of the strings subject to
(\ref{eq:2.3d.11}) is unique
except when $\tau_i=1/2.$ In this case the argument just given also
provides the presentation 
\begin{equation}
(-F_i-1+1/2,\dots ,-f_i-1+1/2).
\label{eq:2.3d.12}
\end{equation}
We choose between (\ref{eq:2.3d.10}) and (\ref{eq:2.3d.12}) the one
whose leftmost term is larger in absolute value. That is, we require
$|f_i+\tau_i|\ge|F_i+\tau_i|$ whenever $\tau_i=1/2.$

\noindent\textbf{Remarks}
\begin{enumerate}
\item A (real) spherical parameter $\chi$ is hermitian if
and only if there is $w\in W(D_n)$ such that $w\chi=-\chi.$ This is the
case if the parameter has a coordinate equal to zero, or if none of the
coordinates are 0, but then $n$ must be even. 
\item Assume the nilpotent orbit $\cCO$ is very even, \ie all the
  parts of the partition are even (and therefore occur an even number
  of times). The  nilpotent orbits labelled (I) and (II) are characterized by
  the fact that $\vm_{BC}$ is of the form
\begin{equation*}
  \begin{aligned}
(I)&\longleftrightarrow gl(a_1)\times\dots\times gl(a_{k-1})\times gl(a_k),\\
(II)&\longleftrightarrow gl(a_1)\times\dots\times gl(a_{k-1})\times gl(a_k)'.
  \end{aligned}
\end{equation*}
The last $gl$ factors differ by which extremal root
of the fork at the end of the diagram for $D_n$ is in the Levi
component. The string for $k$ is 
\begin{equation*}
  \begin{aligned}
(I)&\longleftrightarrow 
(-\frac{a_k-1}{2}+\nu_k,\dots ,\frac{a_k-1}{2}+\nu_k),\\
(II)&\longleftrightarrow 
(-\frac{a_k-1}{2}+\nu_k,\dots \frac{a_k-3}{2}+\nu_k,-\frac{a_k-1}{2}-\nu_k).
  \end{aligned}
\end{equation*}
We can put the parameter in the form (\ref{eq:2.3d.10}) and
(\ref{eq:2.3d.11}), because all strings are even length. In any case
(I) and (II) are conjugate by the outer automorphism, and for
unitarity it is enough to consider the case of (I).

The assignment of a nilpotent orbit (I) or (II) to a parameter is
unambiguous. If a $\chi$ has a coordinate equal to 0, it might be written as
$h_I/2+\nu_I$ or $h_{II}/2+\nu_{II}.$ But then it can also be
written as $h'/2+\nu'$ for a larger nilpotent orbit. For example, in
type $D_2,$ the two cases are $(2,2)_{I}$ and $(2,2)_{II}$, and  we
can write 
\begin{equation*}
  \begin{aligned}
    (I)\longleftrightarrow (1/2,-1/2)+(\nu,\nu),\\
    (II)\longleftrightarrow (1/2,1/2)+(\nu,-\nu).
  \end{aligned}
\end{equation*}
The two forms are not conjugate unless the parameter contains a $0$. But then
 it has to be $(1,0)$ and this
corresponds to $(1,3),$ the larger principal nilpotent orbit.

\item Because we can only change an even number of signs using the 
Weyl group, we might not be able to change all the signs
of a string. We can always do this if the parameter contains a
coordinate equal to 0, or if the length of the string is even. 
If there is an odd length string, and none of the
coordinates of $\chi$ are 0, changing all of the
signs of the string cannot be achieved unless some other 
coordinate changes sign as well. However if $\chi=\vh/2+\nu$ cannot be made to
satisfy (\ref{eq:2.3d.10}) and (\ref{eq:2.3d.11}), then $\chi',$ the
parameter obtained from $\chi$ by applying the outer automorphism, 
can. Since $L(\chi)$ and $L(\chi')$ are either both unitary of both
nonunitary, it is enough to consider just the cases that can be made
to satisfy  (\ref{eq:2.3d.10}) and (\ref{eq:2.3d.11}).
{For example, the parameters 
\begin{equation*}
  \begin{aligned}
&(-1/3,2/3,5/3;-7/4,-3/4,1/4),\\
&(-5/3,-2/3,1/3;-7/4,-3/4,1/4)   
  \end{aligned}
\end{equation*}
in type $D_6$ are of this kind. Both parameters are in a form
satisfying (\ref{eq:2.3d.10}) but only the second one satisfies
(\ref{eq:2.3d.11}). The first one cannot be conjugated by $W(D_6)$ 
to one satisfying (\ref{eq:2.3d.11}).}

\end{enumerate}

\subsection{Relation between infinitesimal characters and strings}
\label{sec:2.4} 
In the previous sections we described for each nilpotent orbit $\cCO$
the parameters of the form $\vh/2+\nu$ with $\nu\in\fk z(\cCO)$
semisimple, along with  condition (\ref{eq:2.3.2}). In this section we
give algorithms to find the data $(\cCO,\nu)$ satisfying (\ref{eq:2.3.1}) and
(\ref{eq:2.3.2}), and the various Levi components  from a
$\chi\in\va.$ The formulation was  suggested by S. Sahi. Given a
$\chi\in\va,$ we need to specify,
\begin{description}
\item[(a)] strings,  same as sequences of coordinates with increment 1,
\item[(b)] a partition, same as a nilpotent orbit $\vO\subset\vg,$
\item[(c)] the centralizer of a Lie triple corresponding to $\fz(\vO)$,
\item[(d)] coordinates of the parameter $\nu,$ coming from the
  decompositon $\chi=\vh/2 +\nu$, 
\end{description}
Furthermore, we give algorithms for 

\begin{description}
\item[(e)] two Levi components $\vm_{BC}$ and $\vm_{KL}$,
\item[(f)] another two Levi components $\vm_e$ and $\vm_o,$
\item[(g)] one dimensional characters $\chi_e$ and $\chi_o$ of the Levi
  components $\vm_e$ and $\vm_o.$
\end{description}
Parts (f) and (g) are described in detail in section
\ref{sec:5.3}. These Levi components are used to compute
multiplicities of relevant $K-$types in $L(\chi).$

\subsubsection*{Algorithms for (a) and (b)}\

\subsubsection*{\textbf{Step 0}} \ 

\noindent{\textit{G of type C}.} Double the number of $0$'s and add one more.

\noindent{\textit{G of type D}.} Double the number of
$0$'s. If there are no coordinates equal to 0 and the rank is odd, the
parameter is not hermitian. If the rank is even, only an even number
of sign changes are allowed in the subsequent steps.

\subsubsection*{\textbf{Step 1}} \ 

\noindent{\textit{G of type C,D}.} Extract maximal strings of the form
$(0,1,\dots )$. Each contributes a part in the partition of size
\textit{2(length of string)-1} to $\vO$.

\noindent{\textit{G of type B}.} Extract maximal strings of type
$(1/2,3/2,\dots ).$ Each contributes a part of size
\textit{2(length)} to $\vO.$ 

\subsubsection*{\textbf{Step 2}}\ 

 For all types, extract maximal strings {from the
   remaining entries after Step 1,} 
changing signs if necessary. Each string contributes two parts of size
\textit{(length of string)} to $\vO.$ In type D, if the rank is odd and no
coordinate of the original $\chi$ is 0, the parameter is not hermitian. If there
are no 0's and the rank of type D is even, only an even number of sign
changes is allowed. In this case, the last string might be $(\dots ,b,-b-1).$
If so, and all strings are of even size, $\vO$ is very even, and is
labelled $II$. If all strings are of the form $(\dots ,b,b+1),$ then
the very even orbit is labelled $I.$

\subsubsection*{Algorithms for (c)}\ 

\medskip
\noindent\textit {$G$ of type C,D.}\quad  $\fz(\vO)=so(m_1)\times sp(m_2)\times
so(m_3)\times \dots$, where $m_i$ are the  number of parts
{of $\check\CO$} equal to $i.$ 

\noindent{\textit G of type B.}\quad  $\fz(\vO)=sp(m_1)\times so(m_2)\dots
sp(m_3)\times\dots $ where again $m_i$ is the number of parts
{of $\check\CO$} equal to $i.$

\subsubsection*{Algorithms for (d)}\ 

\medskip 
The parameter $\nu$ is a vector of size equal to $\rk G.$ For each
$\fk z_i,$ of $\fz(\C O),$ add $\rk \fk z_i$ coordinates each equal to the
average of the string cooresponding to the size $i$ part of $\vO.$ 
For each factor $\fz_i,$ of $\fz(\vO),$ the coordinates are the
averages of the corresponding strings.
The remaining  coordinates of $\nu$ are all zero.

\subsubsection*{Algorithm for (e)}\ 

\medskip
The Levi subgroups are determined by specifying the $GL$
factors. There is at most one other factor {$G_0(n_0)$} of the same
type as the group. 

For $\vm_{BC},$ each pair of parts $(k,k)$ yields a $GL(k).$ If the
corresponding string comes from Step 2, then the character on $GL(k)$
is given by $|\det(*)|^{\text{average of string}}.$ Otherwise it is the
trivial character. The parts of the remaining partition have
multiplicity 1 {corresponding to a distinguished orbit in
  $\check\fg_0(n_0)$}.  

For $\vm_{KL},$ apply the same procedure as for $\vm_{BC},$ except
for pairs coming from Step 1. If originally there was an odd number
of parts, then there is no change. If there was an even number, leave
behind one pair.  The parts in
the remaining partition have multiplicity 1 or 2 {
  corresponding to an even orbit in $\check\fg_0(n_0)$}.

\subsubsection*{Algorithms for (f) and (g)}\ 

\medskip
Both $\vm_e$ and $\vm_o$ acquire a $GL(k)$ factor for each pair of
parts $(k,k)$ in Step 2, with character given by the average of the
corresponding string as before. 

\medskip
For the parts coming from Step 1, write them in decreasing order
$a_r\ge \dots\ge  a_1>0.$

\medskip
\noindent\textit{For $\vm_e$,} there are additional {$GL$} factors
\begin{description}
\item[G of type B] $(a_1+a_2)/2,(a_3+a_4)/2,\dots ,(a_{r-2}+a_{r-1})/2$
  if $r$ is odd, $(a_1)/2,(a_2+a_3)/2,\dots ,(a_{r-2}+a_{r-1})/2$ when
  $r$ is even. The characters are given 
  by the averages of the strings $(-(a_{r-1}-1)/2,\dots ,(a_{r-2}-1)/2)$ and so
  on, and $(-(a_1-1)/2,\dots ,-1/2),\dots$ when $r$ is even. 
Recall that the $a_i$ are all even  because $\vO$ is a
  nilpotent orbit in type C. 
\item[G of type C] $(a_1+a_2)/2,(a_3+a_4)/2,\dots $, with 
  characters given by the averages of the strings $(-(a_1-1)/2,\dots ,
 (a_2-1)/2)$, and so on. In this case $\vO$ is type B, so there are an
 odd number of odd parts.
\item[G of type D] $(a_1+a_2)/2,(a_3+a_4)/2,\dots ,$ with characters
  obtained by the same procedure as in type C. In this case $\vO$ has
  an even number of odd parts. 
\end{description}

\medskip
\noindent\textit{For $\vm_o$,} there are additional {$GL$} factors

\medskip
\begin{description}
\item[G of type B] $(a_2+a_3)/2,(a_4+a_5)/2,\dots ,(a_{r-1}+a_{r})/2$
  leaving $a_1$ out if $r$ is odd, $(a_1+a_2)/2,(a_3+a_4)/2,\dots
  ,(a_{r-1}+a_{r})/2$ if $r$ is even. The characters are given   by
the averages of the strings \newline
  $(-(a_r-1)/2,\dots ,(a_{r-1}-1)/2)\dots $ and so on. 
\item[G of type C] $(a_2+a_3)/2,{(a_4+a_5)/2},\dots $, and
  characters given by the averages of the strings $(-(a_3-1)/2,\dots ,
  (a_2-1)/2),\dots$. In this case $\vO$ has  an odd
  number of odd sized parts. 
\item[G of type D] $(a_1+a_2)/2,\dots
  ,{(a_{r-3}+a_{r-2})/2}, ((a_{r-1}-1)/2)$ with characters 
obtained by the averages of the strings 
   $(-(a_2-1)/2,\dots ,(a_1-1)/2),\dots ,((-a_{r-1}-1)/2,\dots
,-1)$.   
In this case $\vO$ has  an even number of odd sized parts.
\end{description}

\subsection{}\label{sec:2.6}
Let $\chi=\vh/2 +\nu$ be associated to the orbit $\cCO.$ Recall from
\ref{sec:2.3} 
\begin{equation}
  \label{eq:2.6.1}
  I_M(\chi):=Ind_M^G[L_M(\chi)],
\end{equation}
where $L_M(\chi)$ is the irreducible spherical module of $M$ with
parameter $\chi.$ Write the nilpotent orbit in (\ref{eq:2.3.4})
with the $(d_1,\dots ,d_l)$ as in
sections \ref{sec:2.3b}-\ref{sec:2.3d} depending on the Lie
algebra type. Then $\vm_{BC}=gl(a_1)\times\dots\times
gl(a_k)\times \vg_0(n_0)$ is as in (\ref{eq:2.3.5}). 
Thus $\chi$ determines a spherical irreducible module
\begin{equation}
  \label{eq:2.6.2}
  L_{M_{BC}}(\chi)=L_1(\chi_1)\otimes\dots \otimes L_k(\chi_k)\otimes
  L_0(\chi_0), 
\end{equation}
with $\chi_i=(-\frac{a_i-1}{2}+\nu_i,\dots ,\frac{a_i-1}{2}+\nu_i)$, while
$\chi_0=\vh_0/2$ for the nilpotent $(d_i)$. 

\medskip
Let $\vm_{KL}$ be the Levi component attached to
$\chi=\vh/2+\nu$ in sections \ref{sec:2.3b}-\ref{sec:2.3d}. As for
$\vm_{BC}$ we have a parameter $L_{M_{KL}}(\chi)$. In this case
$\cCO=((a_1',a_1'),\dots, (a_{r}',a_{r}');d_1',\dots d_{l}')$ as
described in \ref{sec:2.3b}-\ref{sec:2.3d}. Then (a $W$-conjugate of)
$\chi$ can be written as in (\ref{eq:2.3b.2})-(\ref{eq:2.3d.2})), and
\begin{equation}
    \label{eq:2.6.3}
    \begin{aligned}
&    \vm_{KL}=gl(a_1')\times\dots\times gl(a_{r}')\times \vg_0(n_0'),\\
&    L_{M_{KL}}(\chi)=L_1(\chi_1')\otimes \dots
  \otimes L_r(\chi_r')\otimes L_0(\chi_0').      
    \end{aligned}
\end{equation}
  \begin{theorem}
    \label{t:2.6}
In the p-adic case
\begin{equation*}
  I_{M_{KL}}(\chi)=L(\chi).
\end{equation*}
  \end{theorem}
\begin{proof}
 This is in \cite{BM1}, $\vm_{KL}$ was defined in such a way that this
 result holds.
  \end{proof}
\begin{corollary}\label{c:2.6}
The module $I_{M_{BC}}(\chi)$ equals $L(\chi)$ in the $p-$adic case if
all the $\nu_i\ne 0.$ 
\end{corollary}

\medskip
\section{The Main Result}\label{sec:3}
 
\subsection{}\label{2.1} 
Recall that $\vG$ is the (complex) dual group, and $\vA\subset \vG$ the maximal
torus dual to $A.$ Assuming as we may that the parameter is real, a
spherical irreducible representation corresponds to an 
orbit of a hyperbolic element $\chi\in\va,$ the Lie algebra of $\vA.$
In section \ref{sec:2} we attached a nilpotent orbit $\vO$ in $\vg$ with
partition $(\unb{r_1}{a_1,\dots a_1},\dots,\unb{r_k}{a_k,\dots ,a_k})$
to such a parameter. Let $\{\ve,\vch,\vf\}$ be a Lie triple attached
to $\vO$.  Let  $\chi:=\vh/2+\nu$ satisfy (\ref{eq:2.3.1})-(\ref{eq:2.3.2}).

\begin{definition}
  \label{d:2.1}
A representation $L(\chi)$ is said to be in the complementary series
for $\cCO,$ if the parameter $\chi$ is attached to $\cCO$ in the sense of
satisfying (\ref{eq:2.3.1}) and (\ref{eq:2.3.2}), and is unitary.
\end{definition}
We will describe the complementary series explicitly in coordinates.

The centralizer $Z_\vG(\ch e,\ch h,\vf)$ has Lie algebra $\fk z(\vO)$
which is a product of $sp(r_l,\bC)$ or $so(r_l,\bC),$ $1\le l\le k$,
according to the  rule
\begin{description}
\item[${\mathbf\vG}$ of type B, D] $sp(r_l)$ for $a_l$ even, 
$so(r_l)$ for $a_l$ odd,
\item[${\mathbf \vG}$ of type C]    $sp(r_l)$ for $a_l$ odd, 
$so(r_l)$ for $a_l$ even.
\end{description}
{On $\fk z(\vO)$, $\nu$
determines  a spherical irreducible module $L_{\cCO}(\nu)$ for  the
split group whose dual is $Z_{\vG}(\ve,\vh,\vf)^0.$  Applying
(\ref{eq:2.3.1}-\ref{eq:2.3.2}), we find that $\nu$ is attached to 
the trivial orbit in $\fk z (\cCO).$ If it were not, then
$\nu=\vh'/2+\nu',$ for a triple $\{\ve',\vh',\vf'\},$ with $\ve'\ne 0,$ so 
$\chi=(\vh/2+\vh'/2) +\nu'$,  and $\{\ve+\ve',\vh+\vh',\vf+\vf'\}$ is a
Lie triple such that the orbit of $\ve +\ve'$ is strictly larger than
the orbit of $\ve,$ and contains it in its closure. 
}
\begin{theorem}\label{thm:3.1} The complementary series attached
to $\vO$ coincides with the one attached to the
trivial orbit in $\fk z(\vO).$ For the trivial orbit $(0)$ in each of
the classical cases, the complementary series are
\begin{description}
\item[${\mathbf G}$ of type B]
$$
0\le \nu_1\le\dots\le \nu_k<1/2.
$$  
\item[${\mathbf G}$ of type C, D] 
$$
0\le \nu_1\le \dots\le \nu_k\le 1/2<\nu_{k+1}<\dots
<\nu_{k+l}<1 
$$
so that $\nu_i+\nu_j\ne 1.$ There are 
\begin{enumerate}
\item an even number of $\nu_i$ such that $1-\nu_{k+1}<\nu_i\le 1/2,$
\item {for every $1\le j\le l$,} there is an odd number of $\nu_i$
  such that $1-\nu_{k+j+1}<\nu_i<1-\nu_{k+j}.$   
\item  In type D of odd rank, $\nu_1=0$ or else the parameter is not hermitian. 
\end{enumerate}

\end{description}
\end{theorem}
\subsection*{Remarks}
\begin{enumerate}
\item The complementary series for $\vO=(0)$ consists of
  representations which are both spherical and generic in the sense
  that they have Whittaker models.
\item The condition that $\nu_i+\nu_j\ne 1$ implies that in types C,D there
is at most one  $\nu_k=1/2$.
\item In the case of $\vO\ne (0),$ $\chi=\vh/2 +\nu,$ and each
of the coordinates $\nu_i$ for the parameter on $\fz(\vO)$ comes from a 
 string, \ie each $\nu_i$ comes  from  $(-\frac{a_i-1}{2}+\nu_i,\dots
  ,\frac{a_i-1}{2}+\nu_i)$. { The
parameter does not satisfy (\ref{eq:2.3d.11}). For (\ref{eq:2.3d.11}) to
hold, it suffices to change $\nu_{k+j}$ for types $C,\ D$ to
$1-\nu_{k+j}.$} {More precisely, for $1/2<\nu_{k+j}<1$ the
connection with the  strings in the form (\ref{eq:2.3d.10}) and
(\ref{eq:2.3d.11}) is as follows. Write $(-\frac{a_{k+j}-1}{2}+\nu_{k+j},\dots
  ,\frac{a_{k+j}-1}{2}+\nu_{k+j})$ as $(-\frac{a_{k+j}-3}{2}+(\nu_{k+j}-1),\dots
  ,\frac{a_{k+j}+1}{2}+(\nu_{k+j}-1))$ and then conjugate each entry
to its negative to form $(-\frac{a_{k+j}-3}{2}+\nu_{k+j}',\dots
  ,\frac{a_{k+j}+1}{2}+\nu_{k+j}')$, with $0<\nu_{k+j}'=1-\nu_{k+j}<1/2.$}
\end{enumerate}
\subsection*{An Algorithm}
{For types C and D we give an algorithm, due to S. Sahi, to decide
whether a parameter in types C and D is unitary. This algorithm is for
the complementary series for $\vO=(0).$ For arbitrary $\vO$ it
applies to the parameter for $\fz(\vO)$ obtained as in remark (3)
above.}

\medskip
Order the parameter in dominant form,
\begin{equation}
  \label{eq:3.1.1}
  \begin{aligned}
&0\le \nu_1\le \dots\le \nu_n, \text{ for type C},\\    
&0\le |\nu_1|\le \dots\le \nu_n, \text{ for type D}.
  \end{aligned}
\end{equation}
The first condition is that $\nu_n<1,$ and in addition that if the
type is $D_n$ with $n$ odd, then $\nu_1=0.$ Next replace each
coordinate $1/2< \nu_i$ by $1-\nu_i.$ Reorder the
new coordinates in increasing order 
as in (\ref{eq:3.1.1}).
{Let $F(\nu)$ be the set of new positions of
the $1-\nu_i.$ If any position is ambiguous, the parameter is not
unitary, or is attached to a different nilpotent orbit. 
This corresponds to either a $\nu_i+\nu_j=1,$ or a
$1/2<\nu_i=\nu_j.$ Finally, $L(\nu)$ is unitary
if and only if $F(\nu)$ consists of odd numbers only.}

\subsection{}\label{sec:3.2} We prove the unitarity of the parameters
in the theorem for $\vO=(0)$ for types B,C, and D. First we record
some facts.

\medskip
Let $G:=GL(2a)$  and 
\begin{equation}
  \label{eq:3.2.1}
 \chi:= (-\frac{a-1}{2}-\nu,\dots ,\frac{a-1}{2}-\nu;
   -\frac{a-1}{2}+\nu,\dots ,\frac{a-1}{2}+\nu).
\end{equation}
Let $M:=GL(a)\times GL(a)\subset GL(2a)$. Then the two strings of
$\chi$ determine an irreducible spherical (1-dimensional) representation
$L_M(\chi)$ on $M$. Recall $I_M(\chi):=Ind_M^G[L_M(\chi)].$ 
\begin{lemma}[1]\label{l1:3.2} The representation $I_M(\chi)$ is unitary
  irreducible for $0\le\nu<1/2.$ The irreducible spherical module
  $L(\chi)$ is not unitary for $\nu>\frac12,\ 2\nu\notin \bZ.$  
\end{lemma}
\begin{proof} This is well known and goes back to \cite{Stein} 
(see also  \cite{T} and \cite{V1}). 
\end{proof}
We also recall the following well known result due to Kostant in the
real case, Casselman in the $p-$adic case.
\begin{lemma}[2]\label{l2:3.2}
If none of the $\langle\chi,\al\rangle$ for $\al\in\Delta(\vg,\va)$ is a
nonzero integer, then $X(\chi)$ is irreducible.
In particular, if $\chi=0,$ then 
\begin{equation*}
  L(\chi)=X(\chi)=Ind_A^G[\chi],
\end{equation*}
and it is unitary.
 \end{lemma}
Let $\vm\subset\vg$ be a Levi component, and $\xi_t\in\fk z(\vm)$,
where $\fk z(\vm)$ is the center of $\vm,$ depending continuously on
$t\in [a,b]$.
{\begin{lemma}[3]
  \label{l3:3.2}
Let $\xi_t$ be a character of $M,$ depending continuously on
$t\in\bR.$  Assume that
\begin{equation*}
  I_M(\chi_t):=Ind_M^G[L_M(\eta_0)\otimes\xi_t]
\end{equation*}
is irreducible and hermitian for $a\le t\le b$. 
Then if $I_M(\chi_t)$ (equal to $L(\chi_t)$) is unitary for some $a\le
t_0\le b,$ it is unitary for all $a\le t\le b.$  
\end{lemma}

This is well known, and amounts to the fact that if a hermitian matrix
is nondegenerate and depends continuously on a parameter $a\le t\le
b$, if it is positive definite for some $a\le t_0\le b,$ it is
positive definite throughout the interval; to change from positive to
negative, it would have to go through a zero, \ie become degenerate. 
This happens for example when
$L_M(\eta_0)\otimes\xi_{t_0}$ is unitary for some $t_0,$ and
$I_M(\chi_t)$ is irreducible. I don't know the original reference.
}

When the conditions of lemma 3 are satisfied, we say that
$I_M(\chi_t)$ is a \textit{continuous deformation} of $L(\chi_0).$ 

\medskip
We now start the proof of the unitarity.
\subsubsection*{Type B} In this case there are no
roots $\al\in\Delta(\vg,\va)$ such that $\langle\chi,\al\rangle$ is a nonzero
integer. Thus 
\begin{equation*}
L(\chi)=Ind_A^G[\chi]  
\end{equation*}
as well. When deforming $\chi$ to $0$ continuously, the induced module stays
irreducible. Since $Ind_A^G[0]$ is unitary, so is $L(\chi).$

\subsubsection*{Type C,D} There is no root such that 
$\langle\chi,\al\rangle$ is a nonzero integer, so
$L(\chi)=Ind_A^G[\chi]$. If there are no $\nu_{k+i}>1/2$ the argument
for type B carries over word for word. 
When  there are $\nu_{k+i}>1/2$ we have to be more careful with the
deformation. We will do an induction on the rank.
{Suppose that $\nu_{j-1}=\nu_j$ for some $j$.
  Necessarily, $\nu_j<1/2.$}
Conjugate $\chi$ by the Weyl group so that 
\begin{equation}
  \label{eq:3.2.2}
  \chi=(\nu_1,\dots ,\nu_i,\dots \wht{\nu_{j-1}},\wht{\nu_j},\dots ,
  ;\nu_{j-1};\nu_j):=(\chi_0;\nu_{j-1};\nu_j). 
\end{equation}
Let $\vm:=\vg(n-2)\times gl(2),$ and denote by $M$ the corresponding
Levi component. Then  by induction in stages,
\begin{equation}
  \label{eq:3.2.3}
L(\chi)=Ind_M^G[L_M(\chi)],
\end{equation}
where $L_M(\chi)=L_0(\chi_0)\otimes L_1(\nu_{j-1},\nu_j).$ By lemma (1) of
\ref{l1:3.2}, $L_1(\nu_{j-1},\nu_j)$ is unitary. Thus $L(\chi)$ is unitary
if and only if $L_0(\chi_0)$ is unitary. If $\chi$ satisfies the assumptions of
the theorem, then so does $\chi_0$. By the induction hypothesis,
$L_0(\chi_0)$ is unitary, and therefore so is $L(\chi).$ 
Thus we may assume that 
\begin{equation}
  \label{eq:3.2.4}
0\le\nu_1<\dots <\nu_k\le 1/2< \nu_{k+1}<\dots <\nu_{k+l}.
\end{equation}
If $\nu_k<1-\nu_{k+1}$, then the assumptions imply
$1-\nu_{k+2}<\nu_k$. Consider the parameter
\begin{equation}
  \label{eq:3.2.5}
  \chi_t:=(\dots ,\nu_k,\nu_{k+1}-t,\dots).
\end{equation}
Then
\begin{equation}
  \label{eq:3.2.6}
  L(\chi_t)=Ind_A^G[\chi_t],\qquad\text{ for } 0\le t\le \nu_{k+1}-\nu_k,
\end{equation}
because no $\langle\chi_t,\al\rangle$ is a nonzero integer.
At $t=\nu_{k+1}-\nu_k,$ the parameter is in the case just considered
earlier. By induction we are done. 

If on the other hand $1-\nu_{k+1}<\nu_k,$ the assumptions on the
parameter are such that necessarily $1-\nu_{k+1}<\nu_{k-1}<\nu_{k}.$  
Then repeat the argument with 
\begin{equation}
  \label{eq:3.2.7}
  \chi_t:=(\dots ,\nu_{k-1},\nu_k-t,\dots ),\qquad 0\le t\le \nu_k-\nu_{k-1}.
\end{equation}
This completes the proof of the unitarity of the parameters in theorem
\ref{thm:3.1} when $\vO=(0).$

\subsection{}\label{sec:3.3}
We prove the unitarity of the parameters in theorem \ref{thm:3.1} in
the general case when $\vO\ne (0)$.    

The proof is essentially the same as for
$\vO=(0),$ but special care is needed to justify the irreducibility of
the modules. Recall the notation of the partition of $\vO$ (\ref{eq:2.3.3}).

The factors of $\fk z(\vO)$ isomorphic to $sp(r_j),$
contribute $r_j/2$ factors of the form $gl(a_i)$ to $\vm_{KL}.$ The
factors of type $so(r_j)$ with $r_j$ odd, contribute a $d_i$ (notation
(\ref{eq:2.3.4})) to the
expression (\ref{eq:2.3.4}) of the partition of $\vO,$  and
$\frac{r_j-1}{2}$ $gl(a_i)$. The factors $so(r_j)$ of type D ($r_j$
even)  are more complicated. Write the strings coming from this factor as
in (\ref{eq:2.3.6}), 
$$
(-\frac{a_i-1}{2}+\nu_i,\dots , \frac{a_i-1}{2}+\nu_i)
$$
with the $\nu_i$ satisfying the assumptions of theorem \ref{thm:3.1}.
If $r_j$ is not divisible by 4, then there must be a $\nu_1=0,$
(otherwise the corresponding spherical parameter is not hermitian), and
$\vm_{BC}\ne\vm_{KL}$. Similarly  when $r_j$ is divisible by 4 and
$\nu_1=0,$ $\vm_{BC}\ne\vm_{KL}$. In all situations, we consider
\begin{equation}
  \label{eq:3.3.1}
  I_{M_{KL}}(\chi)
\end{equation}
as for $\vO=(0).$ We aim to show that this module stays irreducible
under the deformations used for $\vO=(0),$ separately for the $\nu_i$
for the same partition size or equivalently simple factor of
$\fz(\vO)$. {It is enough to prove that under these
  deformations the strings stay  strongly nested at all the values of the
parameters. Then proposition \ref{p:3.3} applies.}

Using the conventions of section
\ref{sec:2}, the strings are of the form
\begin{align}
&(-A-1+\nu_1,\dots , A-1+\nu_1)  \label{eq:3.3.1a}\\
&(-B+\nu_2,\dots, B-1+\nu_2)  \label{eq:3.3.1b}\\
&(-C+ \nu_3,\dots ,C+\nu_3).  \label{eq:3.3.1c}  
\end{align}
{The $\nu_i$ satisfy $0\le \nu_i\le 1/2.$ {When
    $\nu_1=1/2,$ the string is $(-A-1/2,\dots , A-1/2)$ 
so it conforms to (\ref{eq:2.3b.11}),
    (\ref{eq:2.3c.11}), (\ref{eq:2.3d.11}). Similarly when $\nu_2=1/2,$ the
    string is $(-B+1/2,\dots ,B-1/2)$. But
when $\nu_3=1/2,$  the string is $(-C+1/2,\dots
    ,C+1/2)$ and it must be replaced by $(-C-1/2,\dots ,C-1/2)$ to
    conform to (\ref{eq:2.3b.11}),
    (\ref{eq:2.3c.11}), (\ref{eq:2.3d.11}).} 
{ The string in (\ref{eq:3.3.1a})  gives a
  $\nu_{1,\fz(\vO)}=1-\nu_1$ {(as explained in remark (3) of section
  \ref{2.1})}.  The other ones $\nu_2$ and $\nu_3$ give
$\nu_{j,\fz(\vO)}=\nu_2$ or $\nu_3$ respectively. 
Suppose first that there is only one size of string. {
    This means that the corresponding entries $\nu_{j,\fz(\vO)}$ belong
    to the same simple factor of $\fz(\vO).$} Then the strings
are either all of the form (\ref{eq:3.3.1a}) and (\ref{eq:3.3.1c})
with $A=C$ or all of the form (\ref{eq:3.3.1b}). Consider the first
case. If there is a string (\ref{eq:3.3.1a}) then
$\nu_{1,\fz(\vO)}=1-\nu_1\ge 1/2,$ and  $\nu_{1,\fz(\vO)}$  is deformed
downward.  By the assumptions $\nu_{1,\fz(\vO)}$ does not equal any
$\nu_{j,\fz(\vO)}$ nor does $\nu_{1,\fz(\vO)}+\nu_{j,\fz(\vO)}=1$ for
any $j.$ The module stays irreducible. When $\nu_{1,\fz(\vO)}$ crosses
$1/2,$ the string becomes $(-A +\nu_1',\dots , A+\nu_1'),$ and
$\nu_1'$ is deformed downward from $1/2$ to either some $\nu_3$ or to
$0.$ Again $\nu_1'+\nu_{j,\fz(\vO)}\ne 1$ for any $j,$ so no
  reducibility occurs. The irreduciblity in the case when $\nu_1'$ reaches $0$  
is dealt with by section \ref{sec:10}. 

Remains to check that in these deformations no reducibility occurs
because the string interacts with one of a different length,
{in other words, when a 
  $\nu_{i,\fz(\vO)}$ for one size string or equivalently factor of
  $\fz(\vO)$ becomes equal to a $\nu_{j,\fz(\vO)}$ from a 
  distinct factors of $\fz(\vO).$ We explain the case when the
  deformation involves a string of type (\ref{eq:3.3.1a}), the others
  are similar and easier.} 
Consider the deformation of $\nu_1$ from $0$ to $1/2.$ If the module
is to become reducible $\nu_1$ must reach a $\nu_3$ so that there is 
a string of the form (\ref{eq:3.3.1c}) satisfying
\begin{equation}
  \label{eq:3.3.1d}
  -A-1<-C,\qquad A-1<C.
\end{equation}
This is the condition that the strings are not nested for some value
of the parameter. 
This  implies that $A-1<C<A+1$ must hold. Thus $A=C,$ and the two strings
correspond to the same size partition or simple factor in $\fz(\vO).$
}

{\medskip
In the p-adic case, the irreducibility of (\ref{eq:3.3.1}) in the case
of strongly nested strings follows from the results of Kazhdan-Lusztig. 
In the case of real groups, the same irreducibility results hold, but
are harder to prove.  Given $\chi$, consider the root system
\begin{equation}
  \label{eq:3.3.5}
 \check \Delta_\chi:=\{\cha\in \check\Delta : \langle\chi,\cha\rangle\in \bZ\}.
\end{equation}
Let $G_\chi$ be the connected split real group whose dual root system
is $\check\Delta_\chi.$ 

The Kazhdan-Lusztig-Vogan conjectures for nonintegral
infinitesimal character relate statements about the character theory
of admissible $(\fk g_c, K)-$modules of $G$ with infinitesimal
character $\chi$  to similar statements for characters of admissible
modules for $G_\chi$ but at integral infinitesimal character, for
which the Kazhdan-Lusztig-Vogan conjectures for regular integral
infinitesimal character can be used.  
This is beyond the scope of this paper (or my competence), I refer to
\cite{ABV}, chapters 16 and 17 for an explanation. 

In short, $\chi$ determines an irreducible spherical representation
$L(\chi),$ and a spherical  $L_{G_\chi}(\chi).$ To prove the
irreducibility results we need, we will prove them for
$L_{G_\chi}(\chi),$ where we may assume that the infinitesimal
character is integral. Since $G_\chi$ is not simple, it is sufficient
to prove the needed irreducibility result for each simple factor. 
This root system is a product of classical systems as follows.
For each $0\le \tau\le 1/2$ let $A_\tau$ be the set of
  coordinates of $\chi$ congruent to $\tau$ modulo $\bZ.$  
Each $A_\tau$ contributes to $\check\Delta_\chi$ as follows.
\begin{description}
\item[$\mathbf G$ of type B] Every $0<\tau<1/2$ contributes a type A
  of size equal to the number of coordinates in $A_\tau$. Every
  $\tau=0,1/2$   contributes a type C of rank equal to the number of
  coordinates in $A_\tau.$
\item[$\mathbf G$ of type C] Every $0<\tau<1/2$ contributes a type
  A as for type B. Every $\tau=0$ contributes a type $B,$ while $\tau=1/2$
  contributes type D. 
\item[$\mathbf G$ of type D]  Every $0<\tau<1/2$ contributes a type A. Every
  $\tau=0,1/2$   contributes a type D.
\end{description}
}

The irreducibility results for $I_{M_{KL}}(\chi_t)$ needed to carry out
the proof are contained in the following proposition.
\begin{proposition}\label{p:3.3}
Let
\begin{equation*}
\chi:=(\dots ;-\frac{a_i-1}{2}+\nu_i,\dots , \frac{a_i-1}{2}+\nu_i;\dots )  
\end{equation*}
be given in terms of strings, and let $\vm=gl(a_1)\times\dots\times
gl(a_k)$ be the corresponding Levi component. Assume that $\chi$ is
integral (\ie $\langle \chi,\cha\rangle\in \bZ$).  
In addition assume that the coordinates of $\chi$
are 
\begin{itemize}
\item in $\bZ,$ in type B,

\item in $1/2+\bZ$ for type D. 
\end{itemize}
\noindent If the strings are strongly nested, then
\begin{equation*}
  L(\chi)=Ind_{M}^G[L_M(\chi)]. 
\end{equation*}
\end{proposition}
The proof of the proposition will be given in section
\ref{sec:10}.
{Unlike in type B and D where $\langle\chi,\cha\rangle\in
  \bZ$ implies that the coordinates are integers or half integers, in
  type C, the condition implies that the coordinates are integers only.}

\subsection*{Remark} For the $\nu_j$ attached to factors of type D in
$\fk z(\vO)$, it is important in the argument that we do not deform to
$(0).$ The next example illustrates why.  

Assume $\vO=(2,2,2,2)\subset sp(4)$. The parameters of the form $\vh/2+\nu$ are 
\begin{equation}
  \label{eq:3.3.2}
  (-1/2+\nu_1,1/2+\nu_1;-1/2+\nu_2,1/2+\nu_2),
\end{equation}
and, because parameters are up to $W-$conjugacy, we may restrict
attention to the region $0\le \nu_1\le\nu_2.$ In 
this case $\fk z(\vO)=so(4),$ and the unitarity region is
$0\le\pm\nu_1+\nu_2<1.$   
Furthermore $\vm_{BC}=gl(2)\times gl(2),$ but $\vm_{KL}=\vm_{BC}$ only if
$0<\nu_1.$ When $\nu_1=0,$ $\vm_{KL}=sp(2)\times gl(2),$ the
nilpotent orbit is rewritten $(2,2; (2,2)),$ and $\vh_0/2=(1/2,1/2).$
For $\nu_1=0,$ the induced representations
\begin{equation}
  \label{eq:3.3.3}
  I_{M_{KL}}(\chi_{\nu_2}):=Ind_{Sp(2)\times
    GL(2)}^{Sp(4)}[L_0(1/2,1/2)\otimes L_1(-1/2+\nu_2,1/2+\nu_2)]
\end{equation}
are induced irreducible in the range $0\le \nu_2<1.$ For $0<\nu_1$ the
representation
\begin{equation}
  \label{eq:3.3.4}
  Ind_{GL(4)}^{Sp(4)}[L((-1/2+\nu_1+t,1/2+\nu_1+t);(-1/2-\nu_1-t,1/2-\nu_1-t))]
\end{equation}
is induced irreducible for $0\le t\le 1/2-\nu_1.$

The main point of the example is that
{$Ind_{GL(2)}^{Sp(2)}[L(-1/2+t,1/2+t)]$ is
  \textbf{reducible} at $t=0.$}
So we cannot conclude that $L(\chi)$ is unitary for a $(\nu_1,\nu_2)$
with $0<\nu_1$ from the unitarity of $L(\chi)$ for a parameter with
$\nu_1=0.$ Instead we conclude that the representation is unitary
in the region $0\le \pm\nu_1+\nu_2<1$ because it is a deformation of
the irreducible module for $\nu_1=\nu_2$ which is unitarily induced
irreducible from a Stein complementary series on $GL(4).$\qed 

\medskip

\section{Relevant $K-$types}\label{sec:4}

\medskip
{In this chapter we define a special set of $K-$types (occuring in
the spherical principal series) which we call relevant. Their first
important property is that the intertwining operators used to compute
the hermitian form are particularly simple. This is the property of
\textit{petite}. A more general notion was introduced and studied
independently by Oda, \cite{oda}.  In the notation of chapter
\ref{sec:5}  these operators only depend on the $W-$structure of
$V^{K\cap B}.$ In this chapter we determine these Weyl group
representations for the relevant $K-$types. In chapter \ref{sec:6} 
we carry out the intertwining operator calculations in terms of the Weyl group
representations. These calculations are used to obtain necessary
conditions for unitarity in chapter \ref{sec:7}. 
But in chapters \ref{sec:9} and \ref{sec:10} they also prove to be
crucial for determining unitarity and irreducibility of certain
modules needded to complete the determination of the unitary dual.
}
\subsection{}
In the real case we will call a $K-$type $(\mu,V)$
\textit{quasi-spherical} if it 
occurs in the spherical principal series.
{In section \ref{sec:4}, we will use the notation $M=K\cap B.$} 
 By Frobenius reciprocity
$(\mu,V)$ is quasi-spherical if and only if $V^{K\cap B}\ne 0$. 
Because the Weyl group $W(G,A)$ may be realized as $N_K(A)/Z_K(A),$
this Weyl group acts naturally  on this space.

\medskip
The representations of $W(A_{n-1})=S_n$ are parametrized by partitions
$(a):=(a_1,\dots, a_k), \ a_i\le a_{i+1},$ of $n,$ and we write
$\sig((a))$ for the corresponding representation.  
The representations of $W(B_n)\cong W(C_n)$ are parametrized as in \cite{L1} by
pairs of partitions, and we write as 
\begin{align}
  \label{eq:4.0.4}
\sig((a_1,\dots , a_r)&,(b_1,\dots,b_s)),\notag \\
  &a_i\le a_{i+1},\quad b_j\le b_{j+1},\quad \sum a_i +\sum b_j=n.
\end{align}
Precisely the representation parametrized by (\ref{eq:4.0.4}) is as
follows. Let $k=\sum a_i,\ l=\sum b_j.$ Recall that $W\cong S_n\ltimes
\bZ_2^n.$ Let $\chi$ be the character of $\bZ_2^n$ which is trivial on
the first $k$ $\bZ_2$'s, sign on the last $l.$ Its centralizer in $S_n$ is
$S_k\times S_l.$ Let $\sig((a))$ and $\sig((b))$ be the representations of
$S_k,\ S_l$ corresponding to the partitions $(a)$ and $(b).$ Then let
$\sig((a),(b),\chi)$ be the unique representation of $(S_k\times
S_l)\ltimes Z_2^n$ which is a multiple of $\chi$ when restricted to
$Z_2^n,$ and $\sig(a)\otimes \sig(b)$ when restricted to $S_k\times
S_l.$ 
The representation in (\ref{eq:4.0.4}), is
\begin{equation}
  \label{eq:4.0.5}
\sig((a),(b))=Ind_{(S_k\times S_l)\ltimes \bZ_2^n}^{W}\
[\sig(a,b,\chi)].
\end{equation}
If $(a)\ne(b),$ the representations $\sig((a),(b))$ and $\sig((b),(a))$ restrict
to the same  irreducible representation of  $W(D_n),$ which we denote
again by the same symbol. When $a=b,$ the restriction is a sum of two
inequivalent representations which we denote $\sig((a),(a))_{I,\ II}.$ 
Let $W_{(a),I}:=S_{a_1}\times\dots\times S_{a_r}$ and $W_{(a),II}:=S_{a_1}\times
\dots\times S_{a_r}',$ be the Weyl groups corresponding to the Levi
components considered in Remark (2) in section \ref{sec:2.3d}. 
Then $\sig((a),(a))_I$ is characterized by the fact that its restriction to
$W_{(a),I}$ contains the trivial representation. Similarly
$\sig((a),(a))_{II}$ is the one that contains the trivial representation
of $W_{(a),II}.$

\medskip
\subsection{Symplectic Groups}\label{sec:4.2}
The group is $Sp(n)$ and the maximal compact subgroup is $U(n).$ The
highest weight of a $K$-type will be written as $\mu(a_1,\dots ,a_n)$
with $a_i\ge a_{i+1}$ and $a_i\in\bZ$, or
\begin{equation}
  \label{eq:4.2.1}
  \mu(a_1^{r_1},\dots ,a_k^{r_k}):=(\unb{r_1}{a_1,\dots ,a_1},\dots
  ,\unb{r_k}{a_k,\dots ,a_k}).
\end{equation}
when we want to emphasize the repetitions. We will repeatedly use the
following restriction formula 
\begin{lemma}
  \label{l:4.2}
The restriction of $\mu(a_1,\dots ,a_n)$ to $U(n-1)\times U(1)$ is 
\begin{equation*}
  \sum \mu(b_1,\dots , b_{n-1})\otimes\mu(b_n), 
\end{equation*}
where the sum ranges over all possible $a_1\ge b_1\ge a_2\ge \dots \ge
b_{n-1}\ge a_n,$ and $b_n=\sum_{1\le i\le n} a_i-\sum_{1\le j\le n-1} b_j.$
\end{lemma}

\begin{definition}\label{def:4.2}
The representations $\mu_e(n-r,r):=\mu(2^r,0^{n-r})$ and
$\mu_o(k,n-k):=\mu(1^k,0^{n-2k},-1^k)$ are called \textbf{relevant}.   
\end{definition}

\begin{proposition}\label{p:4.2}
The relevant $K-$types are quasispherical. 
The representation of $W(C_n)$ on $V^M$ is 
\begin{equation*}
  \begin{aligned}
  &\mu_e(n-r,r)\longleftrightarrow \sig[(n-r), (r)],\\
  &\mu_o(k,n-k)\longleftrightarrow \sig[(k,n-k), (0)],\\    
  \end{aligned}
\end{equation*}
\end{proposition}
The $K$-types $\mu(0^{n-r},(-2)^{r})$, dual to $\mu_e(n-r,r)$  are also
quasispherical, and could be used in the same way {as $\mu_e(n-r,r)$}.
\begin{proof}
We do an induction on $n.$ {When $n=1,$ the only relevant
  representations of $U(1)$ are $\mu_e(1,0)=\mu(0)$ and
  $\mu_e(0,1)=(2)$, which correspond to the trivial and the sign
  representations of $W(C_1)=\bZ/2\bZ$, respectively.}
Consider the case $n=2.$ There are four relevant representations of
$U(2)$ with highest weights $(2,0),$  $(1,-1),$ $(2,2)$ and $(0,0).$ The first
representation is the 
symmetric square of the standard representation, the second one is the
adjoint representation and the fourth one is the trivial
representation. The normalizer of $A$ in $K$ can be
identified with the diagonal subgroup $(\pm 1, \pm 1)$ inside
$U(1)\times U(1)\subset U(2).$ The Weyl group is generated by the elements 
\begin{equation}
  \label{eq:4.2.3}
\begin{bmatrix}i&0\\ 0&1\end{bmatrix},\quad
\begin{bmatrix}1&0\\ 0&i\end{bmatrix},\quad
\begin{bmatrix}0&1\\-1&0\end{bmatrix}.
\end{equation}
The restriction to $U(1)\times U(1)$ of the four representations of
$U(2)$ is 
\begin{align}
  \label{eq:4.2.4}
  (2,0)&\longrightarrow (2)\otimes(0)+(1)\otimes(1) +
  (0)\otimes(2),\notag\\
  (1,-1)&\longrightarrow (1)\otimes(-1)+(0)\otimes(0) +
  (-1)\otimes(1),\notag\\ 
  (2,2)&\longrightarrow (2)\otimes(2),\\
  (0,0)&\longrightarrow (0)\otimes(0).\notag
\end{align}
The space $V^M$ is the sum of all the weight spaces $(p)\otimes (q)$
with both $p$ and $q$ even. For the last one, the representation of
$W$ on $V^M$ is $\sig[(2),(0)].$  The third one is 1-dimensional so
$V^M$ is 1-dimensional; the Weyl group representation is
$\sig((0),(2)).$ The second one has $V^M$ 1-dimensional and the Weyl
  group representation 
is $\sig((11), (0)).$ For the first one, $V^M$ is 2-dimensional and the
Weyl group representation is $\sig((1),(1)).$ These facts can be read
off from explicit realizations of the representations.

\smallskip
Assume  that the claim is proved for $n-1.$ Choose  a parabolic
subgroup  so that its Levi component is 
$M'=Sp(n-1)\times GL(1)$ and $M$ is contained in it. Let
$H=U(n-1)\times U(1)$ be such that $M\subset M'\cap K\subset H.$

Suppose that $\mu$ is relevant. The cases when $k=0$ or $r=0$ are
1-dimensional and are straightforward. So we only consider $k,\ r>0.$
The K-type $\mu(2^r,0^{n-r})$ restricts to the sum of
\begin{align}
&\mu(2^r,0^{n-r-1})\otimes\mu(0)
\label{eq:4.2.5}\\
&\mu(2^{r-1},1,0^{n-r-1})\otimes\mu(1)
\label{eq:4.2.6}\\
&\mu(2^{r-1},0^{n-r})\otimes\mu(2).
\label{eq:4.2.7}
\end{align}
Of the representations appearing, only
$\mu(2^r,0^{n-r-1})\otimes\mu(0)$ and\\
$\mu(2^{r-1},1,0^{n-r-1})\otimes\mu(2)$ are quasispherical. 
So the restriction of $V^M$ to $W(C_{n-1})\times W(C_1)$ is the sum of
\begin{align}
&\sig[(n-r-1),(r)]\otimes\sig[(1),(0)]\label{eq:4.2.8}\\
&\sig[(n-r),(r-1)]\otimes\sig[(0),(1)]\label{eq:4.2.9}
\end{align}
The only irreducible representations of $W(C_n)$ containing (\ref{eq:4.2.8}) in
their restrictions to $W(C_{n-1})$ are
\begin{align}
&\sig[(1,n-r-1), (r)]\label{eq:4.2.10}\\
&\sig[(n-r),(r)].\label{eq:4.2.11}
\end{align}
But the restriction of $\sig[(1,n-r-1),(r)]$ to $W(C_{n-1})\times W(C_1)$ 
contains \newline $\sig[(1,n-r-1),(r-1)]\otimes\sig[(0),(1)],$ and
this does not appear in  (\ref{eq:4.2.8})-(\ref{eq:4.2.9}).
Thus the representation of $W(C_n)$ on $V^M$ for (\ref{eq:4.2.10})
must be (\ref{eq:4.2.6}), and the claim is proved in this case. 

Consider the case $\mu(1^k,0^{l},-1^{k})$ for $k>0,\ 2k+l=n.$
The restriction of this
K-type to $U(n-1)\times U(1)$ is the sum of
\begin{align}
&\mu(1^k,0^l,-1^{k-1})\otimes\mu(-1) \label{eq:4.2.12}\\
&\mu(1^{k-1},0^l,-1^k)\otimes\mu(1)\label{eq:4.2.13}\\
&\mu(1^{k-1},0^{l+1},-1^{k-1})\otimes\mu(0)\label{eq:4.2.14}\\
&\mu(1^k,0^{l-1},-1^{k})\otimes\mu(0)\label{eq:4.2.15}
\end{align}
Of the representations appearing, only (\ref{eq:4.2.14}) and
(\ref{eq:4.2.15}) are quasispherical.  So the restriction of $V^M$ to
$W(C_{n-1})\times W(C_1)$ is the sum of  
\begin{align}
  &\sig[(k-1,k+l),(0)]\otimes\sig[(1),(0)],\label{eq:4.2.16}\\
  &\sig[(k,k+l-1),(0)]\otimes\sig[(1),(0)].\label{eq:4.2.17}
\end{align}
The representation (\ref{eq:4.2.17}) can only occur in  the
restriction to $W(C_{n-1})\times W(C_1)$ of $\sig[(1,k,k+l-1),(0)]$ or
$\sig[(k,k+l),(0)].$ If $k>1$, 
the first one contains $\sig[(1,k-1,k+l-1),(0)]$ in its restriction,
which is not in the sum of (\ref{eq:4.2.16}) and (\ref{eq:4.2.17}). 
If $k=1$ then (\ref{eq:4.2.16}) can only occur in  the restriction of 
$\sig[(0,l+2),(0)],$ or
$\sig[(1,l+1),(0)]$. But  $V^M$ cannot consist of $\sig[(0,l+2),(0)]$
alone, because (\ref{eq:4.2.16}) does not occur in its restriction. If it
consists of both $\sig[(0,l+2),(0)]$ and $\sig[(1,l),(0)],$ then the
restriction is too large. The claim is proved in this case. 
\end{proof}

\subsection{Orthogonal groups}\label{sec:4.3}
Because we are dealing with the spherical case, we can use the groups
$O(a,b)$, $SO(a,b),$ or the connected component of the identity,
$SO_e(a,b)$. {The corresponding $K$'s are $O(a)\times
  O(b)$, $S(O(a)\times O(b))$, and $SO(a)\times SO(b)$, respectively.}
We will use $O(a,b)$ {for the calculation of relevant
  $K$-types.}   For $SO(a),$ an
irreducible representation will be
identified by its highest weight in coordinates, {$\mu(x_1,\dots
,x_{[a/2]})$, or if there are repetitions,
$\mu(x_1^{n_1},\dots,x_k^{n_k})$. For $O(a)$ we use the
parametrization of Weyl, \cite{Weyl}. Embed $O(a)\subset U(a)$ in the
standard way. Then we denote by $\mu(x_1,\dots ,x_k,0^{[a/2]-k};\ep)$
the irreducible $O(a)-$component generated by the highest weight of
the representation 
$$
\mu(x_1,\dots,x_k,1^{(1-\ep)(a/2-k)},0^{a-k-(1-\ep)(a/2-k)})
$$ 
of $U(a)$. In these formulas,
$\ep=\pm,$ is often written as $+$ for $1,$ and $-$ for
$-1$.}

\subsection{}\label{sec:4.4}
We describe the \textbf{relevant} $K$-types for the orthogonal groups
$O(a,a).$ 
{
\begin{definition}[even orthogonal groups]
  \label{def:4.4}
The \textbf{relevant} $K-$types for $O(a,a)$ 
\begin{align}
&\mu_e([a/2]-r,r):=\mu(0^{[a/2]};+)\otimes \mu(2^r,0^l;+)
\label{eq:4.4.1}\\
&\mu_o(r,[a/2]-r):=\mu(1^r,0^l;+)\otimes\mu(1^r,0^l;+).\label{eq:4.4.2}
\end{align}
where $r+l=[a/2]$.
\end{definition}
}
\begin{proposition}
  \label{p:4.4}
The relevant $K-$types are quasispherical. The representation of
$W(D_a)$  of $O(a,a)$ on $V^M$ is 
 \begin{align}
&\sig[(r,a-r), (0)]&&\longleftrightarrow&& 
\mu(0^{[a/2]};+)\otimes \mu(2^r,0^l;+),
\label{eq:4.4.3}\\
&\sig[(a-k), (k)],&&\longleftrightarrow&&
\mu(1^k,0^l;+)\otimes\mu(1^k,0^l;+),\label{eq:4.4.4}
\end{align}
When $l=0,$ and $a$ is even,
\begin{align}
  \label{eq:4.4.5}
&\sig[(a/2,a/2),(0)]&&\longleftrightarrow
&&\mu(0^{a/2})\otimes\mu(2^{a/2-1},\pm 2), \\ 
&\sig[(a/2),(a/2)]_{I,II} &&\longleftrightarrow&&
\mu(1^{a/2-1},\pm 1)\otimes\mu(1^{a/2-1},\pm 1).  
\end{align}
\end{proposition}
\noindent We will prove this together with the corresponding proposition for
$O(a+1,a)$ in section \ref{sec:4.6}.

\subsection{}\label{sec:4.5}
We describe the relevant $K$-types for $O(a+1,a)$
\begin{definition}[odd orthogonal groups]
  \label{def:4.5}
The \textbf{relevant} $K$-types for  $O(a+1,a)$ are
\begin{align}
&\mu_e(a-r,r):=\mu(0^{[(a+1)/2]};+)\otimes\mu(2^r,0^l;+)
\label{eq:4.5.1}\\
&\mu_o(a-k,k):=\mu(1^k,0^l;+)\otimes\mu(1^k,0^s;+)
\label{eq:4.5.2}\\
&\mu_o(k,a-k):=\mu(1^{k+1},0^l;+)\otimes\mu(1^k,0^s;+)
\label{eq:4.5.3}
\end{align}
where $r+l=[a/2]$ in (\ref{eq:4.5.1}), $k+l=[(a+1)/2],\ k+s=[a/2]$ in
(\ref{eq:4.5.2}), and $k+1+l=[(a+1)/2],\ k+s=[a/2]$ in (\ref{eq:4.5.3}).
\end{definition}
{
\begin{proposition}
  \label{p:4.5}
The representations of $W(B_a)$ on $V^M$ for the relevant $K$-types are
\begin{align}
&\sig[(r,a-r),(0)]&&\longleftrightarrow&&
\mu(0^{[(a+1)/2]};+)\otimes\mu(2^r,0^l;+) \label{eq:4.5.4}\\
&\sig[(a-k),(k)]&&\longleftrightarrow&& 
\mu(1^k,0^{[(a+1)/2]-k};+)\otimes\mu(1^k,0^{[a/2]-l};+)
,\label{eq:4.5.5}\\
&\sig[(k),(a-k)]&&\longleftrightarrow&& 
\mu(1^{k+1},0^{[(a+1)/2]-k-1};+)\otimes\mu(1^k,0^{[a/2]-k};+).\label{eq:4.5.6}
\end{align}
When $a$ is even, 
\begin{align}
&\sig[(a/2),(a/2)]&&\longleftrightarrow&& 
\mu(1^{a/2})\otimes\mu(1^{a/2-1},\pm 1).  \label{eq:4.5.7}
\end{align}
When $a$ is odd,
\begin{align}
&\sig[(\frac{a-1}{2}),(\frac{a-1}{2})]&&\longleftrightarrow&& 
\mu(1^{(a-1)/2},\pm 1)\otimes\mu(1^{(a-1)/2}).  \label{eq:4.5.8}
\end{align}
\end{proposition}
}
The proof will be in section \ref{sec:4.6}.
\subsection{Proof of propositions \ref{p:4.4} and
  \ref{p:4.5}}\label{sec:4.6} 
We use the standard realization of the orthogonal groups $O(a+1,a)$
and $O(a,a).$  Let 
\begin{equation}
  \label{eq:4.6.1}
  \wti{M}:=\{(\eta_0,\eta_1,\dots ,\eta_a,\ep_1,\dots ,\ep_a)\ :\ 
\eta_i,\ \ep_j=\pm 1,\prod\eta_i=\prod\ep_j=1\},
\end{equation}
viewed as the subgroup of $O(a+1)\times O(a)$ with the $\eta_i,\
\ep_j$ on the diagonal. With the appropriate choice of $\fk a\cong
\bR^a,$ $\wti M\subset N_K(\fk a),$ and the action is
\begin{equation}
  \label{eq:4.6.2}
  (\eta_i,\ep_j)\cdot (\dots ,x_k,\dots )=(\dots ,\eta_k\ep_k
  x_k,\dots ).
\end{equation}
Then $M:=K\cap B$ is the subgroup of $\wti{M}$
determined by the relations $\eta_j=\ep_j,\ j=1,\dots ,a.$ Similarly
for $O(a)\times O(a)$ but there is no $\eta_0.$  

We do the case $O(a+1,a)$, $O(a,a)$ is similar. 
The representations $\mu_o(a-k,k)$ and $ \mu_o(k,a-k)$ can be realized as 
$\bigwedge^k\bC^{a+1}\otimes\bigwedge^k\bC^a,$ respectively
$\bigwedge^{k+1}\bC^{a+1}\otimes \bigwedge^k(\bC^a).$  
Let $e_i$ be a basis of $\bC^{a+1}$ and $f_j$ a basis of $\bC^a.$ 
The space $V^M$ is the span of the vectors $e_{i_1}\wedge\dots\wedge
e_{i_k}\otimes f_{i_1}\wedge\dots\wedge f_{i_k},$ and $e_0\wedge
e_{i_1}\wedge\dots\wedge e_{i_k}\otimes f_{i_1}\wedge\dots\wedge
f_{i_k}.$ The elements of $W$ corresponding to short root reflections
all have representatives of the form $\eta_0=-1,\eta_j=-1,$ the rest
zero. The action of $S_a\subset W$ on the space $V^M$ is by permuting the $e_i,
f_j$ diagonally. Claims (\ref{eq:4.4.4}-\ref{eq:4.4.5})
and (\ref{eq:4.5.5}-\ref{eq:4.5.6}) follow from these considerations,
we omit further details.

\medskip
For cases (\ref{eq:4.4.3}) and (\ref{eq:4.5.4}) we do an induction on
$r.$ We do the case $O(a,a)$ only.  
The claim is clear for $r=0.$ Since the first factor of
$\mu_e([a/2]-r,r)$ is the trivial representation, we only concern ourselves
with the second factor. Consider
$\bigwedge^r\bC^a\otimes\bigwedge^r\bC^a.$  The space of $M-$fixed
vectors has dimension $\binom{a}{r}$, and a basis is
\begin{equation}
  \label{eq:4.6.3}
  e_{i_1}\wedge\dots\wedge e_{i_r}\otimes e_{i_1}\wedge\dots\wedge e_{i_r}
\end{equation}
 As a module of $S_a,$ this is
\begin{equation}
  \label{eq:4.6.4}
  Ind_{S_r\times S_{a-r}}^{S_a}[triv\otimes triv]=\sum_{1\le j\le r} (j,a-j)
\end{equation}
On the other hand, the tensor product
$\bigwedge^r\bC^r\otimes\bigwedge^r\bC^a$ consists of representations with
highest weight $\mu(2^\al,1^\beta,0^\gamma).$ From the explicit
description of $\bigwedge^k\bC^a,$ and the action of $M,$ we can infer
that $V^M$ for $\beta\ne 0$ is $(0).$ This is because the
representation occurs in
$\bigwedge^{\al+\beta}\bC^a\otimes\bigwedge^{\al}\bC^a,$ which has no
$M-$fixed vectors. But
$\mu(2^j,0^l)$ for $j\le r$ occurs (for example by the P-R-V
conjecture).  By the induction hypothesis, $(j,a-j)$
occurs in $\mu(2^j,0^l),$ for $j<r,$ and so only $(r,a-r)$ is
unaccounted for. Thus $V^M$ for $\mu_e([a/2]-r,r)$ 
be $(r,a-r).$ The claim now follows from the fact that the action of
the short root reflections is trivial, and the description of the
irreducible representations of $W(B_a).$

\subsection{General linear groups}\label{sec:4.7} The maximal compact
subgroup of 
$GL(a,\bb R)$ is $O(a),$ the Weyl group is $W(A_{a-1})=S_a$ and
$M\cong\uset{a}{\ub{O(1)\times\dots\times O(1)}}.$ We list
the case of the connected component $GL(a,\bR)^+$ (matrices with
positive determinant) instead, because its maximal compact group is
$K=SO(a)$ which is connected, and irreducible representations are
parametrized by their highest weights. 
\begin{definition}
  \label{def:4.7}
The \textbf{relevant} $K-$types are the ones with highest weights 
\begin{equation*}
 \mu(2^k,0^l).
\end{equation*}
The corresponding Weyl group representations on $V^M$ are
$\sig[(k,a-k)].$
\end{definition}
We omit the details, the proof is essentially the discussion about the
representation of $S_a$ on $\bigwedge\bC^a\otimes\bigwedge\bC^a$ for the
orthogonal groups. 

\subsection{Relevant $W-$types}\label{sec:5.8}
\begin{definition} Let $W$ be the Weyl group of type B,C,D.
  \label{d:5.8}
The following $W-$types will be called \textbf{relevant}.
\begin{equation}
  \label{eq:4.8.1}
  \sig_e(n-r,r):=\sig[(n-r),(r)],\qquad \sig_o(k,n-k):=\sig[(k,n-k),(0)]
\end{equation}
In type D for $n$ even, and $r=n/2$ there are two $W-$ types,\newline
$\sig_e[(n/2),(n/2)]_{I,II}:=\sig[(n/2),(n/2)]_{I,II}.$ 
If the root system is not simple, the relevant $W-$types are tensor
products of relevant $W-$types on each factor.
\end{definition}

\section{Intertwining Operators}\label{sec:5} 

\medskip
\subsection{}\label{sec:5.1} 
{The notation from chapter \ref{sec:4} is in effect.} Recall that
$X(\nu)$ denotes the spherical principal series. Let $w\in W.$ 
Then there is an intertwining operator 
\begin{equation}
  \label{eq:5.1.1}
  I(w,\nu): X(\nu)\longrightarrow X(w\nu) 
\end{equation}
defined as follows. Recall $B=HN$ is a Borel subgroup, and let
$\ovr{\cdot}{w}$ be a representative of $w.$  Realize $X(\nu)$ on the
space of $K-$finite functions

$$
\{f:G\longrightarrow\bC\ :\
f(gb)=\chi(b^{-1})\delta_B(b)^{-1/2}f(g)\}\qquad \pi(g)f(x):=f(g^{-1}x).
$$ 
If $f\in X(\nu),$ then 
$$
I(w,\nu)f(g):=\int_{N/N\cap wNw^{-1}} f(g\ovr{\cdot}{w}^{-1}n)dn.
$$
{These intertwining operators make
sense in both the real and $p-$adic case. They have the properties that
the integral converges for $Re\langle \nu,\al\rangle >>0,$ and 
$I(w,\nu)$ has an analytic continuation for
$\langle Re\nu,\al\rangle \ge 0.$ 
}

\medskip
If $(\mu,V)$ is a $K-$type, ($K$ the maximal compact subgroup in the
real case, $G(\C R)$ in the $p-$adic case) then $I(w,\nu)$ induces a map 
\begin{equation}
  \label{eq:5.1.2}
  I_V(w,\nu): \Hom_K [V,X(\nu)]\longrightarrow \Hom_K[ V,X(w\nu)].
\end{equation}
By Frobenius reciprocity, we get a map
\begin{equation}
  \label{eq:5.1.3}
R_{V}(w,\nu): (V^*)^{K\cap B}\longrightarrow  (V^*)^{K\cap B}.
\end{equation}  

In case $(\mu,V)$ is trivial the spaces are 1-dimensional and $R_V(w,\nu)$ is
a scalar. We normalize $I(w,\nu)$ so that this scalar is 1. The
$R_V(w,\nu)$ are meromorphic functions in $\nu,$ and the 
$I(w,\nu)$ have the following additional properties.
\begin{enumerate}
\item If $w=s_{\al_1}\cdots s_{\al_k}$ is a reduced decomposition
  ($\al_i$ simple roots and $s_{\al_i}$ the corresponding root
  reflections) let $\ell(w):=k.$   If $w=w_1\cdot w_2$ with
  $\ell(w)=\ell(w_1)+\ell(w_2),$ then 
  $I(w,\nu)=I(w_1,w_2\nu)\circ I(w_2,\nu).$ In particular, 
$I(w,\nu)$ factors into a product of intertwining operators
$I_j$, one for each $s_{\al_j}.$ These operators are
\begin{equation}
  \label{eq:5.1.4}
 I_j\ :\ X(s_{\al_{j+1}}\dots s_{\al_k}\cdot\nu)\longrightarrow
X(s_{\al_j}\dots s_{\al_k}\cdot\nu)  
\end{equation}
\item Let $P$ be a standard parabolic subgroup with levi component
  $M$ so that $A\subset M$, and $w\in
  W(M,A).$ Because $X(\nu)=\Ind_P^G[X_M(\nu)],$  
we can induce the intertwining operator
  $I_M(w,\nu):X_M(\nu)\longrightarrow X_M(w\nu)$ to form the
  \textit{induced intertwining operator}
\[
Ind_M^G[I_M(w,\nu)]:X(\nu)=Ind_P^G[X_M(\nu)]\longrightarrow 
X(w\nu)=Ind_P^G[X_M(w\nu)],
\]
Then  $I(w,\nu)=Ind_M^G[I_M(w,\nu)].$
\item If $Re\langle \nu,\al\rangle \ge 0$ for all positive roots $\al,$
  then $R_V(w_0,\nu)$ has no poles, and the image of $I(w_0,\nu)$ ($w_0\in
  W$ is the long element) is   $L(\nu).$ 
\item If $-\ovl{\nu}$ is in the same Weyl group orbit as $\nu,$ let
  $w$ be the shortest element so that $w\nu=-\ovl{\nu}.$ Then the
  hermitian dual of $X(\nu)$ is $X(\nu)^h=X(w\nu).$ {Let $\langle\
  ,\ \rangle_h$ be the hermitian pairing between $X(\nu)$ and
  $X(w\nu).$} It follows that $L(\nu)$
  is hermitian with inner product 
\[
\langle v_1,v_2\rangle:=\langle v_1, I(w,\nu)v_2\rangle_h.
\]
\end{enumerate}
{
For the remainder of the section we consider the real case only.}

Let $\al$ be a simple root and $P_\al=M_\al N_\al$ be the standard
parabolic subgroup {so that the derived part of the Lie algebra
$M_\al$ is isomorphic to the $sl(2,\bb R)$ generated by the root
vectors $E_{\pm \al}.$} We assume that $\theta E_\al=-E_{-\al},$ {where
  $\theta$ is   the Cartan involution corresponding to $K$}. Let
$D_\al=\sqrt{-1}(E_\al-E_{-\al})$ and $s_\al=e^{\sqrt{-1}\pi
  D_\al/2}.$ {Here by $s_\al,$ we actually mean the representative
in $N_K(A)$ of the Weyl group reflection.} 
Then
$s_\al^2=m_\al$ is in $K\cap B\cap M_\al.$ Since the square of any element
in $K\cap B$ is in the center, and $K\cap B$ 
normalizes the the root vectors, $\Ad
m (D_\al) =\pm D_\al.$ Grade $V^*=\oplus V^*_i$ according to the
absolute values of the eigenvalues of $D_\al$ (which are
integers). Then $K\cap B$ preserves this grading and 
\[
(V^*)^{K\cap B}=\bigoplus_{i\ even} (V^*_i)^{K\cap B}.
\]
The map $\psi_\al: sl(2,\bb R)\longrightarrow \fk g$ determined by 
\[
\psi_\al\begin{bmatrix} 0&1\\0&0\end{bmatrix}=E_\al,\qquad 
\psi_\al\begin{bmatrix} 0&1\\0&0\end{bmatrix}=E_{-\al}
\]
determines a map
  \begin{equation}\label{eq:5.1.5}
    \Psi_\al\ :\ SL(2,\bb R)\longrightarrow G
  \end{equation}
with image $G_\al,$ a connected group with Lie algebra isomorphic to
$sl(2,\bb R).$ 
\begin{proposition}\label{p:5.1} On $(V^*_{2m})^{K\cap B},$
$$
R_{V}(s_\al,\nu)=\begin{cases} Id &\text{ if } m=0,\\
\prod_{0\le j< m} \frac{2j+1-<\nu,\cha>}{2j+1+<\nu,\cha>}\ Id\quad 
                &\text{ if } m\ne 0.
\end{cases}
$$
In particular, $I(w,\nu)$ is an isomorphism unless $\langle
\nu,\cha\rangle\in 2\bZ +1.$  
\end{proposition}
\begin{proof}
The formula is well known for $SL(2,\bb R).$ The second assertion
follows from this and  the listed properties of intertwining operators.
\end{proof}

\begin{corollary}\label{c:5.1} 
For relevant $K-$types the formula is
$$
R_{V}(s_\al,\nu)=\begin{cases} Id &\text{ on the +1 eigenspace of } s_\al,\\ 
             \frac{1-<\nu,\cha>}{1+<\nu,\cha>}\ Id
                &\text{ on the -1 eigenspace of } s_\al.
\end{cases}
$$
When restricted to $(V^*)^{K\cap B},$ the long intertwining operator is the
product of the $R_V(s_\al,*)$ corresponding to the reduced decomposition of
$w_0$ and depends only on the Weyl group structure of $(V^*)^{K\cap B}.$
\end{corollary}
\begin{proof}
Relevant $K-$types have the property that the
{even}
eigenvalues of $D_\al$ are $0, \pm 2$ only. The element $s_\al$ acts
by $1$ on the zero eigenspace of $D_\al$ and by $-1$ on the $\pm 2$
eigenspace. The claim follows from this.  
\end{proof}
\subsection{}\label{sec:5.2}
We now show that the formulas in the previous section  coincide with
corresponding ones in the $p-$adic case. In the split $p$-adic case,
spherical representations are a subset of representations with $\C
I$-fixed vectors, where $\C I$ is an Iwahori subgroup. 
As explained in \cite{Borel}, the category of representations with $\C I$
fixed vectors is equivalent to the category of finite dimensional
representations of the Iwahori-Hecke algebra $\C H:=\C H(\CI\backslash G/\C
I)$. The equivalence is
\begin{equation}
  \label{eq:5.2.1}
  \C V\longrightarrow \C V^{\C I}.
\end{equation}
The papers \cite{BM1} and \cite{BM2} show that the
problem of the  determination of  the unitary dual of representations
with $\C I$ fixed vectors, is equivalent to the problem of the
determination of the unitary irreducible representations of $\C H$ 
with real infinitesimal character. In fact it is the affine graded
Hecke algebras we will need to consider, and they are as follows.

{
Let $\bA:=S(\va),$ and define the affine graded Hecke algebra to be
$\bH:=\bC[W]\otimes \bA$ as a vector space, and usual algebra
structure for  $\bC[W]$ and $\bA$. 
}

{The generators of
  $\bC[W]$ are denoted by $t_\al$ corresponding to the simple
  reflections $s_\al,$ while the generators of $\bA$ are
  $\omega\in\va.$} Impose the additional relation 
\begin{equation}\label{eq:5.2.2}
\om t_\al=s_\al(\om)t_\al + <\om,\cha>,\qquad \om\in \va,
\end{equation}
where $t_\al$ is the element in $\bC[W]$ corresponding to the simple
root $\al.$ If $X(\chi)$ is the standard (principal series) module
determined by $\chi,$  then 
\begin{equation}
  \label{eq:5.2.3}
  X(\chi)^\C I=\bH\otimes_\bA \bC_\chi.
\end{equation}
The intertwining operator  $I(w,\chi)$ is a product of operators
$I_{\al_i}$ according to a reduced decomposition of
$w=s_{\al_1}\cdot\dots\cdot s_{\al_k}.$ If $\al$ is a simple root, 
\begin{equation}\label{eq:5.2.4}
r_{\al}:=(t_\al\cha - 1)\frac{1}{\cha - 1},\qquad I_\al\ :\
x\otimes \one_\chi \mapsto xr_{\al}\otimes \one_{s_{\al}\chi}.
\end{equation}

The $I(w,\chi)$ have the same properties as in the real case. The
$r_\al$ are multiplied on the right, so we can replace $\cha$ with 
$-\langle\nu,\cha\rangle$ in the formulas. Furthermore, 
\[
\bb C[W]=\sum_{\sig\in \widehat W} V_\sig\otimes V_\sig^*.
\]
Since $r_\al$ acts as multiplication on the right, it gives rise to an operator 
\[
r_\sig(s_\al,\nu):V^*_\sig\longrightarrow V^*_\sig.
\]
\begin{theorem}\label{t:5.2}
The $R_{V}(s_\al,\nu)$  for the real case on relevant $K-$types coincide with
the $r_\sig(s_\al,\nu)$ on the  $V^*_\sig\cong (V^*)^{K\cap B}$ 
\end{theorem}
\begin{proof}
These operators act the same way:
\begin{equation}
  \label{eq:5.2.5}
  r_{\sig}(s_\al,\nu)=
  \begin{cases}
    Id &\text{ on the } +1 \text{ eigenspace of } s_\al,\\
    \frac{1-\langle \nu,\cha\rangle}{1+ \langle \nu,\al\rangle}Id
       &\text{ on the } -1 \text{ eigenspace of } s_\al.
  \end{cases}
\end{equation}
The assertion is now clear from corollary (\ref{c:5.1})
and formula (\ref{eq:5.2.2}).
\end{proof} 
\subsection*{Remarks}

  \begin{enumerate}
  \item The Hecke algebra for a $p$-adic group $G$ is typically defined using
the dual root system of the complex group $\vG$. For example the
formulas for $r_\al$  in the literature, \eg \cite{BM1}-\cite{BM3},
have  roots instead of coroots.  
\item The intertwining operator $I_\al$ is coincides with the
  intertwining operators obtained by grading the one coming from the
  p-adic group $SL(2)$ applied to $\C I-$fixed vectors. For more
  details on the properties of the intertwining operators in the
  affine graded Hecke algebra case, the reader
  may consult \cite{BM3} and \cite{BC2} section 2.   
  \end{enumerate}


\subsection{}\label{sec:5.3} The main point of section \ref{sec:5.2}
is that  for the real case, and a relevant $K$-type $(V,\mu)$, the intertwining
operator calculations coincide with the intertwining operator
calculations for the affine graded Hecke algebra on the space $V^{K\cap
B}.$
Thus we will deal with the Hecke algebra calculations exclusively, but
the conclusions hold for both the real and $p$-adic case. 
Recall from section \ref{sec:2.3} that to each $\chi$ we have
associated a nilpotent orbit $\vO$, and Levi components $\vm_{BC}$ and
$\vm_{KL}$. These are special instances of the following situation. 
Assume that $\vO$ is written as in (\ref{eq:2.3.4}) (\ie
$((a_1,a_1),\dots ,(a_k,a_k);(d_i))$ with
\begin{description}
\item[$\vg$ of type B] $(d_i)$ all odd; they are relabelled
  $(2x_0+1,\dots ,2x_{2m}+1)$, 
\item[$\vg$ of type C] $(d_i)$ all even; they are relabelled
  $(2x_0,\dots ,2x_{2m})$,
\item[$\vg$ of type D] $(d_i)$ all odd; they are relabelled
  $(2x_0+1,\dots , 2x_{2m-1}+1).$
\end{description}
Similar to (\ref{eq:2.3.5}), let 
\begin{equation}
  \label{eq:5.3.1}
  \vm :=gl(a_1)\times \dots\times gl(a_k)\times \vg(n_0),\quad
  n_0=n-\sum a_i.
\end{equation}
We consider parameters of the form $\chi=\vh/2+\nu.$ 

Write $\chi_0$ for the parameter $\vh/2,$ and
$\chi_i:=(-\frac{a_i-1}{2}+\nu_i,\dots ,\frac{a_i-1}{2}+\nu_i).$
 We
focus on $\chi_0$ as a parameter on  $\vg(n_0)$. We attach two Levi
components   
\begin{equation}\label{eq:5.3.2}
\begin{array}{ll}
{\vg_e:} &\\
B & gl(x_{2m-1}+x_{2m-2}+1)\times\dots\times gl(x_1+x_0+1)\times \vg(x_{2m})\\
C & gl(x_{2m-1}+x_{2m-2})\times\dots\times gl(x_1+x_0)\times\vg(x_{2m})\\
D &gl(x_{2m-1}+x_{2m-2}+1)\times\dots\times gl(x_1+x_0+1)\\
&\\
{\vg_o:} &\\
B & gl(x_{2m}+x_{2m-1}+1)\times\dots\times gl(x_2+x_1+1)\times \vg(x_0)\\
C & gl(x_{2m}+x_{2m-1})\times\dots\times gl(x_2+x_1)\times \vg(x_0)\\
D & gl(x_{2m-3}+x_{2m-4}+1)\times\dots\times
gl(x_{2m-2})\times \vg(x_{2m-1}+1).
\end{array}
\end{equation}
There are 1-dimensional representations $L(\chi_e)$ and
$L(\chi_o)$ such that the spherical irreducible representation
{ $L(\chi_0)=\ovl X(\chi_0)$ with infinitesimal character $\chi_0$}
is the spherical irreducible subquotient of $X_e:=Ind_{P_e}^G(L(\chi_e))$  and
$X_o:=Ind_{P_o}^G(L(\chi_o))$ respectively}

The parameters $\chi_e$ and $\chi_o$ are written in terms of strings as follows:
\begin{description}
\item[$X_e$]
\begin{equation}
  \label{eq:5.3.3}
  \begin{aligned}
 &B:\ \dots (-x_{2i-1},\dots ,
 x_{2i-2})\dots (-x_{2m},\dots ,-1)    \\
&C:\ \dots (-x_{2i-1} +1/2,\dots ,
 x_{2i-2}-1/2)\dots (-x_{2m}+1/2,\dots ,-1/2)    \\
&D:\ \dots (-x_{2i-1},\dots , x_{2i-2})\dots     \\
  \end{aligned}
 \end{equation}
 \item[$X_o$]
\begin{equation}
\begin{aligned}
&B:\   \dots (-x_{2i},\dots , x_{2i-1})\dots  (-x_{0},\dots , -1)  \\
&C:\ \dots(-x_{2i}+1/2,\dots , x_{2i-1}-1/2)\dots  (-x_{0}+1/2,\dots , -1/2)
\\
&D:\ \dots (-x_{2i},\dots , x_{2i-1})\dots (-x_{2m-2},\dots ,-1)
(-x_{2m-1}+1,\dots , 0)\\ 
\end{aligned}
\end{equation}
\end{description}
\begin{theorem}
\label{th:5.3}  
For the Hecke algebra, $p$-adic groups,
\begin{equation*}
\begin{aligned}
&[\sig[(n-r),(r)]\ :\ X_e]=[\sig[(n-r),(r)]\ :\ L(\chi_0)],\\
&[\sig[(k,n-k),(0)]\ :\ X_o]=[\sig[(k,n-k),(0)]\ :\ L(\chi_0)] 
\end{aligned}
\end{equation*}
hold.
\end{theorem}
\noindent The proof is in section {\ref{sec:6.8}.} 

\medskip
For a general parameter {$\chi=\chi_0+\nu$}, the strings defined in
section \ref{sec:2} and the above construction define parabolic
subgroups with Levi components $ gl(a_1)\times\dots\times
gl(a_k)\times \vg_e$ and $gl(a_1)\times\dots\times gl(a_r)\times \vg_o,$ and
corresponding $L_e(\chi)$ and $L_o(\chi)$. 
We denote these induced modules by $X_e$ and $X_o$ as well. 

\begin{corollary}
The relations
\begin{equation*}
\begin{aligned}
&[\sig[(n-r),(r)]\ :\ X_e]=[\sig[(n-r),(r)]\ :\ L(\chi)],\\
&[\sig[(k,n-k),(0)]\ :\ X_o]=[\sig[(k,n-k),(0)]\ :\ L(\chi)] 
\end{aligned}
\end{equation*}
hold in general. For real groups, in the notation of sections
\ref{sec:4.2}-\ref{sec:4.4}, 
\begin{equation*}
  \begin{aligned}
&[\mu_e(r,n-r)\ :\ X_e]=[\mu_e(r,n-r)\ :\ L(\chi)],\\
&[\mu_o(k,n-k)\ :\ X_o]=[\mu_o(k,n-k)\ :\ L(\chi)].   
  \end{aligned}
\end{equation*}
\end{corollary}

{
\begin{proof}

The results in section \ref{sec:5.2} show that the intertwining
operators on $\sig_e(k,n-k)$ for the $p-$adic group equal the
intertwining operators for $\mu_e(k,n-k)$ for the real group, and
similarly for $\sig_o$ and $\mu_o.$ Thus the multiplicities of the
$\sig_e/\sig_o$ in $L(\chi)$ for the $p-$adic case equal the
multiplicities of the corresponding $\mu_e/\mu_o$ in $L(\chi)$ in the real case.
We do the $p-$adic case first. Recall theorem \ref{t:2.6} which states
that $I_{M_{KL}}(\chi)=L(\chi).$ The Levi 
subgroup $M_{KL}$ is a product of $GL$ factors, which we will denote
$M_A$, with a factor $G(n_0).$ So for $W-$type multiplicities 
we can replace $I_{M_{KL}}(\chi)$ by $Ind_{M_A\times
  G(n_0)}[\otimes triv\otimes L(\chi_0)].$ 
We explain the case of $\sigma_e=\sigma_e(k,n-k)$,
that of $\sigma_o$ being identical. By Frobenius reciprocity, 
\begin{align}\label{eq:5.3.4}
\Hom_W[\sigma_e:L(\chi)]&=\Hom_W [\sigma_e:I_{M_{KL}}(\chi)]\\
                       &=\Hom_{W(M_{A})\times W(G(n_0))} [\sigma_e:
triv\otimes L(\chi_0)].
\end{align}
Using the formulas for restrictions of representations for Weyl groups
of classical types, it follows that
\begin{align}\label{eq:5.3.5}
\dim\Hom_W[\sigma_e:L(\chi)]&=\sum_{k'}
\dim\Hom_{W(G(n_0))}[\sigma_e(k',n_0-k'):L(\chi_0)]\notag\\ 
&\\
                           &=\sum_{k'}
\dim\Hom_{W(G(n_0))}[\sigma_e(k',n_0-k'):X_{e,G(n_0)}], \notag
\end{align}
where the last step is theorem \ref{sec:5.3}. Since
$X_e=Ind_{M_A\times G(n_0)}^G[\otimes L(\chi_i)\otimes
X_{e,G(n_0)}(\chi_0)]$, again by Frobenius reciprocity, one can show
that this is also equal to $\dim\Hom_W[\sigma_e:X_e]$. This proves the
claim for the $p$-adic case.

 In the real case the
proof is complete once we observe that in all the steps for the
$p-$adic case, the multiplicity of
$triv\otimes \mu_e(k',n_0-k')$ in the restriction of $\mu_e(k,n-k)$
matches the multiplicity of  $triv\otimes \sig_e(k',n_0-k')$ in the
restriction of $\sig_e(k,n-k)$. Similarly for $\mu_o$ and $\sig_o$. 

\end{proof}
}

\section{Hecke algebra calculations}\label{sec:6} 

\medskip
\subsection{}\label{sec:6.1} The proof of the results in
\ref{sec:5.3} is by a computation of intertwining operators on the
relevant $K-$types. It only depends on the $W-$type of $V^{K\cap B}$,
so we work in 
the setting of the Hecke algebra. The fact that we can deal
exclusively with $W-$types, is a big advantage. In particular we do
not have to worry about disconnectedness of Levi components.
We will write $GL(k)$ for the Hecke algebra of type $A$ and
$G(n)$ for the types $B,\ C$ or $D$ as the case may be. This is
so as to emphasize that the results are about groups, real or $p-$adic. 

The intertwining operators will be decomposed into products of simpler
operators induced from operators coming from maximal Levi subgroups. 
We introduce these first.

\medskip
Suppose $M$ is a Levi component of the form 
\begin{equation}
  \label{eq:6.1.1}
  GL(a_1)\times\dots\times GL(a_l)\times G(n_0).
\end{equation}
Let $\chi_i$ be characters for $GL(a_i).$ We simplify the notation
somewhat by writing
\begin{equation}
  \label{eq:6.1.2}
\chi_i\longleftrightarrow (\nu_i):= (-\frac{a_i-1}{2}+\nu_i,\dots ,
\frac{a_i-1}{2} +\nu_i). 
\end{equation}
The parameter is antidominant, and so $L(\chi_i)$ occurs as a
submodule of the principal series $X((\nu_i)).$ The module is spherical
1-dimensional, and the action of  $\mathfrak a$    is
\begin{equation}
  \label{eq:6.1.3}
  \chi_i(\omega)=\langle \omega , (\frac{a_i-1}{2}+\nu_i,\dots ,
  -\frac{a_i-1}{2} +\nu_i)\rangle,\qquad \omega\in\fk a,
\end{equation}
while $W$ acts trivially.
The trivial representation $\chi_0$ of $G(n_0)$ corresponds to the string
$(-n_0+\ep,\dots ,-1+\ep)$ where
\begin{equation}\label{eq:6.1.4}
\ep:=\begin{cases}
  0\qquad &\bH \text{ of type    B},\\
  1/2,\qquad &\bH \text{ of type    C},\\
  1,\qquad        &\bH \text{ of type   D}.
\end{cases}
\end{equation}
We abbreviate this as $(\nu_0)$.
Again $L(\chi_0)$ is the trivial representation, and
because $\chi_0$ is antidominant, it
appears as a submodule of the principal series $X(\chi_0).$
We abbreviate 
\begin{equation}
  \label{eq:6.1.5}
  X_M(\dots(\nu_i)\dots ):= Ind_{\prod GL(a_i)\times G(n_0)}^G[\otimes
  \chi_i\otimes triv].
\end{equation}
The module $X_M(\dots (\nu_i)\dots )$ is a submodule of the standard
module $X(\chi)$ with parameter corresponding to the strings
\begin{equation}
  \label{eq:6.1.6}
\chi:=(\dots ,-\frac{a_i-1}{2}+\nu_i,\dots,\frac{a_i-1}{2}+\nu_i,\dots
  ,-n_0 +\ep,\dots ,-1 +\ep).
\end{equation}
In the setting of the Hecke algebra, the induced modules
(\ref{eq:6.1.5}) is really 
$X_M(\dots(\nu_i)\dots)=\bH\otimes_{\bH_M}[\bigotimes\chi_i\otimes triv].$  

Let $w_{i,i+1}\in W$ be the shortest Weyl group element which interchanges the 
strings $(\nu_i)$ and $(\nu_{i+1})$  in $\nu,$ and
fixes all other coordinates. The intertwining operator $I_{w_{i,i+1}}:
X(\nu)\longrightarrow X(w_{i,i+1}\nu)$ restricts to an intertwining operator
\begin{equation}
  \label{eq:6.1.7}
  \begin{aligned}
  I_{M,i,i+1}&(\dots(\nu_i)(\nu_{i+1})\dots ):\\
&X_M(\dots (\nu_i)(\nu_{i+1})\dots)\longrightarrow
X_{w_{i,i+1}M}(\dots (\nu_{i+1})(\nu_{i})\dots).    
  \end{aligned}
\end{equation}
This operator is induced from the same kind for $GL(a_i+a_{i+1})$
where $M_{i,i+1}=GL(a_i)\times GL(a_{i+1})\subset GL(a_i+a_{i+1})$ is the Levi
component of a maximal parabolic subgroup. 

\medskip
Let  $w_l\in W$ be  the
 shortest element which changes {$(\nu_l)$ to $(-\nu_l),$} 
and fixes all  other coordinates. It induces an intertwining operator
\begin{equation}
  \label{eq:6.1.8}
  I_{M,l}(\dots(\nu_l)(\nu_0)):X_M(\dots (\nu_l),(\nu_0))\longrightarrow
X_{w_lM}(\dots (-\nu_l),(\nu_0)).
\end{equation}
In this case, $w_lM=M$, so we will not always include it in the notation. 
In type D, if $n_0=0,$ 
the last entry of the resulting string might have to stay 
$-\frac{a_l-1}{2}+\nu_l $ instead of $\frac{a_l-1}{2}-\nu_l.$
This operator is induced from the same kind on $G(a_l+n_0)$ with
$M_l=GL(a_l)\times G(n_0)\subset G(a_l+n_0)$ the Levi component of a
maximal parabolic subgroup. 
\begin{lemma}
The operators $I_{M,i,i+1}$ and $I_{M,l}$ are meromorphic in $\nu_i$
in both the real and $p$-adic case.
\begin{enumerate}
\item $I_{M,i,i+1}$ has poles only 
if {$\frac{a_i-1}{2}+\nu_i-\frac{a_{i+1}-1}{2}-\nu_{i+1}\in\bZ.$}
If so,  a pole only occurs  if
\begin{equation*}
-\frac{a_i-1}{2}+\nu_i<-\frac{a_{i+1}-1}{2}-\nu_{i+1},\qquad 
\frac{a_i-1}{2}+\nu_i<\frac{a_{i+1}-1}{2}+\nu_{i+1}.  
\end{equation*}
\item $I_{M,l}$ has a pole only if $\frac{a_l-1}{2}+\nu_l\equiv \ep
  (mod\  \bZ)$. In that case,  a pole only occurs if
  \begin{equation*}
 -\frac{a_l-1}{2}+\nu_l<0.   
  \end{equation*}
\end{enumerate} 
\end{lemma}
\begin{proof}
We prove the assertion for $I_{M,i,i+1}$,  the other one is similar.  
The fact that the integrality condition is necessary is clear. For the second
condition, it is sufficient to consider the case $M=GL(a_1)\times
GL(a_2)\subset GL(a_1+a_2).$   If the strings are strongly nested,
then the operator cannot have any pole because $X_M$ is
irreducible. Remains  to show there is no pole in  the case when
$-\frac{a_{2}-1}{2}+\nu_{2} \le -\frac{a_1-1}{2}+\nu_1,$ and
  $\frac{a_1-1}{2}+\nu_1>\frac{a_{2}-1}{2}+\nu_{2}.$ Let 
    \begin{equation}
      \label{eq:6.1.9}
      \begin{aligned}
&M':=GL(\frac{a_1+a_2}{2}+\nu_2+\nu_1)\times
GL(\frac{a_1-a_2}{2}+\nu_1-\nu_2)\times GL(a_2),\\      
&(\nu'_1)=(-\frac{a_1-1}{2}+\nu_1,\dots ,\frac{a_2-1}{2}+\nu_2)\\
&(\nu'_2)=(\frac{a_2-1}{2}+1+\nu_2,\dots ,\frac{a_1-1}{2}+\nu_1)\\
&(\nu_3')=(\nu_2)=(-\frac{a_2-1}{2}+\nu_2,\dots ,\frac{a_2-1}{2}+\nu_2). 
      \end{aligned}
    \end{equation}
Then $X_M((\nu_1)(\nu_2))\subset X_{M'}((\nu_1')(\nu_2')(\nu_3')),$
and $I_{M,1,2}$ is the restriction of
$I_{w_{2,3}M',1,2},((\nu_1')(\nu_3')(\nu_2')\circ
I_{M',2,3}((\nu'_1)(\nu'_2)(\nu'_3))$ to $X_M.$ Because the strings
$(\nu_1')(\nu_3')$  are strongly nested, $I_{w_{2,3}M',1,2}$ has no pole, and
$I_{M',2,3}$ has no pole because it is a restriction of operators
coming from $SL(2)$'s which do not have poles. The claim follows.
\end{proof}
Let $\sig$ be a $W-$type. We are interested in computing 
$r_\sig(w,\dots(\nu_i)\dots)$, where $w$ changes all the $\nu_i$
for $1\le i$ to $-\nu_i.$   The operator  can be
factored into a product of $r_\sig(w_{i,i+1},*)$ of the type (\ref{eq:6.1.7})
and $r_\sig(w_l,*)$ of the type (\ref{eq:6.1.8}). These operators are
more tractable. Here's a more precise explanation. Let $M$ be the Levi
component   
\begin{align} 
&GL(a_1)\times\dots\times GL(a_i +a_{i+1})\times\dots
\text{ in case (\ref{eq:6.1.7})}\label{eq:6.1.11}\\
&GL(a_1)\times\dots\times G(a_l+n_0)\text{ in case (\ref{eq:6.1.8})}\label{eq:6.1.12}  
\end{align}
Since $X_M$ is induced from the trivial $W(M)$ module,
\begin{align}
\Hom_{W}&[\sig, X_M((\nu_i))]=Hom_{W(M)}[\sig|_{W(M)}\ :\ triv\otimes
X_{M_{i,i+1}}((\nu_i),(\nu_{i+1}))\otimes triv]\notag \\
&\text{ in case (\ref{eq:6.1.7}) }\label{eq:6.1.13}\\
\Hom_{W}&[\sig, X((\nu_i))]=Hom_{W(M)}[\sig|_{W(M)}\ :\ triv\otimes
X_{M_l}((\nu_l),(\nu_0))]\notag\\
 &\text{ in case (\ref{eq:6.1.8}) }\label{eq:6.1.14} 
\end{align}
where $M_{{i,i+1}}=GL(a_i)\times GL(a_{i+1})$ is a maximal Levi
component of $GL(a_i+a_{i+1})$ and {$M_l=GL(a_l)\times G(n_0)$} 
is a maximal Levi component of $G(a_l+n_0).$
To compute the $r_{\sig}(w_{i,i+1},*)$ and $r_{\sig}(w_l,*)$, it is
enough to compute the corresponding $r_{\sig_j}$ for the $\sig_j$
ocuring in the restriction $\sig\mid_{W(M)}$ in the cases $GL(a_i)\times
GL(a_{i+1})\subset GL(a_i+a_{i+1})$ and $GL(a_l)\times G(n_0)\subset
G(a_l+n_0)$. 
The restrictions of relevant $W-$types to Levi components consist of
relevant $W-$types of the same kind, \ie $\sig[(n-r),(r)]$ restricts to
a  sum of representations of the kind $\sig_e,$ and $\sig[(k,n-k),(0)]$
restricts to a sum of $\sig_o$. Typically the multiplicities of the
factors are 1. 

We also note that 
\begin{equation}
  \label{eq:6.1.15}
  X_M\mid_W=\sum_{\sig\in\wht W} V_\sig\otimes (V_\sig^*)^{W(M)}.
\end{equation}
So the  $r_\sig(w,*)$ map $(V_\sig^*)^{W(M)}$ to $(V_\sig^*)^{W(wM)}.$

\subsection{Maximal Levi components, summary of results}\

{
In the next sections we will compute the cases of
Levi components of maximal parabolic subgroups. We summarize the
results. A typical Levi component will be denoted $GL(k)\times
G(n)\subset G(k+n).$ The type refers to the type of $G(n),$ 

\begin{theorem}\ Let $\ep$ be as in \ref{eq:6.1.4}.

  \begin{description}
  \item[Type A] The interwtining operator $I_{M,1,2}((\nu_1)(\nu_2))$
    restricted to $\sig(m,k+n-m)$ is multiplication by the scalar
    \begin{equation}
      \label{eq:mpa}
 r_{\sig(m,k+n-m)}((\nu_1),(\nu_2))=\prod_{0\le j\le  m-1}
\frac{(\nu_1-\frac{k-1}{2})-(\frac{n-1}{2}+\nu_2+1)+j}
{(\nu_1+\frac{k-1}{2})-(-\frac{n-1}{2}+\nu_2-1)-j}.     
    \end{equation}
\item[Type B] For $\sig=\sig_e(m,k+n-m)$, 
  \begin{equation}
    \label{eq:mteb}
 r_{\sig_e(m,k+n-m)}((\nu))= \und{0\le j\le m-1}{\prod}
\frac{n+1/2-(-\frac{k-1}{2}+\nu)-j}{n+1/2+(\frac{k-1}{2}+\nu)-j} 
  \end{equation}
\item[Type C] For $\sig=\sig_e(m,k+n-m)$,
\begin{equation}\label{eq:mtec}
r_{\sig_e(m,k+n-m)}((\nu))=\und{0\le j\le m-1}{\prod}
\frac{n+1/2-(-\frac{k-1}{2}+\nu)-j}{n+1/2+(\frac{k-1}{2}+\nu)-j}
\end{equation}
\item[Type D] For $\sig=\sig_e(m,k+n-m)$,
\begin{equation}\label{eq:mted}
r_{\sig_e(m,k+n-m)}((\nu))=\und{0\le j\le m-1}{\prod}
\frac{n-(-\frac{k-1}{2}+\nu)-j}{n+(\frac{k-1}{2}+\nu)-j}
\end{equation}
  \end{description}
For types B,C the formulas hold for $n=0$ as well. 
In type D when $n=0,$ there are two cases $GL(k)\subset G(k)$ and $GL(k)'\subset
G(n)$. The formula is the same for $r_{\sig_e(n-m,m)}((\nu))$ and
$r_{\sig_e(n-m,m)}((\nu)')$:
\begin{equation}
  \label{eq:mted0} 
\prod_{0\le j<m}\frac{(\frac{k-1}{2}-\nu)-j}{(\frac{k-1}{2}+\nu)-j} .
\end{equation}
For $2m=n$ there are two representations with subscript I and II; the
formula above is  the same.

\medskip
For $\sig=\sig_o(m,k+n-m)$ and type B,C and D, the scalar\newline
$r_{\sig_o(m,k+n-m)}((\nu)(\nu_0))$ equals
\begin{equation}
  \label{eq:mto}
\prod_{0\le j\le
  m-1}\frac{(\nu-\frac{k-1}{2})-(1-\ep)+j}{(\nu+\frac{k-1}{2})
              -(-n -\ep)-j}\cdot
\frac{(-n-\ep )-(-\nu+\frac{k-1}{2})+j}{(1-\ep )-(-\nu-\frac{k-1}{2})-j}
\end{equation}
In this case $n>0$ for $\sig_o$ to occur in the induced module.
\end{theorem}
}
\subsection{$GL(k)\times GL(n)\subset GL(k+n)$}\label{6.2} 
This is the case of $I_{i,i+1}$ with $i<l.$ 
The module $X_M((\nu_1),(\nu_2))$
induced from the characters corresponding to 
\begin{equation}\label{eq:6.2.1}
  (-\frac{k-1}{2}+\nu_1,\dots
  ,\frac{k-1}{2}+\nu_1),(-\frac{n-1}{2}+\nu_2,\dots, \frac{n-1}{2}+\nu_2)
\end{equation}
has the following $S_{k+n}$ structure. Let $s:=\min(k,n)$ and write
$\sig(m,k+n-m)$ for the module corresponding to the partition
$(m,k+n-m),\ 0\le m\le s.$ Then 
\begin{equation}
  \label{eq:6.2.2}
  X_M((\nu_1),(\nu_2))\mid_{W}=\bigoplus_{0\le m\le s} \sig(m,k+n-m).
\end{equation}
\begin{lemma}\label{l:6.2} For $1\le m\le s,$ the intertwining operator 
$I_{M,1,2}((\nu_1)(\nu_2))$ restricted to $\sig$ gives 
\begin{equation*}
r_{\sig(m,k+n-m)}(\nu_1,\nu_2)=\prod_{0\le j\le  m-1}
\frac{(\nu_1-\frac{k-1}{2})-(\frac{n-1}{2}+\nu_2+1)+j}
{(\nu_1+\frac{k-1}{2})-(-\frac{n-1}{2}+\nu_2-1)-j}. 
\end{equation*}
\end{lemma}
\begin{proof}
The proof is an induction on $k,\ n$ and $m.$ We omit most details but give
the general idea. Assume $0<m<s,$ the case $m=s$ is simpler. 
Embed $X_M((\nu_1),(\nu_2))$ into $X_{M'}((\nu'),(\nu''),(\nu_2))$,
{where $M'=GL(k-1)\times GL(1)\times GL(n)$,}
corresponding to the strings 
\begin{equation}
  \label{eq:6.2.3}
  (-\frac{k-1}{2}+\nu_1,\dots,\frac{k-3}{2}+\nu_1),
  (\frac{k-1}{2}+\nu_1),(-\frac{n-1}{2}+\nu_2,\dots, \frac{n-1}{2}+\nu_2).
\end{equation}
The intertwining operator $I_{M,1,2}(\nu_1,\nu_2)$ is the restriction 
of 
\begin{equation}
  \label{eq:6.2.4}
I_{M',1,2}(\nu',\nu_2,\nu'')\circ I_{M',2,3}(\nu';\nu'',\nu_2)
\end{equation}
to $X_M((\nu_1),(\nu_2))\subset X_{M'}((\nu'),(\nu''),(\nu_2))$.
By an  induction on $k+n$ we can assume that these operators are known.
The $W-$type $\sig(m,k+n-m)$ occurs with multiplicity 1 in
$X_M((\nu_1),(\nu_2))$ and with multiplicity  2 in
$X_{M'}((\nu'),(\nu ''),(\nu_2)).$ The  restrictions are
\begin{align}
\sig(m,k+n-m)\mid_{W(M')}&= 
triv\otimes\sig(m-1,n+1-m) + triv\otimes\sig(m,n-m)\\
&\text{ for } I_{M',1,2}\label{eq:6.2.5}\\
\sig(m,k+n-m)\mid_{W(M')}&=\sig(1,n) + \sig(0,n+1) \text{ for } I_{M',2,3}
\label{eq:6.2.6}
\end{align}
The representation $\sig(m,k+n-m)$ has a realization as harmonic
polynomials in $S(\fk a )$ spanned by
\begin{equation}
  \label{eq:6.2.7}
  \prod_{1\le \ell\le m} (\ep_{i_\ell}-\ep_{j_\ell})
\end{equation}
{where $(i_1,j_1),\dots ,(i_k,j_k)$ are $m$ pairs of integers
satisfying $1\le i_\ell,\ j_\ell\le k+n,$ and $i_\ell\ne j_\ell.$} 
We apply the intertwining operator to the $S_k\times S_n-$fixed vector
\begin{equation}
  \label{eq:6.2.8}
e:=\sum_{x\in S_k\times S_n} 
x\cdot[(\ep_1 -\ep_{k+1})\times\dots\times (\ep_m-\ep_{k+m})].
\end{equation}
The intertwining operator $I_{M',2,3},$ has a simple form on the vectors
\begin{align}
&e_1:=\sum_{x\in S_{k-1}\times S_{n+1}}
x\cdot[(\ep_1 -\ep_{k+1})\times\dots\times (\ep_m-\ep_{k+m})],
\text{ in } \sig(0,n+1)\label{eq:6.2.9}\\
&e_2:=\sum_{x\in S_{k-1}\times S_1\times S_{n}}
x\cdot[(\ep_1 -\ep_{k+1})\times\dots\times (\ep_{m-1}-\ep_{k+m-1})
(\ep_k-\ep_{k+m})],\text{ in } \sig(1,n)\label{eq:6.2.10}
\end{align}
which appear in (\ref{eq:6.2.6}). They are mapped into
scalar multiples (given by the lemma) of the vectors $e_1',\ e'_2$
which are invariant under $S_{k-1}\times S_n\times S_1,$ and transform
according to $triv\otimes\sig(0,n+1)$ and $triv\otimes\sig(1,n).$ 
We choose
\begin{equation}
  \label{eq:6.2.11}
  \begin{aligned}
&e_1'=e_1,\\
&e_2':=\sum_{x\in S_{k-1}\times S_{n}\times S_1}
x\cdot[(\ep_1 -\ep_{k})\times\dots\times
(\ep_{m-1}-\ep_{k+m-2})(\ep_{k+n}-\ep_{k+m-1})] 
  \end{aligned}
\end{equation}
The intertwining operator  $I_{M',1,2}$ has a simple form on the vectors
invariant under $S_{k-1}\times S_n\times S_1$ transforming according
to $\sig(m,k+n-m-1)$ and $\sig(m-1,k+n-m).$ We can choose
multiples of
\begin{align}
&f_1:=\label{eq:6.2.12}\\
\sum_{x\in S_{k-1}\times S_{n}\times S_1}&x\cdot[(\ep_1-\ep_{k})
\times\dots\times (\ep_{m-1}-\ep_{k+m-2})(\ep_m-\ep_{k+m-1})],\notag\\
&\text{ in }\sig(m-1,k+n-m)\notag \\
&f_2:=\label{eq:6.2.13}\\
\sum_{x\in S_{k-1}\times S_{n}\times S_1}
&x\cdot[(\ep_1 -\ep_{k})\times\dots\times (\ep_{m-1}-\ep_{k+m-2})
\cdot\notag\\
\cdot(e_m+\dots+\ep_{k-1}+&\ep_k+\ep_{k+m}+\dots +
\ep_{k+n-1}-(k+n-2m+1)\ep_n)]\notag\\
&\text{ in }\sig(m,k+n-m-1) \notag
\end{align}
The fact that $f_1$  transforms according to $\sig(m,k+n-1)$
follows from (\ref{eq:6.2.7}). The fact that $f_2$ transforms
according to $\sig(m-1,k+n)$  is slightly more complicated. The product
$\prod (\ep_1-\ep_{k})\times\dots\times (\ep_{m-1}-\ep_{k+m-2})$
transforms according to $\sig(m-1,m-1)$ under $S_{2m-2}.$ The vector 
$$
\big[e_m+\dots+\ep_{k-1}+\ep_k+\ep_{k+m}+\dots+\ep_{k+n-1}-
(k+n-2m+1)\ep_{k+n}\big]
$$ 
is invariant under the $S_{k+n-2m-1}$ acting
on the coordinates\newline $\ep_m,\dots \ep_k,\ep_{k+m},\dots, \ep_{k+n-1}.$
Since $\sig(m,k+n-m-1)$ does not have such invariant vectors,
the product inside the sum in (\ref{eq:6.2.13}) must transform
according to $\sig(m-1,k+n-m).$ {The average under $x$ in
(\ref{eq:6.2.13}) is nonzero.} 
The operator $I_{M',2,3}$ maps $f_1$ and $f_2$ into multiples (using the
induction hypothesis) of the vectors $f_1',\ f_2'$ which are the
$S_n\times S_{k-1}\times S_1$ invariant vectors transforming according
to $\sig(m,k+n-1)$ and $\sig(m-1,k+n-m).$ 
The composition $I_{M',1,2}\circ I_{M',2,3}$ maps $e$ into a multiple of 
\begin{equation}
  \label{eq:6.2.14}
e':=\sum_{\sig\in S_n\times S_k} 
\sig\cdot[(\ep_1 -\ep_{n+1})\times\dots\times (\ep_m-\ep_{n+m})].
\end{equation}
The multiple is computable by using the induction hypothesis and the
expression of
\begin{itemize}
\item[] $e$ in terms of $e_1,\ e_2,$
\item[] $e'_1,\ e'_2$ in terms of $f_1,\ f_2,$ and
\item[] $e'$ in terms of $f'_1,\ f'_2.$
\end{itemize}
For example for the case $k=1,$ we get the following formulas.
\begin{equation}
  \label{eq:6.2.14a}
\begin{aligned}
&e=n(\ep_1+\dots +\ep_k)-k(\ep_{k+1}+\dots +\ep_{k+n}),\\
&e_1=(n+1)(\ep_1+\dots +\ep_{k-1})-(k-1)(\ep_k+\dots +\ep_{k+n}),\\
&e_2=n\ep_k-(\ep_{k+1}+\dots +\ep_{k+n}),\\
&f_1=n(\ep_1+\dots +\ep_{k-1})-(k-1)(\ep_k+\dots +\ep_{k+n-1}),\\   
&f_2=(\ep_1+\dots +\ep_{k-1})+(\ep_k+\dots +\ep_{k+n-1})-(k+n-1)\ep_{k+n},\\ 
&e'=-k(\ep_1+\dots +\ep_n) - n(\ep_{n+1}+\dots +\ep_{k+n}),\\ 
&e'_1=(n+1)(\ep_1+\dots +\ep_{k-1})-(k-1)(\ep_k+\dots +\ep_{k+n}),\\
&e'_2=-(\ep_k+\dots +\ep_{k+n-1})+n(\ep_{k+n}),\\
&f'_1=-(k+1)(\ep_1+\dots +\ep_{n})+n(\ep_{n+1}+\dots +\ep_{k+n-1}),\\
&f'_2=(\ep_1+\dots +\ep_{n})+(\ep_{n+1}+\dots +\ep_{k+n-1})-(k+n-1)\ep_{k+n}.
\end{aligned}
\end{equation}
Then
\begin{equation}
  \label{eq:6.2.15}
  \begin{aligned}
&e=\frac{k-1}{n+1}e_1 - \frac{k+n}{n+1}e_2,\\
&e'_1=\frac{k+n}{k+n-1}f_1 + \frac{k-1}{k+n-1}f_2,\\
&e'_2=\frac{1}{k+n-1}f_1 - \frac{n}{k+n-1}f_2,\\
&e'=\frac{k+n}{k+n-1}f'_1 - \frac{n}{k+n-1}f'_2.
  \end{aligned}
\end{equation}


\end{proof}

\subsection{$GL(k)\times G(n)\subset G(n+k)$}\label{sec:6.3} 
In the next sections we prove theorem \ref{sec:5.3} in the case of a
parabolic subgroup with Levi component $GL(k)\times G(n)$ for the
induced module  \begin{equation}\label{eq:6.3.1}
X_M((\nu_1)(\nu_0))= Ind_{M}^G[L(\chi_1)\otimes L(\chi_0)].
\end{equation}
{The    notation is as in section \ref{sec:6.1}.}

{The strings are
\begin{equation}
  \label{eq:6.3.2}
(-\frac{k-1}{2}+\nu,\dots , \frac{k-1}{2}+\nu)(-n+1 +\ep,\dots,-1 +\ep). 
\end{equation}}
Recall that $\ep=0$ when the Hecke algebra is type B, $\ep=1/2$ for type C,
and $\ep=1$ for type D, and 
\begin{equation}\label{eq:6.3.3}
  r_\sig(\nu):(V_\sig^*)^{W(M)}\longrightarrow (V_\sig^*)^{W(M)}.
\end{equation}
We will compute  $r_\sig(w_1,(\nu)(\nu_0))$ by induction on $k.$ In this
case the relevant $W-$types have multiplicity $\le 1$ so $r_\sig$ is a
scalar. 
\subsection{}\label{sec:6.4}
We start with  the special case $k=1$
when the maximal parabolic subgroup $P$ has Levi component  
$M= GL(1)\times G(n) \subset G(n+1).$ In type $D$ we assume $n\ge 1.$
Then 
\begin{equation}\label{eq:6.4.1}
X_M\mid_W=\sig[(n+1),(0)]\ +\ \sig[(1,n), (0)]\ +\ \sig[(n), (1)],
\end{equation}
and all the $W-$types occuring are relevant.
{In types B,C the operator $r_\sig(\nu)$ is the restriction to
$(V_\sig^*)^{W(M)}$ of the product  
\begin{equation}\label{eq:6.4.2}
r_{1,2}\circ\dots \circ r_{n,n+1}\circ r_{n+1}\circ r_{n,n+1}\circ
\dots \circ r_{1,2}
\end{equation}
as an operator on $V_\sig.$ Here $r_{i,j}$ is the $r_\sig(w,*)$
corresponding to the root $\ep_i-\ep_j$ and $r_{n+1}$  is the $r_\sig$
corresponding to  $\ep_{n+1}$ or $2\ep_{n+1}$ in types B and C. In
type D, the operator is
\begin{equation}\label{eq:6.4.2d}
r_{1,2}\circ\dots \circ r_{n,n+1}\circ\wti{r_{n,n+1}}\circ
\dots \circ r_{1,2}
\end{equation}
where $r_{i,i+1}$ are as before, and $\wti{r_{n,n+1}}$ corresponds to
$\ep_n+\ep_{n+1}.$ }
Since the multiplicities are 1, this is a
scalar.
\begin{proposition}\label{p:6.4} The scalar $r_\sig(w_1,((\nu)(\nu_0)))$ is   
\begin{equation}\label{eq:6.4.3}
\begin{array}{lll}
\qquad&{\sig_e(1,n)=\sig[(n),(1)]}\qquad&{\sig_o(1,n)=\sig[(1,n),
  (0)]}\\
&&\\
B&\frac{n+1-\nu}{n+1+\nu}&-\frac{n+1-\nu}{n+1+\nu}\\
&&\\
C&\frac{1/2+n-\nu}{1/2+n+\nu}
&\frac{1/2+n-\nu}{1/2+n+\nu}\cdot
\frac{1/2-\nu}{1/2+\nu}\\
&&\\
D&\frac{n-\nu}{n+\nu}&\frac{n-\nu}{n+\nu}\frac{1-\nu}{1+\nu}
\end{array}
\end{equation}
\end{proposition}
\begin{proof}
We do an induction on $n.$  

The reflection representation $\sig[(n),(1)]$ has dimension $n+1$ and
the usual basis $\{\ep_i\}.$ The $W(M)-$fixed vector
is $\ep_1.$ The representation $\sig[(1,n), (0)]$ has a basis
$\ep_i^2-\ep_j^2$ with the symmetric square action. 
The $W(M)-$fixed vector is
$\ep_1^2-\frac{1}{n}(\ep^2_2+\dots +\ep^2_{n+1}).$ 

The case $n=0$ for type $C$ is
clear; the intertwining operator is $1$ on $\mu_o(1,0)=triv$ and
$\frac{1/2-\nu}{1/2+\nu}$ on $\mu_e(0,1)=sgn.$ We omit the details for
type B. In type for $n=1,$ \ie  $D_2$, the
middle $W-$type in  (\ref{eq:6.4.1}) decomposes further
\begin{equation}
  \label{eq:6.4.4}
  \sig[(2), (0)] +\sig[(1),(1)]_I +\sig[(1), (1)]_{II}+\sig [(0), (2)].
\end{equation}
The representations $\sig[(1),(1)]_{I,II}$ are 1-dimensional with
bases $\ep_1\pm\ep_2.$ The result is straightforward in this case as well.

We now do the induction step. We give details for type B. 
In the case $\sig_e(1,n),$ 
embed $X_M$ in the induced module from the characters corresponding to
\begin{equation}
  \label{eq:6.4.5}
{(\nu)(-n)(-n+1 ,\dots ,-1 ).}  
\end{equation}

Write $M'=GL(1)\times GL(1)\times G(n-1)$ for the Levi component
corresponding to these three strings. Then the intertwining operator
$I:X_M((\nu)(\nu_0))\longrightarrow X_M((-\nu)(\nu_0))$ is the restriction of   
\begin{equation}\label{eq:6.4.6}
I_{M',1,2}((-n),(-\nu)(\nu_0))\circ I_{M',2}((-n)(\nu)(\nu_0))\circ
I_{M',1,2}(\nu,(-n),(\nu_0) ). 
\end{equation}
The $r_\sig$ have a corresponding decomposition
\begin{equation}
  \label{eq:6.4.7}
  (r_\sig)_{M',1,2}((-\nu),(-n)(\nu_0) )\circ
  (r_\sig)_{M',2}((-n)(\nu)(\nu_0))\circ
  (r_\sig)_{M',1,2}((\nu)(-n)(\nu_0)). 
\end{equation}
We need the restrictions of $\mu_e(1,n)$ and $\mu_o(1,n)$ to $W(M')$. We have 
{
\begin{align}\label{eq:6.4.8}
&Ind_{W(B_{n-1})}^{W(B_{n+1})}[\sig[(n-1), (0)]]= 
\sig[(n+1), (0)] + 2\sig[ (n),(1)] + 2\sig[(1,n),(0)] \notag \\ 
&\qquad  +2\sig[(1,n-1),(1)]
+\sig[(n-1), (2)] + \sig[(n-1),(1,1)] +\notag\\
&\quad +\sig[(2,n-1),(0)]+\sig[(1,1,n-1),(0)],\tag{a}\\ 
&\notag\\
&Ind_{W(B_{n})}^{W(B_{n+1})}[\sig[(n), (0)]]= \sig[(n+1), (0)]
+\sig[(n),(1)]+\sig[(1,n),(0)],\tag{b}\\
&\\
&Ind_{W(B_1)W(B_{n})}^{W(B_{n+1})}[\sig[(1),(0)]\otimes \sig[(n), (0)]]=
\sig[(n+1), (0)] +\sig[(1,n),(0)] \tag{c}\\
&Ind_{W(B_1)W(B_{n})}^{W(B_{n+1})}[\sig[(0),(1)]\otimes \sig[(n), (0)]]=
\sig[(n-1),(1)]\tag{d} 
\end{align}
}
Thus $\mu_e(1,n)$ occurs with multiplicity 2 in $X_{M'}.$ The $W(M')$ fixed
vectors are the linear span of $\ep_1,\ \ep_2.$ The intertwining
operators $I_{M',1,2}$ and $I_{M',2}$ are induced from maximal parabolic subgroups
whose Levi components we label $M_1$ and $M_2.$ Then
$\ep_1+\ep_2$ transforms like $triv\otimes triv$ under $W(M_{1})$ and
$\ep_1-\ep_2$ transforms like $sgn\otimes triv$. The vector $\ep_1$ is
fixed under $W(B_n)$ (which corresponds to $M_2$) and the vector
$\ep_2$ is fixed under $W(B_{n-1})$ and transforms like $\mu_o(1,n)$
under $W(B_n).$ The matrix $r_\sig$ is, according to (\ref{eq:6.4.7}),
{
\begin{equation}
  \label{eq:6.4.9}
  \begin{bmatrix}
    \frac{1}{2+\nu -n}&\frac{\nu-n+1 }{2+\nu-n }\\\frac{\nu-n+1
    }{1+\nu-n+1 }& \frac{1}{2+\nu-n } 
  \end{bmatrix}
\cdot 
\begin{bmatrix}
  1&0\\0&\frac{n-\nu}{n+\nu}
\end{bmatrix}
\cdot
  \begin{bmatrix}
    \frac{1}{1+\nu+n}&\frac{\nu+n}{1+\nu+n}\\\frac{\nu+n}{1+\nu+n}& 
\frac{1}{2+\nu+n}
  \end{bmatrix}.
\end{equation}
}

\medskip
So the vector $\ep_1$ is mapped into $\frac{n+1-\nu}{n+1+\nu}\ep_1$ as
claimed. 
\medskip
For $\sig_o(1,n)$ we apply the same method.
In this case the operator $I_{M',2}$ is the identity  because
in  the representation $\mu_o(1,n)$ the element $t_n$ corresponding to the
short simple root acts by 1. 

The calculation for type $D$ is analogous, we sketch some
details. We decompose the strings into 
\begin{equation}
  \label{eq:6.4.10}
  (\nu)(-n+1,\dots ,-1)(0),
\end{equation}
and $M'=GL(1)\times GL(n-1)\times GL(1).$ Then

  \begin{align}
  I_{M,1}((\nu)(\nu_0))&=  \label{eq:6.4.11}\\
I_{M',1,2}((-n+1,\dots,-1)(-\nu)(0))&\circ
  I_{M',1}((-n+1,\dots ,-1)(\nu)(0))\circ\notag\\
& I_{M',1,2}((\nu)(-n+1,\dots ,-1)(0)).\notag
  \end{align}

\end{proof}

\subsection{}\label{sec:6.5}
In this section we consider (\ref{eq:6.3.2}) for $k>1,$ $n\ge 1$ and    
the $W-$types $\sig_e(m,n+k-m)$ for $0\le m\le k$
{(notation as in definition \ref{d:5.8})}. These are the
$W-$types which occur in $X_M$, {with $M=GL(k)\times
  G(n)\subset G(k+n).$} 
\begin{proposition}\label{p:6.5}
The $r_\sig(w_1,((\nu)(\nu_0))$ for  $\sig=\sig_e(m,n+k-m)$ are
scalars.  They equal
\begin{description}
\item[Type B] 
\begin{equation}\label{eq:6.5.1}
\und{0\le j\le m-1}{\prod}
\frac{n+1-(-\frac{k-1}{2}+\nu)-j}{n+1+(\frac{k-1}{2}+\nu)-j}
\end{equation}
\item[Type C]
\begin{equation}\label{eq:6.5.2}
\und{0\le j\le m-1}{\prod}
\frac{n+1/2-(-\frac{k-1}{2}+\nu)-j}{n+1/2+(\frac{k-1}{2}+\nu)-j}
\end{equation}
\item[Type D]
\begin{equation}\label{eq:6.5.3}
\und{0\le j\le m-1}{\prod}
\frac{n-(-\frac{k-1}{2}+\nu)-j}{n+(\frac{k-1}{2}+\nu)-j}
\end{equation}
\end{description}
\end{proposition}
\begin{proof}
The proof is by induction on $k$.
The case $k=1$  was done in section \ref{sec:6.4} so we only need to do
the induction step. For types B,C  factor the intertwining operator as
follows. Decompose the string 
{\begin{equation}
  \label{eq:6.5.4}
((\nu')(\frac{k-1}{2}+\nu)(\nu_0)):=((-\frac{k-1}{2}+\nu,\dots
,\frac{k-3}{2}+\nu)(\frac{k-1}{2}+\nu)(\nu_0))   
\end{equation}
}
and let $M':=GL(k-1)\times GL(1)\times G(n),$ and $M''=GL(1)\times
GL(k-1)\times G(n).$ Thus
\begin{equation}\label{eq:6.5.5}
  \begin{aligned}
I_{M,1}&=I_{M'',2}((-\frac{k-1}{2}-\nu)(\nu')(\nu_0))\circ\\
&I_{M',1,2}((\nu')(-\frac{k-1}{2}-\nu)(\nu_0))\circ\\
&I_{M',2}((\nu')(\frac{k-1}{2}+\nu)(\nu_0))
  \end{aligned}
\end{equation}
$I_{M',1,2}$ and $I_{M',2}$ were computed earlier, while $I_{M'',2}$
is known by induction.   
Then
{
\begin{align}  
&\sig_e(m,n+k-m)\mid_{W(GL(k-1)\times W(G(n+1))}=\notag\\
&triv\otimes[\sig_e(1,n)+\sig_e(0,n+1)+\dots
\label{eq:6.5.6}\\
&\sig_e(m,n+k-m)\mid_{W(GL(k)\times W(G(n+k-1))}=\notag\\
&[(k)\otimes triv + (1,k-1)\otimes triv]+\dots
\label{eq:6.5.7}\\
&\sig_e(m,n+k-m)\mid_{W(GL(1)\times W(G(n+k-1))}=\notag\\
&triv\otimes[\sig_e(m-1,n+k-m)+\sig_e(m,n+k-1-m)]+\dots\label{eq:6.5.8}
\end{align} 
}
where $\dots$ denote $W-$types which are not spherical for $W(M),$ so
do not matter for the computations.

The $W-$type $\sig_e(m,n+k-m)\cong\bigwedge^m\sig_e(1,n+k-1).$ It occurs with
multiplicity 2 in $X_{M'}$ for $0<m<\min(k,n)$ and multiplicity 1 for
$m=\min(k,n).$  We will write out an explicit
basis for the invariant $S_1\times S_{k-1}\times W(B_n)$ vectors. Formulas
(\ref{eq:6.5.2})-(\ref{eq:6.5.4}) then 
come down to a computation with $2\times 2$ matrices as in the case
$k=1.$ Let
\begin{equation}
  \label{eq:6.5.9}
e:=\frac{1}{m!(k-m)!}\sum_{x\in  S_k}x\cdot[\ep_1\wedge\dots\wedge\ep_m].
\end{equation}

This is the $S_k\times W(B_n)$ fixed vector of $\sig_e(m,n+k-m).$ It
decomposes as 
\begin{equation}
  \label{eq:6.5.10}
  e=e_0+e_1=f_0+f_1
\end{equation}
where
\begin{equation}  \label{eq:6.5.11}
\begin{aligned}
  &e_0=\frac{1}{m!(k-1-m)!}\sum_{x\in S_{k-1}\times
  S_1}x\cdot[\ep_1\wedge\dots\wedge \ep_m],\\
  &e_1=\frac{1}{(m-1)!(k-m)!}\sum_{x\in S_{k-1}\times
  S_1}x\cdot[\ep_1\wedge\dots\wedge \ep_{m-1}]\wedge \ep_k,  \\ 
  &f_0=\frac{1}{m!(k-1-m)!}\sum_{x\in S_1\times S_{k-1}}
  x\cdot[\ep_2\wedge\dots\wedge \ep_{m+1}],\\
  &f_1=\frac{1}{(m-1)!(k-m)!}\sum_{x\in S_1\times S_{k-1}} 
  \ep_1\wedge x\cdot[\ep_2\wedge\dots\wedge \ep_m].
\end{aligned}
\end{equation}
Let also 
\begin{equation}
  \label{eq:6.5.12}
\begin{aligned}
&e'_0=e''_0=\frac{1}{(m-1)!(k-m)!}\sum_{x\in S_k}x\cdot
[\ep_1\wedge\dots\wedge\ep_m],\\
&e'_1=\sum_{x\in S_{k-1}\times S_1}
x\cdot[\ep_1\wedge\dots\wedge\ep_{m-1}\wedge(\ep_m-\ep_k)],\\
&e''_1=\sum_{x\in S_1\times S_{k-1}}
x\cdot[(-\ep_1+\ep_{m+1})\wedge\ep_2\wedge\dots\wedge\ep_{m+1}].
\end{aligned}
\end{equation}
Then
\begin{equation}\label{eq:6.5.13}
\begin{aligned}
&e_0=\frac{k-m}{k}e'_0 +\frac{m}{k}e'_1,&&\qquad e_1=\frac{m}{k}e'_0-\frac{m}{k}e'_1,\\
&e''_0=f_0+f_1, &&\qquad e''_1=f_0-\frac{k-m}{m}f_1.
\end{aligned}
\end{equation}
We now compute the action of the intertwining operators. The following
relations hold:
\begin{equation}
  \label{eq:6.5.14}
\begin{aligned}
&I_{M',2}(e_0)=e_0,
\qquad I_{M',2}(e_1)= \frac{n+\ep-(\frac{k-1}{2}+\nu)}{n+\ep+(\frac{k-1}{2}+\nu)}e_1,\\
&I_{M',12}(e'_0)=e''_0,\qquad
I_{M',12}(e'_1)=\frac{2\nu-1}{2\nu+k-1}e''_1,\\
&I'_{M'',2}(f_0)=\prod_{0\le j\le m-2}\frac{n+\ep
  -(-\frac{k-1}{2}+\nu)-j}{n+\ep+(\frac{k-3}{2}+\nu)-j}f_0,\\
&I_{M'',2}(f_1)=\prod_{0\le j\le m-1}
\frac{n+\ep-(-\frac{k-1}{2}+\nu)-j}{n+\ep+(\frac{k-3}{2}+\nu)-j}f_1, &&
\end{aligned}  
\end{equation}
where $\ep=1$ in type B, $\ep=1/2$ in type C, and $\ep=0$ in type D.
Then 
\begin{equation}
  \label{eq:6.5.15}
  I_{M',2}(e_0+e_1)=e_0+ \frac{n+\ep-(\frac{k-1}{2}+\nu)}{n+\ep+(\frac{k-1}{2}+\nu)}e_1.
\end{equation}
Substituting the expressions of $e_0,\ e_1$ in terms of $e'_0,\ e'_1,$  we get
\begin{equation}
  \label{eq:6.5.16}
  [\frac{k-m}{k}+\frac{m}{k}
\frac{n+\ep-(\frac{k-1}{2}+\nu)}{n+\ep+(\frac{k-1}{2}+\nu)}]e_0'
  + \frac{m}{k}[1-\frac{n+\ep-(\frac{k-1}{2}+\nu)}{n+\ep+(\frac{k-1}{2}+\nu)}]e'_1.
\end{equation}
Applying $I_{M,2}$ to this has the effect that $e'_0$ is sent
to $e''_0$ and the term in $e'_1$ is
multiplied by $\frac{2\nu-1}{2\nu+k-1}$ and $e'_1$ is replaced by
$e''_1.$ Substituting the formulas for $e''_0$ and $e''_1$ in terms of
$f_0,\ f_1,$ and applying $I_{M'',2},$ we get the claim of the proposition.
\end{proof}

\subsection{}\label{sec:6.6}
We now treat the case $\sig=\sig_o(m,n+k-m).$ We assume $n>0$ or else
these $W-$types do not occur in the induced module $X_M.$ 
\begin{proposition}\label{p:6.6}
The $r_\sig(w_1,((\nu)(\nu_0))$ are scalars. They equal 
{
\begin{equation}
  \label{eq:6.6.1}
\prod_{0\le j\le
  m-1}\frac{(\nu-\frac{k-1}{2})-(1-\ep)+j}{(\nu+\frac{k-1}{2})
              -(-n -\ep)-j}\cdot
\frac{(-n-\ep )-(-\nu+\frac{k-1}{2})+j}{(1-\ep )-(-\nu-\frac{k-1}{2})-j}
\end{equation}
}
\end{proposition}
\begin{proof}
The intertwining operator $I_M(\nu)$
decomposes in the same way as (\ref{eq:6.5.5}). Furthermore,
{$\sig_o(m,n+k-m)=\bigwedge^m\sig_o(1,n+k-1).$} 
The difference from the cases $\sig_e$ is that while $\sig_e(1,n+k-1)$
is the reflection representation, and therefore realized as the
natural action on $\ep_1,\dots \ep_{n+k},$ $\sig_o(1,n+k-1)$ occurs in
$S^2\sig_e(1,n+k-1)$, generated by $\ep_i^2-\ep_j^2$ with $i\ne j.$ We
can apply the same technique as for $\sig_e(m,n+k-m)$, and omit the details. 
\end{proof}
{
\subsection*{$\mathbf{GL(k)\subset G(k)}$ in types B, C}
The formulas in proposition \ref{p:6.5} and \ref{p:6.6} hold with $n=0.$ 
The proof is the same, but because $n=0,$ $(\nu_0)$ is not
present. The operator $I_{M',2}$ is an intertwining operator in
$SL(2)$ and therefore simpler.

}

\subsection{$\mathbf{GL(k)\subset G(k)}$ in type D}\label{sec:6.7} 
In this section we consider the maximal
Levi components $M:=GL(k)\subset G(k)$ and $M':=GL(k)'\subset G(k)$ for
type $D_n$. The parameter corresponds to the string
$(\nu):=(-\frac{k-1}{2}+\nu,\dots ,\frac{k-1}{2}+\nu)$ or
$(\nu'):=(-\frac{k-1}{2}+\nu,\dots ,-\frac{k-1}{2}-\nu).$ 
\begin{description}
\item[$k$ even] The $W-$structure of $X_M((\nu))$ and $X_{M'}((\nu)')$
is $\sig_e[(n-r),(r)]$ for $0\le r<k/2,$ and
{$\sig_e[(k/2),(k/2)]_{I},$ or $\sig_e[k/2),(k/2)]_{II}$}
respectively, with multiplicity 1.  There are intertwining operators
\begin{equation}
  \label{eq:6.7.1}
  \begin{aligned}
&I_M((\nu)):X_M((\nu))\longrightarrow X_M((-\nu)),\\    
&I_{M'}((\nu)'):X_{M'}((\nu)')\longrightarrow  X_{M'}((-\nu)').    
\end{aligned}
\end{equation}
corresponding to the shortest Weyl group element changing $((\nu))$ to
$((-\nu)).$ They determine scalars $r_\sig((\nu))$ and $r_\sig((\nu)').$ 
\item[$k$ odd] The $W-$structure in this case is $\sig_e[(n-r),(r)]$
  with $0\le r\le [k/2]$ for both $X_M$ and $X_{M'},$ again with
  multiplicity 1. In this case there is a shortest Weyl group element
  which changes $((\nu))$ to $((-\nu)'),$ and one which changes
  $((\nu)')$ to $((-\nu)).$ These  elements give rise to intertwining
  operators
 \begin{equation}
  \label{eq:6.7.2}
  \begin{aligned}
&I_M((\nu)):X_M((\nu))\longrightarrow X_{M'}((-\nu)'),\\    
&I_{M'}((\nu)):X_{M'}((\nu)')\longrightarrow  X_{M}((-\nu)).
\end{aligned}    
\end{equation} 
Because the $W-$structure of $X_M$ and $X_{M'}$ is the same, and
$W-$types occur with multiplicity 1, these intertwning operators
define scalars $r_\sig(\nu)$ and $r_\sig((\nu)').$
\end{description}
\begin{proposition}
  \label{p:6.7}
The scalars $r_\sig((\nu))$ and $r_\sig((\nu)')$ are
\begin{equation}
  \label{eq:6.7.3}
r_{\sig_e[(n-r),(r)]}((\nu))=  
\prod_{0\le j<r}\frac{(\frac{k-1}{2}-\nu)-j}{(\frac{k-1}{2}+\nu)-j} .
\end{equation}
These numbers are the same for $((\nu))$ and $((\nu)'),$ and
representations with subscripts $I,II$ they depend only on $r.$
\end{proposition}

\subsection{Proof of theorem \ref{sec:5.3}}\label{sec:6.8} 
We use the
results in the previous sections to prove the theorem in general.
We give the details in  the case of the group of type B and
{$W-$types} $\sig_e.$  Thus the Hecke algebra is type C.
There are no significant changes in the proof for the other
cases. Recall the notation from 
section \ref{sec:2.3}. 
{Conjugate $\nu$ the middle element for the
nilpotent orbit with partition $(2x_0,\dots , 2x_{2m})$ 
so that it is dominant (\ie
the coordinates are in decreasing order),
$$
\nu=(x_{2m}-1/2,\dots ,x_{2m}-1/2,\dots, x_0-1/2,\dots ,x_0-1/2,\dots
,1/2,\dots , 1/2) 
$$
}
Then $\nu$ is dominant, so $X(\nu)$ has a unique irreducible quotient $L(\nu).$
We factor the long intertwining  operator so that 
\begin{equation}\label{eq:6.8.1}
X(\nu)\ovr{I_1}{\longrightarrow} X_e(\nu)\ovr{I_2}{\longrightarrow} X(-\nu),
\end{equation}
{where $X_e$ was defined in section \ref{sec:5.3}.}
The claim will follow if the decomposition has the property that the
operator $I_1$ is onto, and $I_2$ is into, when restricted to
the $\sig_e$ isotypic component.

\smallskip

\noindent{\textbf{Proof that $\mathbf{I_1}$ is onto.}} 
The operator $I_1$ is a composition of several operators. First take
the long intertwining operator induced from the Levi component $GL(n),$
\begin{equation}\label{eq:6.8.2}
X(x_{2m}-1/2,\dots, 1/2)\longrightarrow X(1/2,\dots ,x_{2m}-1/2),
\end{equation}
corresponding to the shortest Weyl group element that permutes the entries
of the parameter from {decreasing order to  increasing
order.}
The image is the induced from the
corresponding irreducible spherical module $L(1/2,\dots ,x_{2m}-1/2)$  on
$GL(n).$ In turn this is induced irreducible from 1-dimensional 
spherical characters on a $GL(x_0)\times\dots\times GL(x_{2m})$ Levi component
corresponding to the strings
$$
(1/2,\dots ,x_{0}-1/2)\dots (1/2,\dots,x_{2m}-1/2)
$$
or any permutation thereof. This is well known by results of
Bernstein-Zelevinski in the $p-$adic case, \cite{V1} for the real
case. 

Compose with the intertwining operator 
\begin{equation}\label{eq:6.8.3}
X(\dots (1/2,\dots,x_{2m}-1/2))\longrightarrow 
X(\dots (-x_{2m}+1/2,\dots ,-1/2)),
\end{equation}
all other entries unchanged. This intertwining operator is induced
from the standard long intertwining operator on $G(x_{2m})$ which has image
equal to the trivial representation. The image is an induced module
from characters on $GL(x_0)\times\dots \times GL(x_{2m-1})\times
G(x_{2m}).$ Now compose with the intertwining
operator
\begin{align}
  \label{eq:6.8.4}
 X&(\dots (1/2,\dots, x_{2m-1}-1/2)(-x_{2m}+1/2,\dots ,-1/2)) 
\\
&\longrightarrow 
X(\dots (-x_{2m-1}+1/2,\dots ,-1/2)(-x_{2m}+1/2,\dots ,-1/2))\notag
\end{align}
(again all other entries unchanged). 
This is $I_{M,2m-1}$ {defined in (\ref{eq:6.1.8})}, 
so its restriction of (\ref{eq:6.8.4}) to the $\sig_e$
isotypic component is an isomorphism. Now compose this operator with
the one corresponding to 
\begin{align}
  \label{eq:6.8.5}
 X&(\dots (1/2,\dots, x_{2m-2}{-1/2})
(-x_{2m-1}+1/2,\dots ,1/2)\dots) 
\\
&\longrightarrow 
X(\dots (-x_{2m-1}{+1/2},\dots ,x_{2m-2}-1/2)\dots )\notag
\end{align}
with all other entries unchanged. 
This is induced from 
$$
GL(x_0)\times\dots\times GL(x_{2m-3})
\times GL(x_{2m-2}+x_{2m-1})\times G(x_{2m})
$$
and the image is the
representation induced from the  character corresponding to the string 
$$
(-x_{2m-1}-1/2,\dots,-1/2,1/2,\dots ,x_{2m-2}) 
\text{ on } GL(x_{2m-2}+x_{2m-1}).
$$ 
Now compose further with the intertwining operator 
\begin{align}\label{eq:6.8.6}
X&(\dots (-x_{2m-1}+1/2,\dots,x_{2m-2}-1/2)(-x_{2m}-1/2,\dots ,-1/2))\\
&\longrightarrow
X((-x_{2m-1}+1/2,\dots,x_{2m-2}-1/2)\dots (-x_{2m}-1/2,\dots ,-1/2))\notag
\end{align}
from the representation induced from
$$
GL(x_0)\times\dots\times
GL(x_{2m-3})\times GL(x_{2m-2}+x_{2m-1})\times G(x_{2m})
$$
to the induced from  
$$
GL(x_{2m-2}+x_{2m-1})\times GL(x_0)\times\dots\times 
GL(x_{2m-3})\times G(x_{2m}).
$$

By lemma \ref{6.2},  this intertwining operator is an isomorphism on
any $\sig_e$ isotypic component. In fact, because the strings are
strongly nested, the irreducibility results for $GL(n)$, \ref{p:3.3} imply that
the induced modules are isomorphic. 

We have constructed a composition of intertwining operators 
from the standard module $X(\nu)$ where the
coordinates of $\nu$ are positive and in decreasing order (\ie
dominant) to a module induced from 
$$
GL(x_{2m-2}+x_{2m-1})\times GL(x_0)\times\dots\times
GL(x_{2m-3})\times G(x_{2m})
$$
corresponding to the strings
\begin{align*}
((-x_{2m-1}+1/2,&\dots ,x_{2m-2}-1/2)(1/2,\dots ,x_{0}-1/2),\dots \\
&\dots,(1/2,\dots,x_{2m-3}-1/2) (-x_{2m}+1/2,\dots ,-1/2))
\end{align*}
so that the restriction to any $\sig_e$ isotypic component is onto. We
can repeat the procedure with $x_{2m-4},x_{2m-3}$ and so on to get an
intertwining operator from $X(\nu)$ to the induced from
\[
GL(x_{2m-1}+x_{2m-2})\times \dots\times GL(x_1+x_0)\times G(x_{2m})
\] 
corresponding to the strings
\begin{align*}
((-x_{2m-1}+1/2,&\dots ,x_{2m-2}-1/2)\dots (-x_1+1/2,\dots ,x_{0}-1/2),\\
&(-x_{2m}+1/2,\dots ,-1/2)).
\end{align*}
This is the operator $I_1,$ and it is onto on the $\sig_e(*)$ isotypic
components. 

\medskip

\noindent{\bf Proof that $I_2$ is into.} We now deal with $I_2.$ Consider the group
$G(x_{1}+x_{0}+x_{2m})$ and the Levi component $M=GL(x_{1}+x_{0})\times
G(x_{2m}).$ Let $M'$ be the Levi
component 
\begin{equation}
  \label{eq:6.8.7}
M':=GL(x_{2m-1}+x_{2m-2})\times\dots\times GL(x_3+x_2)\times GL(x_{1})\times
GL(x_{0})\times G(x_{2m}).   
\end{equation}
Then $X_e$ embeds in $X_{M'}(\dots (-x_1+1/2,\dots -1/2)(1/2,\dots
,x_0-1/2)(-x_{2m}+1/2,\dots ,-1/2)).$
The intertwining operator $I_{M',m+1}$ which changes the string
$(x_0-1/2,\dots ,1/2)$ to $(-x_0+1/2,\dots ,-1/2)$ is an isomorphism
on the $\sig_e$ $W-$types, by the results in sections
\ref{sec:6.1}-\ref{sec:6.5}. Since the strings are strongly nested,
the operators $I_{M,i,i+1}$ are all isomorphisms, so we can construct
an intertwining operator to an induced module $X_{M''}(\nu'')$ where 
\begin{align}
  \label{eq:6.8.8}
&M''=GL(x_1)\times GL(x_0)\times
    GL(x_{2m-1}+x_{2m-2})\times\dots\times G(x_{2m}),\notag\\
&\nu''=(-x_1+1/2,\dots -1/2)(-x_0+1/2,\dots ,-1/2)\dots )
\end{align}
which is an isomorphism on the $\sig_e$ isotypic components. Repeating
this argument for $x_3,\ x_2$ up to $x_{2m-1},x_{2m-2}$ we get an
intertwining operator from $X_e$ to an induced module
$X_{M^{(3)}}(\nu^{(3)})$ where 
\begin{equation}
  \begin{aligned}\label{eq:6.8.9}
&M^{(3)}:=GL(x_1)\times GL(x_0)\times\dots\times GL(x_{2m-1})\times
  GL(x_{2m-2})\times G(x_{2m})\\
&\nu^{(3)}:=(-x_{1}+1/2,\dots ,1/2)(-x_{0}+1/2,\dots ,-1/2)\dots\\ 
  &(-x_{2m-1}+1/2,\dots ,-1/2)(-x_{2m-2}+1/2,\dots,-1/2)
   (-x_{2m}+1/2,\dots ,-1/2).    
  \end{aligned}
\end{equation}
which is an isomorphism on the $\sig_e$ isotypic components. Let
{
$$
M^{(4)}:=GL(x_1)\times GL(x_{0})\times\dots\times 
GL(x_{2m-1})\times  G(x_{2m})
$$
and let $\nu^{(4)}$ be the same as $\nu^{(3)}$ but the last string is
viewed as giving a parameter of a 1-dimensional representation on
$GL(x_{2m}).$ Then $M^{(4)}\subset M^{(3)},$ and
$X_{M^{(3)}}(\nu^{(3)})$ is a submodule of $X_{M^{(4)}}(\nu^{(4)})$.
The induced module from $M^{(4)}$ to $M^{(5)}:=GL(x_{2m}+\dots +x_0)$
is irreducible because the strings are strongly nested on the $GL$ factors.  
Thus the intertwining operator which takes $\nu^{(4)}$ to $-\nu$ is an
isomorphism on $X_{M^{(4)}}(\nu^{(4)}).$   So the  induced intertwining
operator to $G$ is therefore injective and maps to  $X(-\nu).$ The
composition of all these operators is $I_1$, and is therefore
injective on the $\sig_e$-isotypic components. The proof is complete
in this case. }

The case of $\sig_o$ is similar, and we omit the details.

\section{Necessary conditions for unitarity}\label{sec:7}

\subsection{}\label{sec:7.1}

We will need the following notions.

\begin{definition}
We will say a spherical irreducible module $L(\chi)$ is \textbf{r-unitary} if
the form is positive on all the relevant $W-$types. Similarly,
an induced module {$I_M(\chi):=Ind_M^G[L_M(\chi)]$} is
\textbf{r-irreducible} if all relevant $W-$types  occur with the same
multiplicity in {$I_M(\chi)$} as in $L(\chi)$. 
  
\end{definition}

\subsection{}\label{sec:7.2}
We recall (\ref{eq:6.1.4}),
\begin{equation*}
  \ep=
  \begin{cases}
    1/2 &\text{ G of type } B,\quad (\bH \text{ of type } C)\\
     0  &\text{ G of type } C,\quad (\bH \text{ of type } B)\\
     1  &\text{ G of type } D.
  \end{cases}
\end{equation*}

\begin{definition}
  \label{d:7.2}
A string of the form $(f+\nu,\dots ,F+\nu)$ with $f,F\in\ep+\bZ$ is
called \textrm{adapted}, if it is

\qquad of \textrm{even} length for $G$ of  type B,

\qquad of \textrm{odd } length for $G$ of types C,D.

Otherwise we say the string is \textrm{not} adapted. 
\end{definition}
We will consider the following case. Let $\vO\subset\vg$ correspond to
the partition
\begin{equation}
  \label{eq:7.2.1}
  \vO\longleftrightarrow ((a_1,a_1),\dots (a_r,a_r);d_1,\dots ,d_l)
\end{equation}
so that $\vO$ meets the Levi component $\vm=gl(a_1)\times\dots \times
gl(a_r)\times \vg(n_0)$, {with
  $2n_0+[{1-\ep}]=d_1+\dots+d_l$, where $[x]$ is the integer part of
  $x$}. The intersection of $\vO$ 
with each $gl(a_i)$  is the principal nilpotent, and the intersection with
$\vg(n_0)$ is the \textbf{even} nilpotent orbit {$\vO_0$}
with partition $(d_1,\dots ,d_l).$ Let 
\begin{equation}
  \label{eq:7.2.2}
  \begin{aligned}
&  \chi_i=(f_i+\nu_i,\dots ,F_i+\nu_i),\qquad 1\le i\le r,\\
&  \chi_0=\vh_0/2, {\text{where $\vh_0$ is a neutral
    element for $(d_1,\dots ,d_l).$ }}
  \end{aligned}
\end{equation}
and $\chi$ be the parameter obtained by concatenating the $\chi_i.$ 
Then $L(\chi)$ is the spherical subquotient of
\begin{equation}
  \label{eq:7.2.3}
  Ind_{M}^G[\bigotimes_{1\le i\le r} L(\chi_i)\otimes L(\chi_0)]
\end{equation}
The next theorem gives necessary conditions for the unitarity of
$L(\chi).$ 
\begin{theorem}
  \label{t:7.2}
{In types B,C, the nilpotent orbit $\vO_0$ is arbitrary.
In type D assume that either
$\vO_0\ne (0)$, or else that the rank is even.} The representation
$L(\chi)$ is unitary \textbf{only if} 
\begin{enumerate}
\item Any string that is \textrm{not adapted} can be written in  the form
\begin{equation}
  \label{eq:7.2.4}
  (-E+\tau,\dots ,E-1+\tau)\qquad 0<\tau\le 1/2,\ E\equiv \ep (mod\ \bZ).
\end{equation}
\item Any string that is adapted can be written in the form
\begin{equation}
  \label{eq:7.2.5}
  \begin{aligned}
  &(-E+\tau,\dots ,E+\tau)\qquad 0<\tau\le 1/2,\ E\equiv\ep (mod\ \bZ),\\
&\text{ or } \\
  &(-E-1+\tau,\dots ,E-1+\tau)\qquad 0<\tau\le 1/2,\ E\equiv\ep (mod\ \bZ).
  \end{aligned}
\end{equation}
\end{enumerate}
\end{theorem}
This is simply the fact that the $\nu_j$ satisfy $0< \nu_j<1/2$ or
$1/2<\nu_j<1$ in theorem \ref{thm:3.1}. The proof will be given in the
next sections. It is by induction on
the dimension of $\vg,$ the number of strings with coordinates in an
$A_\tau$ ({definition after (\ref{eq:3.3.5})}) with
$\tau\ne 0,$ and by downward induction on the dimension 
of $\vO.$  The unitarity of the representation when there are no
coordinates in any $A_\tau$ with $\tau\ne 0$ is done in section \ref{sec:9}. 
\subsection{}\label{sec:7.3}
{
The proposition in this section is a restatement of theorem \ref{t:7.2} for
the case of a parameter of the form (\ref{eq:7.3.1}).  
It is the first case in the initial step of the induction proof of
theorem \ref{t:7.2}. I have
combined the two cases in (2) into a single string $(-E +\nu,\dots ,
E+\nu)$ with $0<\nu<1$ by changing $(-E-1+\tau,\dots ,E-1+\tau)$ into
$(-E +(1-\tau),\dots , E+(1-\tau)).$ This notation seemed more
convenient for the case when there is a single such string present.

Consider the representation $L(\chi)$  corresponding to the strings  
\begin{equation}\label{eq:7.3.1}
(a+\ep+\nu,\dots ,A+\ep+\nu)(-x_0+\ep,\dots ,-1+\ep),\quad
|a|\le A,\ 0<\nu<1,
\end{equation}
where $a,A\in\bZ,$ and $\ep$ is as in \ref{d:7.2}. This is the case
when $L(\chi)$ is the spherical subquotient of an induced from a
character on a maximal parabolic subalgebra of the form
$gl(A-a+1)\times \fg(x_0).$ The second string may not be present;
these are the cases $x_0=-1$ for $\fg$ of type $B,\ D,$ $x_0=0$ for
type C.

\begin{proposition}\label{p:7.3} In type D, assume that if there is no
  string $(-x_0+\ep,\dots ,-1+\ep),$ then 
  $A-a+1$ is even. Let $L(\chi)$ correspond to 
(\ref{eq:7.3.1}). Then $L(\chi)$ is r-unitary if and only if
$a+\ep=-A-\ep,$ and the following hold.
\begin{enumerate}
\item {Assume that
  $(a+\ep+\nu,\dots ,A+\ep +\nu)$ is adapted. If $x_0=A-a+1$, then
  $0\le\nu<1,$ otherwise $\nu=0.$}
\item If $(a+\ep+\nu,\dots ,A+\ep +\nu)$ is not adapted, then $0\le\nu<1/2.$
\end{enumerate}
\end{proposition}
\begin{proof}
This is a corollary of the formulas in section \ref{6.2}.
\end{proof}
}
\subsection{Initial Step}\label{sec:7.4}
We do the case when there is a single $A_\tau$ with $0<\tau<1/2,$ and the
coordinates form a single string. We write the string as in
(\ref{eq:7.3.1}), $(a+\ep+\nu,\dots ,A+\ep +\nu)$ with $0<\nu<1.$ We let
$\vO$ be the nilpotent orbit with partition
$((A-a+1,A-a+1);d_1,\dots ,d_l)$. Let  
$\vm:=gl(A-a+1)\times \vg(n_0),$  and let $\vO_0$ be the
intersection of $\vO$ with $\vg(n_0).$ In type D, either $\vO_0\ne
(0)$ or else $A-a+1$ is even. The statement of theorem \ref{t:7.2} is
equivalent to the following proposition.
\begin{proposition}
  \label{p:7.4}
Assume $\vO_0$ is even, and $\chi$ is attached to $\vO.$ 
Then $L(\chi)$ is r-unitary only if $a+\ep=-A-\ep,$ and the following hold. 
\begin{enumerate}
\item If $(a+\ep+\nu,\dots ,A+\ep +\nu)$ is adapted, then $\nu=0$,
  unless there is $d_j=A-a+1,$   in which case $0\le\nu<1.$ 
\item If $(a+\ep+\nu,\dots ,A+\ep +\nu)$ is not adapted, then $0\le\nu<1/2.$
\end{enumerate}
\end{proposition}
\begin{proof}
We do the case of $G$ of type $C$ only, the others are similar. So $\ep=0,$
and adapted means the length of the string is odd, not adapted means
the length of the string is even. 
The nilpotent orbit $\cCO_0$ corresponds to the partition
$(2x_0+1,\dots ,2x_{2m}+1)$ and the parameter has strings 
\[
(1,\dots x_0)(0,1,\dots ,x_1)\dots (1,\dots ,x_{2m}).
\]
The partition of $\cCO$ is $(A-a+1,A-a+1,2x_0+1,\dots ,2x_{2m}+1).$

We want to show that if $A+a>0,$ or if $A+a=0$ and there is no $x_i=A,$
then $L(\chi)$ is \textit{not} r-unitary. We do an upward
induction on the rank of $\vg$ and a downward induction on the
dimension of $\cCO$.  {In the
  argument below with deformations of strings, we use implicitly the
  irreducibility results from section \ref{sec:2.3c}.} So the first case
is when $\cCO$ is maximal, \ie  the principal nilpotent ($m=0$). The
claim follows from proposition \ref{p:7.3}. So we assume that $m$
is strictly greater than 0. 

\noindent\textbf{Assume} $\mathbf{x_{2i}<A\le x_{2i+1}}$ {for some
  $i.$} This case 
includes the possibility $x_{2m}<A.$ We will
show by induction on rank of $\vg$ and dimension of $\vO$ 
that the form is negative on a {$W$-type of the form}
$\sig[(n-r),(r)].$ So we use the module $X_e$ {(notation
  as in \ref{sec:5.3})}. If there is any pair $x_{2j}=x_{2j+1},$ the
module $X_e$ is unitarily 
induced from $GL(2x_{2j}+1)\times G(n-2x_{2j}-1)$ and all $W-$types
$\sig[(n-r),(r)]$ have the same multiplicity in $L(\chi)$ as in $X_e.$ We can
\textit{remove} the string corresponding to $(x_{2j}x_{2j+1})$ in $X_e$ as
explained in section \ref{sec:3.2}, lemma (3). By induction on
rank we are done. Similarly we can remove any pair $(x_{2j},x_{2j+1})$
such that either $x_{2j+1}\le |a|$ or $A\le x_{2j}$ as follows. Let
$M:=GL(x_{2j}+x_{2j+1}+1)\times G(n-x_{2j}-x_{2j+1}-1).$ There is
$\chi_M$ such that $L(\chi)$ is the spherical subquotient of
\begin{equation}
  \label{eq:7.4.1}
  Ind_M^G[L(-x_{2j+1},\dots , x_{2j})\otimes L(\chi_M)].
\end{equation}
{Precisely, $\chi_M$ is obtained from $\chi$ by
  removing the entries \newline $(1,\dotsc,x_{2j}),(0,1,\dotsc,x_{2j+1}).$}
Write 
\begin{equation}\label{eq:7.4.2}
\chi_t:=(-x_{2j+1}+t,\dots, x_{2j}+t;\chi_M).  
\end{equation}
The induced module
\begin{equation}
  \label{eq:7.4.3}
 X_e(\chi_t):=  Ind_M^G[L(-x_{2j+1}+t,\dots , x_{2j}+t)\otimes L(\chi_M)].
\end{equation}
has  $L(\chi_t)$ as its irreducible spherical subquotient. 
{For $0\le t\le
\frac{x_{2j+1}-x_{2j}}{2},$ the multiplicities of $\sig[(n-r,r)]$ in
$L(\chi_t)$ and $X_e(\chi_t)$ coincide.}
Thus the signatures on the $\sig[(n-r),(r)]$ in $L(\chi_t)$ are constant for
$t$ in the above interval. 
At $t=\frac{x_{2j+1}-x_{2j}}{2},$ $X_e(\chi_t)$ is unitarily induced
from $triv\otimes X_e'$ on $GL(x_{2j}+x_{2j+1}+1)\times
G(n-x_{2j}-x_{2j+1}-1)$ and we can \textit{remove} the string
corresponding to $(x_{2j}x_{2j+1}).$ The induction hypothesis
applies to $X_e'.$ 

\medskip
When $A+a=0,$ by the above argument, we are reduced to the case
\begin{equation}
  \label{eq:7.4.4}
  \cCO_0\longleftrightarrow (2x_0+1,2x_1+1,2x_2+1),\qquad x_0<A< x_1\le x_2.
\end{equation}
We reduce to (\ref{eq:7.4.4}) when $A+a>0$ as well. 
{We assume $2m=2i+2,$} 
since pairs $(x_{2j},x_{2j+1})$ with $A\le x_{2j}$
can be \textit{removed}. Suppose there is a pair
$(x_{2j},x_{2j+1})$ such that  $|a|<x_{2j+1},$ and $j\ne i.$ 
The assumption is that
$x_{2i}<A\le x_{2i+1}$ so $x_{2j+1}\le x_{2i}<A.$

We consider the deformation   $\chi_t$ in (\ref{eq:7.4.2})  with  
\begin{align*}
0\le t<\nu,&\qquad a<0,\\
-\nu<t\le 0,&\qquad a\ge 0.  
\end{align*}
In either case $X_e(\chi_t)=X_e(\chi),$ so the multiplicities of the
$\sig[(n-r),(r)]$ do not change until $t$  reaches $\nu$ in the
first case, $-\nu$ in the second case. If the signature on some
$\sig[(n-r),(r)]$ isotypic component is positive semidefinite on
$L(\chi)$, the same has to hold when $t=\nu$ or $-\nu$
respectively. The corresponding nilpotent orbit for this parameter is
strictly larger, but it has two strings with coordinates which are not
integers {(so the induction hypothesis does not
  apply yet)}.  For example, if $a<0,$ the strings for  $X_e(\chi_\nu)$
are (aside from the ones that were unchanged)
\begin{equation}
  \label{eq:7.4.5}
(-x_{2j+1}+\nu,\dots ,A+\nu),\qquad (a+\nu,\dots ,x_{2j}+\nu).  
\end{equation}
We can deform the parameter further by replacing the second string by
$(a+\nu-t',\dots , x_{2j}+\nu-t')$ with $0\le t'<\nu.$ The 
strings of the corresponding $X_e$ do not change until $t'$ reaches
$\nu.$ At $t'=\nu$ the corresponding nilpotent orbit $\cCO'$ has partition 
\begin{equation}
  \label{eq:7.4.6}
(\dots,2|a|+1,\dots,\widehat{2x_{2j+1}+1},\dots,
A+x_{2j+1}+1,A+x_{2j+1}+1,\dots)  
\end{equation}
which contains $\cCO$ in its closure. Since $x_{2j+1}<A,$ the
induction hypothesis applies. The form is indefinite on a
$W-$type $\sig[(n-r),(r)]$, so this holds for the original $\chi$ as
well.  

\medskip
We have reduced to case {(\ref{eq:7.4.4})}, \ie the
partition of $\vO_0$ has just three terms $(2x_0+1,2x_1+1,2x_2+1).$
We now reduce further to the case 
\begin{equation}
  \label{eq:7.4.7}
 \cCO_0\longleftrightarrow (2x_0+1),\qquad x_0<A.
\end{equation}
which is the initial step.  

Let $I(t)$ be the induced module coresponding to the strings
{
\begin{equation}
  \label{eq:7.4.8}
  (-x_2+t,\dots,x_1+t)(a+\nu,\dots,A+\nu)(-x_0,\dots ,-1).
\end{equation}
}
\ie induced from
\begin{equation}
  \label{eq:7.4.9}
GL(x_1+x_2)\times GL(-a+A+1)\times G(x_0).  
\end{equation}
Consider the irreducible spherical module for the last two strings in
(\ref{eq:7.4.8}), inside the induced module from the Levi component
$GL(-a+A+1)\times G(x_0)\subset G(-a+A+1+x_0).$ 
By section \ref{sec:7.1}, the form is negative on
$\sig[(x_0-a+A),(1)]$ if $x_0<a,$ 
negative on $\sig[(A),(x_0+1-a)]$ if $a\le x_0.$ In
the second case the form is positive on all $\sig[(A+r),(x_0+1-a-r)]$ for
{$1<r<x_0+1-a.$} 
So let $r_0:=1$ or $x_0+1-a$ depending on these two cases.
The multiplicity formulas from section \ref{6.2} imply that 
$$
[\sig[(n-r_0),(r_0)]\ :\ I(t)] = [\sig[(n-r_0),(r_0)]\ :\
L(\chi)]\quad \text{ for }\quad 0\le t\le \frac{x_2-x_1}{2}. 
$$
Thus signatures do not change when we  deform $t$ to
$\frac{x_2-x_1}{2},$ where $I(t)$ is 
unitarily induced. We conclude that the form on $L(\chi)$ is negative
on $\sig[(n-r_0),(r_0)]$. 

\noindent\textbf{Assume} $\mathbf{x_{2i-1}<A\le x_{2i}}.$ In this case
we can do the same arguments using $X_o$ and $\sig[(k,n-k),(0)].$ We
omit the details. 
\end{proof}

\subsection{Induction step}\label{7.5}
{The case when the parameter has a single string with
coordinates in an $A_\tau$ with $0<\tau<1/2$ was done in section
\ref{sec:7.4}. }
So we assume there is more than one string. Again we do the case $G$
of type C, and omit the details for the other ones. 

Write the two strings as in (\ref{sec:2.3c}), 
\begin{equation}\label{eq:7.5.1}
(e+\tau_1,\dots ,E+\tau_1),\qquad (f+\tau_2,\dots , F+\tau_2).
\end{equation}
where $0<\tau_1\le 1/2$ and $0<\tau_2\le 1/2.$ 
Recall that because we are in type C, $e,E,f,F\in\bZ,$ and $\ep=0.$

\medskip
We need to show that if $F+f>0$ or $F+f<-2$ when $F+f$ is even, or
$F+f<-1$ when $F+f$ is odd, then the form is negative on a relevant $W-$type.
{Because $\tau_1,\tau_2>0,$ and since r-reducibility and r-unitarity are not affected by small
deformations, we may as well assume that $(f+\tau_2,\dots ,F+\tau_2)$
is the only string with coordinates in $A_{\tau_2},$ and
$(e+\tau_1,\dots ,E+\tau_1)$ the only one with coordinates in $A_{\tau_1}.$}

\medskip

The strategy is as follows. Assume that $L(\chi)$ is r-unitary. We
deform (one of the strings of) $\chi$ to a $\chi_t$ in such a way that the
coresponding induced module is r-irreducible over a finite interval,
but is no longer so at the endpoint, say $t_0.$ Because of the continuity in 
$t$, the module $L(\chi_{t_0})$ is still r-unitary. {The
deformation is such that $L(\chi_{t_0})$ belongs to a larger nilpotent
orbit than $L(\chi)$,} so the induction hypothesis applies, and we get a
contradiction. Sometimes we have to repeat the procedure before we
arrive at a contradiction.

So replace the first string by
\begin{equation}
  \label{eq:7.5.2}
(e+\tau_1+t,\dots ,E+\tau_1+t).
\end{equation}
If $\chi=(e+\tau_1,\dots E+\tau_1;\chi_M),$  then
\begin{align*}
&\chi_t=(e+\tau_1+t,\dots ,E+\tau_1+t;\chi_M), \\
&X(\chi_t):=Ind_M^G[L(e+\tau_1+t,\dots , E+\tau_1+t)\otimes
L(\chi_M)], 
\end{align*}
where  the Lie algebra of $M$ is $fk m=gl(E-e+1)\times \fk g(n-E+e-1)$ 
(so $\vm =gl(E-e+1)\times\vg(n-E+e-1)$). 

\medskip
If $E < |e|,$ we deform  $t$ in the negative direction, otherwise
in the positive direction.  If $t+\tau_1$ reaches $0$ or $1/2,$ before
the nilpotent orbit changes, we should rewrite the string to conform to
the conventions (\ref{eq:2.3c.10}) and (\ref{eq:2.3c.11}). This means
that we rewrite the string as
$(e'+\tau_1',\dots ,E'+\tau_1')$ with $0\le \tau'_1\le 1/2,$ 
and continue the deformation with a
$t$ going in the  direction $t<0$ if $E'<|e'|$, and
$t>0$ if $E'\ge |e'|.$ This is not essential for the argument. 
We may as well assume that the following cases occur.
\begin{enumerate}
\item The nilpotent orbit changes at $t_0=-\tau_1.$
\item the nilpotent orbit does not change, and at $t_0=-\tau_1$,\newline
either   $e,E> x_{2m}+1$ or $-e,-E>x_{2m}+1.$ 
This is the \textit{easy} case when $t$ can be deformed to $\infty$
without any {r-reducibility occuring.} 
\item The nilpotent orbit changes at a $t_0$ such that $0<\tau_1+t_0\le 1/2.$ 
\end{enumerate}
In the first case, the induction hypothesis applies, and since the
string $(f+\tau_2,\dots ,F+\tau_2)$ is unaffected, we conclude
that the signature is negative on a relevant $W-$type. 
{In the second
case we can deform the string so that either $e+\tau_1+t=x_{2m}+1$ or
$E+\tau_1+t=-x_{2m}-1.$ The induction hypothesis applies, and the form
is negative definite on a $W$-relevant type.}
{In the third
case, the only way the nilpotent orbit can change is if the string
$(e+\tau_1+t_0,\dots , E+\tau_1+t_0)$ can be combined with another string to
form a strictly longer string. } 
If $\tau_1+t_0\ne \tau_2,$ the induction
hypothesis applies, and since the string $(f+\tau_2,\dots ,F+\tau_2)$ is
unaffected, the form is negative on a relevant $W-$type. 
If the nilpotent does not change at $t=\tau_2-\tau_1,$
continue the deformation in the same direction. Eventually either (1)
or (2) are satisfied, or else we are in case (3), and the strings in
(\ref{eq:7.5.1}) combine to give a larger nilpotent. There are four
cases:
\begin{equation}
  \label{eq:7.5.3}
  \begin{aligned}
    &(1)\quad e< f\le E\le F,\ e\le f\le E< F\\
    &(2)\quad f\le e\le F<E,\ f< e\le F\le E\\
    &(3)\quad e\le E=f-1< F,\\
    &(4)\quad f\le F=e-1<E. 
  \end{aligned}
\end{equation}

\noindent\textbf{Assume $|e|\le E$}. {Then $t$ is deformed in the
positive direction so $\tau_1<\tau_2.$}
If $e\le 0$, we look at the deformation
(\ref{eq:7.5.2}) for $-\tau_1\le t\le 0.$
If the nilpotent changes for some $-\tau_1<t<0,$
the string $(f+\tau_2,\dots ,F+\tau_2)$ is not involved, the
induction hypothesis applies, so the parameter is not r-unitary.
Otherwise at $t=-\tau_1$  there is one less string with coordinates in
an $A_\tau$ with $\tau\ne 0$, and again the induction hypothesis
applies so the original parameter is not r-unitary.
Thus we are reduced to the case $0<e< E.$ Then consider the nilpotent
orbit for the parameter with $t=-\tau_1+\tau_2.$ 
{
In cases (1), (2) and (3) of (\ref{eq:7.5.3}), the new nilpotent is
larger}, and  one of the strings is 
\begin{equation}
  \label{eq:7.5.4}
  (e+\tau_2,\dots ,F+\tau_2),
\end{equation}
instead of (\ref{eq:7.5.1}), and  $e+F>0.$ The induction hypothesis
applies, so the parameter is not r-unitary, nor is the original one.

In case (4) of (\ref{eq:7.5.3}), the new nilpotent corresponds to the
strings
\begin{equation}
  \label{eq:7.5.7}
  (f+\tau_2,\dots ,E+\tau_2)
\end{equation}
The induction hypothesis applies, so $f+E=0,-2$ if $f+E$ is even or
$f+E=-1$ if it is odd. 
{If this is the case, 
consider a new deformation in \ref{eq:7.5.2}, this time
$-1+\tau_2<t\le 0.$} We may as well assume
that the parameter is r-irreducible in this interval,
or else the argument from before gives the desired conclusion. 
So we arrive at the case when $t=-1 + \tau_2.$ 
The new nilpotent corresponds to the strings
\begin{equation}
  \label{eq:7.5.8}
  (f+\tau_2,\dots, E-1+\tau_2),\quad (e-1+\tau_1),\quad (F+\tau_2).
\end{equation}
Write the parameter as $(\chi';e-1+\tau_2,F+\tau_2).$ Since $e-1=F,$ the
induced module
\begin{equation}
  \label{eq:7.5.9}
  I=Ind_{GL(2)\times G(n-2)}^G[L(e-1+\tau_2,F+\tau_2)\otimes L(\chi')]
\end{equation}
is unitarily induced from a module which is hermitian and
r-irreducible. But the parameter on $GL(2)$ is not unitary unless
$e-1=F=0.$ Furthermore $f+E-1=0, -2$ if $f+E$ is odd or $f+E-1=-1$ if
$f+E$ is even.  
{So the original parameter (\ref{eq:7.5.1}) is
\begin{equation}
  \label{eq:7.5.10}
  \begin{aligned}
&(1+\tau_1,\dots , E+\tau_1),\quad (1-E+\tau_2,\dots ,\tau_2)\quad f+E=0,\\
&(1+\tau_1,\dots ,E+\tau_1),\quad (-E+\tau_2,\dots ,\tau_2)\quad f+E=-2,\\
&(1+\tau_1,\dots ,E+\tau_1),\quad (-E-1+\tau_2,\dots ,\tau_2)\quad f+E=-1.
  \end{aligned}
\end{equation}
Apply the deformation $t+\tau_2$ in the second string with
$-\tau_2<t\le 0.$ We may as well assume that the parameter stays
r-irreducible in this interval.  But then the induction hypothesis
applies at $t=-\tau_2$ because there is one less string with coordinates
in $A_\tau$ with $\tau\ne 0.$ However the first string does no satisfy
the induction hypothesis.}

\noindent\textbf{Assume $|e|>E$}. The same argument applies, but this
time it is $e<E<0$ that requires extra arguments, and in case (3) instead
of case (4) of (\ref{eq:7.5.3}) we have to consider several deformations.  

\subsection{Proof of necessary condition for unitarity in theorem
  \ref{thm:3.1}}\label{sec:7.6}
We first reduce to the case of theorem \ref{t:7.2}. The difference is
that the coordinates in $A_0$ may not form a $\vh/2$ for an even
nilpotent orbit. 
{However 
because of theorem \ref{t:2.6}, and 
properties of relevant $K-$types, r-reducibility and r-unitary are
unaffected by small deformations of the $\chi'_1,\dots ,\chi'_r$
(notation as in (\ref{eq:2.6.3})).}  
So we can deform the strings  corresponding to
$\chi_1',\dots ,\chi_r'$ with coordinates in $A_0$, so that their
coordinates are no longer in $A_0.$ Then the assumptions in  theorem
\ref{t:7.2} are satisfied. 

The argument now proceeds by analyzing each size of strings
separately. In the deformations that we will consider, strings of
different sizes cannot combine so that the nilpotent orbit attached to
the parameter changes.

Fix a size of strings with coordinates not in $A_0.$ If
the strings are not adapted, they can be written in the form
\begin{equation}
  \label{eq:7.6.1}
 (-E-1+\tau_i,\dots ,E+\tau_i)\qquad 0<\tau_i\le 1/2,\ E\equiv\ep( mod\ \bZ). 
\end{equation}
So there is nothing to prove. Now consider a size of strings that are
adapted. Suppose there are \textbf{two} strings of the form
\begin{equation}
  \label{eq:7.6.2}
(-E-1+\tau_i,\dots ,E-1+\tau_i), \quad 0<\tau_i\le 1/2,\ E\equiv\ep (
mod\ \bZ).
\end{equation}
Let $\fk m:=gl(2E+1)\times \fk g(n-2E-1),$ {(recall that $\fk g(a)$
  means a subalgebra/Levi component of the same type as $\fk g$ of
  rank $a$)}   and write
\begin{equation}
  \label{eq:7.6.3}
  \chi:=((-E-1+\tau_i,\dots ,E-1+\tau_i;-E-1+\tau_i,\dots ,E-1+\tau_i);\chi_M).
\end{equation}
 The module
\begin{equation}
  \label{eq:7.6.4}
  Ind_M^G[L(-E-1+\tau_i,\dots ,E-1+\tau_i;-E-1+\tau_i,\dots
  ,E-1+\tau_i)\otimes L(\chi_M)]
\end{equation}
is r-irreducible, and unitarily induced from a hermitian module on
$M$ where the module on $GL(2E+1)$ \textbf{is not unitary}. Thus $L(\chi)$ is
not unitary either. So $L(\chi)$ is unitary only if for each $\tau_i$
there is at most one string of the form $(-E-1+\tau_i,\dots ,E-1+\tau_i).$

Suppose there are two  strings as in (\ref{eq:7.6.1}) with
$\tau_1<\tau_2.$ If there is no string $(-E+\tau_3,\dots , E+\tau_3)$
with $\tau_1<\tau_3<\tau_2,$ then when we deform
$(-E-1+\tau_1+t,\dots ,E-1+\tau_1+t)$ for $0\le t\le \tau_2-\tau_1,$
$X(\chi_t)$ stays r-irreducible. At $t=\tau_2-\tau_1$ we are in case
(\ref{eq:7.6.2}), so the parameter is not unitary. 

On the other hand suppose that there are \textbf{two} strings of the form
\begin{equation}
  \label{eq:7.6.5}
(-E+\tau_i,\dots ,E+\tau_i),\qquad \text{ same } \tau_i.
\end{equation}
Let $\fk m$ be as before, and write 
\begin{equation}
  \label{eq:7.6.6}
  \chi:=((-E+\tau_i,\dots ,E+\tau_i;-E+\tau_i,\dots ,E+\tau_i);\chi_M).
\end{equation}
 The module
\begin{equation}
  \label{eq:7.6.7}
  Ind_M^G[L(-E+\tau_i,\dots ,E+\tau_i;-E+\tau_i,\dots
  ,E+\tau_i)\otimes L(\chi_M)]
\end{equation}
is irreducible, and unitarily induced from a hermitian module on
$M$ where the module on $GL(2E+1)$ \textbf{is unitary}. Thus $L(\chi)$ is
unitary if and only if $L(\chi_M)$ is unitary.  

So we may assume that for each $\tau_i$ there is at most one string of
the form $(-E+\tau_i,\dots ,E+\tau_i).$  

Similarly if there are two strings of the form
$(-E+\tau_1,\dots,E+\tau_1)$ and $(-E+\tau_2,\dots ,E+\tau_2),$ such
that there is no string of the form $(-E-1+\tau_3,\dots, E-1+\tau_3)$
with {$\tau_1<\tau_3<\tau_2$} we
reduce to the case (\ref{eq:7.6.5}).

Let $\tau_k$ be the largest such that a string of the form
$(-E+\tau_k,\dots ,E+\tau_k)$ occurs, and $\tau_{k+1}$ the smallest such
that a string $(-E-1+\tau_{k+1},\dots , E-1+\tau_{k+1})$ occurs. If
$\tau_k>\tau_{k+1},$ we can deform $(-E+\tau_k+t,\dots ,E+\tau_k+t)$ with
$0\le t\le 1-\tau_k-\tau_{k+1}.$ No r-reducibility occurs, and we are
again in case (\ref{eq:7.6.2}). The module is not unitary. If on the other hand
$\tau_k<\tau_{k+1}$, the deformation $(-E-1+\tau_{k+1}+t,\dots,
E-1+\tau_{k+1}-t)$ for $0\le t\le 1-\tau_{k}-\tau_{k+1}$ brings us to
the case (\ref{eq:7.6.5}).

\medskip
Together the above arguments show that conditions (1) and (2) of
theorem \ref{thm:3.1} in types C,D must be satisfied. Remains to check
that for the case of adapted strings, if there is an odd number of a
given size $2E+1,$ then there is a $d_j=2E+1.$  
{This is condition (3)
in theorem \ref{thm:3.1}.}

The arguments above (also the unitarity proof in the case $\vO=(0)$)
show that an $L(\chi)$ is unitary only if it is of the following
form. There is a Levi component $\fk m=gl(a_1)\times \dots\times
gl(a_r)\times \fk g(n-\sum a_i),$  and parameters $\chi_1,\dots
,\chi_r,\chi_0$ such that,
\begin{equation}
  \label{eq:7.6.8}
  L(\chi)=Ind_M^G[\bigotimes L(\chi_i)\otimes L(\chi_0)],
\end{equation}
with the following additional properties:
{
\begin{enumerate}
\item The $\chi_i$ for $i>0$ are as in lemma (1) of
  section \ref{sec:3.2}, with $0<\nu<1/2,$ in particular unitary.
\item $\chi_0$ is such that there is at most one string for every
  $A_\tau$ with $\tau\ne 0,$ and the strings are of different sizes.  
\end{enumerate}
In addition, conditions (1) and (2) of theorem \ref{thm:3.1} are
satisfied for the strings. 
To complete the proof we therefore only need to consider the case of
$L(\chi_0).$ We can deform the parameters of the strings in the
$A_\tau$ with $\tau\ne 0$ to zero without r-reducibility occuring. If
$L(\chi)$ is unitary, then so is the parameter where we 
deform all but one $\tau\ne 0$ to zero. But for a parameter with a
single string belonging to an $A_\tau$ with $\tau\ne 0,$ the
necessary conditions for unitarity are given in section \ref{sec:7.4}.
}


\section{Real nilpotent orbits}\label{sec:8}
\medskip
{In this chapter we review some well known results for
  real nilpotent orbits. Some additional details and references can be
  found in \cite{CM}. The notation in this section differs from the
  previous sections. Fix a real form $\fk g$ of a complex 
semisimple Lie algebra $\fk g_c.$ Let  $\theta_c$ be the
complexification of the Cartan involution $\theta$ of  
$\fk g,$ and write $\ovl{\phantom{x}}$ for the conjugation. Let
{$G_c$} be the adjoint group with Lie algebra $\fk g_c,$ and let
\begin{equation}
  \label{eq:8.1.1}
  \fk g_c=\fk k_c + \fk s_c,\qquad \fk g=\fk k +\fk s
\end{equation}
be the Cartan decomposition. Write $K_c\subset G_c$ for the subgroup
corresponding to $\fk k_c,$ and $G$ and $K$ for the real Lie groups
corresponding to $\fk g$ and $\fk k.$ 

\medskip
We refer to \cite{CM} for standard results about the classification of
$\Ad G_c$ nilpotent orbits in $\fk g_c,$ $\Ad G$ nilpotent orbits on
$\fk g,$ and $\Ad K_c$ nilpotent orbits in $\fk s_c,$ particularly the
Kostant-Sekiguchi correspondence.   

\subsection{}\label{sec:9.1} 
To motivate the need for these results we start with a review of
properties the asymptotic associated support/cycle, and wave front set. 

Let $\pi$ be an admissible $(\fk g_c, K)$ module. We review some facts
from \cite{BV1}. The distribution character $\Theta_\pi$ lifts to an
invariant eigendistribution  $\theta_\pi$ in a neighborhood of the
identity in the Lie algebra. For $f\in C_c^\infty(U),$ where $U\subset\fk
g$ is a small enough neighborhood of $0,$ let 
$f_t(X):=t^{-\dim \fk g_c}f(t^{-1}X).$ Then 
\begin{equation}
  \label{eq:9.1.1}
  \theta_\pi(f_t)=t^{-d}\sum_j c_j\widehat{\mu_{\C O_j(\bR)}}(f) + \sum
  _{i>0}t^{d+i}D_{d+i}(f) ].
\end{equation}
The $D_i$ are homogeneous invariant distributions (each $D_i$ is
tempered and the support of its Fourier transform is contained in the
nilpotent cone). The $\mu_{\CO_j}$ are invariant measures supported on
real forms $\CO_j$ of a single complex orbit $\C O_c,$
and $\mu_{\CO_j(\bR)}$ is the Liouville measure on the nilpotent orbit
associated to the symplectic form induced by the Cartan-Killing form. 
Furthermore $d=\dim_\bC \CO_c/2,$ and the number $c_j$ is called the
multiplicity of $\C O_j(\bR)$ in the leading term of the expansion. 
The closure of the union of the supports of the Fourier transforms of
all the terms occuring in (\ref{eq:9.1.1})  is called the
\textit{asymptotic support}, denoted $AS(\pi).$ The leading term in
(\ref{eq:9.1.1}) will be called $AC(\pi).$ We will use the fact 
that the nilpotent orbits in the leading term are contained in the
\textit{wave front set} of $\theta_\pi$ at the origin, denoted $WF(\pi).$

\medskip
Alternatively, \cite{V3} attaches to each $\pi$ a combination of
$\theta$-stable orbits with integer coefficients 
\begin{equation}
  \label{eq:9.1.2}
  AV(\pi)=\sum a_j \CO_j,
\end{equation}
where $\CO_j$ are nilpotent $K_c-$orbits in $\fk s_c.$ The main result of 
\cite{SV} is that $AC(\pi)$ in (\ref{eq:9.1.1}) and $AV(\pi)$ in
(\ref{eq:9.1.2}) are \textit{the same}. 
Precisely, the leading term in formula
(\ref{eq:9.1.1}), and (\ref{eq:9.1.2}) are the same, when we identify
real and $\theta$ stable nilpotent orbits via the Kostant-Sekiguchi
correspondence. We will use this when we need to compare  $AC(\pi)$
and $AV(\pi)$.

We will need results on how $AC$ and $AV$ behave under real and
cohomological induction. 
\subsection{Complex induction}\label{sec:8.4} 
Let $\fk p_c=\fk m_c +\fk n_c$ be a parabolic subalgebra of $\fk g_c.$  Let
$\fk c_c:=\Ad M_c\cdot e$ be the orbit of a nilpotent element $e\in\fk m_c.$
According to \cite{LS}, 
the induced orbit from $\fk c_c$ is the unique $G_c$ orbit $\fk C_c$
which has the property that $\fk C_c\cap [\fk c_c+\fk n_c]$ is dense (and
open) in {$\fk c_c + \fk n_c.$} 
\begin{proposition}[1]\label{p:cxind}
Let $E=e+n\in e+\fk n_c\subset \fk c_c +\fk n_c$ be a representative
of $\fk C_c.$  
  \begin{enumerate}
  \item $\dim Z_{M_c} (e)=\dim Z_{G_c}(E).$
  \item $\fk C_c \cap [\fk c_c+\fk n_c]$ is a single $P_c$ orbit.
  \end{enumerate}
\end{proposition}
This is theorem 1.3 in \cite{LS}. 
In particular, an element $E'=e'+n'\in \fk c_c +\fk n_c$ is in $\fk
C_c$ if and only if the map
\begin{equation}
  \label{eq:8.4.1}
  \ad E':\fk p_c\longrightarrow T_{e'}\fk c_c +\fk n_c,\qquad \ad E' (y)=[E',y]
\end{equation}
is onto. 

Another characterization of the induced orbit is the following.
\begin{proposition}
The orbit $\fk C_c$ is the unique open orbit in $\Ad G_c(e+\fk
n_c)=\Ad G_c(\fk c_c +\fk n_c),$ as well as in the closure $\ovl{\Ad
  G_c(e+\fk n_c)}=\Ad G_c(\ovl{\fk c_c} +\fk n_c).$
\end{proposition}
We omit the proof, but note that the statements about the closures
follow from the fact that $G_c/P_c$ is compact.

\begin{proposition}
  [2] 
The orbit $\fk C_c$ depends on $\fk c_c\subset \fk m_c,$ but not on
$\fk  n_c.$
\end{proposition}
\begin{proof}
This is proved in section 2 of \cite{LS}. We give a different
proof which generalizes to the real case. 
Let $\xi\in \fk h_c\subset \fk m_c$ be an element in the center of
$\fk m_c$
such that $\langle \xi,\al\rangle\ne 0$ for all
roots $\al\in\Delta(\fk n_c,\fk h_c).$ Then by a standard argument,
\begin{equation}
  \label{eq:8.4.2}
\Ad P_c(\xi + e)=\xi+\fk c_c+\fk n_c.  
\end{equation}
Again because $G_c/P_c$ is compact, 
\begin{equation}
\label{eq:8.4.3}
\ovl{\bigcup_{t>0}\Ad G_c(t\xi+e)}\backslash \bigcup_{t>0}\Ad
G_c(t\xi+e)=\Ad G_c(\ovl{\fk c_c}+\fk n_c).  
\end{equation}
Formula \ref{eq:8.4.3} is valid for any parabolic subgroup with Levi
component $M_c.$ The claim follows because
the left hand side of (\ref{eq:8.4.3}) only depends on $M_c$ and the
orbit $\fk c_c$.
\end{proof}
We write $ind_{\fk l_c}^{\fk g_c}(e)$ or $ind_{\fk p_c}^{\fk g_c}(e)$ 
for the orbit of $E.$
\subsection{Real induction}\label{sec:rhoind}
Let $\fk p=\fk m + \fk n$ be a real
parabolic subalgebra,  $e\in\fk m$  a nilpotent element, and $\fk c:=\Ad M e$. 
\begin{definition}\label{d:8.4}
The $\rho-$induced set from $\fk c$ to $\fk g$ is the finite union
of orbits $\fk C_i:=\Ad G E_i$ such that one of the following
equivalent conditions hold.
\begin{enumerate}
\item $\fk C_i$ is open in  $\Ad G(\fk c+\fk n)$
 and $\ovl{\bigcup\fk C_i}= \ovl{\Ad G(e+\fk n)}$.
\item The intersection $\fk C_i\cap[\fk c+\fk n]$ is open in $\fk c
  +\fk n,$ and the union of the intersections is dense in $\fk c +\fk n.$  
\end{enumerate}
We write
\begin{equation}\label{eq:8.4.4}
ind_{\fk p}^{\fk g} (\fk c)=\bigcup\fk C_i.
\end{equation}
and we say that each $E_i$ is real or $\rho-$induced from $e.$ Sometimes we will write $ind_{\fk p}^{\fk g}(e)$.
\end{definition}
We omit the details of the proof of the equivalence of the two statements.
\begin{proposition}
  \label{p:8.4}
The $\rho-$induced set depends on the orbit $\fk c$ of $e$ and the
Levi component $\fk m,$ but not on $\fk n.$
\end{proposition}
\begin{proof}
The proof is essentially identical to the one in the complex case.
We omit the details.
\end{proof}
The Kostant-Sekiguchi correspondence establishes a 1-1 correspondence
between $\Ad G$ nilpotent orbits in $\fk g$ and $\Ad K_c$ nilpotent
orbits in $\fk s_c.$ Denote by $e\longleftrightarrow \wti e$ this
correspondence.  Then $\rho-$induction is computed in \cite{BB}. This is as
follows. Let $\fk h_c\subset\fk m_c$ be the complexification of a
maximally split real Cartan 
subalgebra $\fk h$, and $\xi\in \C Z(\fk m_c)\cap \fk s_c$ an element
of $\fk h$ such that   
\[
\al\in \Delta(\fk n_c,\fk h_c) \textit{ if and only if } \al(\xi)>0.
\]
Then 
\begin{equation}
  \label{eq:8.4.5}
  \ovl{\bigcup\Ad K_c(\wti E_i)} =\ovl{\bigcup_{t>0} \Ad K_c(t\xi+\wti e)}
\backslash \bigcup_{t>0} \Ad K_c(t\xi+\wti e).
\end{equation}

{
Let $\fk p=\fk m +\fk n$ be a real
parabolic subalgebra. Suppose $AC(\pi)=\sum c_j\CO_{j,\fk m}$, for a
representation $\pi$ of $M.$ Let $v_j\in\CO_{j,\fk m}$ be representatives, and
write $v_{ij}=v_j+X_{ij}$ for representatives of the induced orbits
$\CO_{ij}$ from $\CO_{j,\fk m}.$ Then $\CO_{i,j}\cap(\CO_{j,\fk m}+\fk
n)$ is a finite union of $P-$orbits. Let $v_{ijk}$ be
representatives. Then  
\begin{equation}
  \label{eq:9.1.3}
AC(\Ind_P^{G}(\pi))=\sum_{i,j}c_j
\bigg(\sum_k\bigg|\frac{C_G(v_{ijk})}{C_P(v_{ijk})}\bigg|\bigg)\CO_{ij}. 
\end{equation}
The argument is elementary and in \cite{B4}, except the sum over $k$
is missing. I would like to thank B.L. Harris for pointing out this
gap. I believe that $k=1,$  but I do not know a general
proof. In sections \ref{sec:8.10}, \ref{sec:8.13},
\ref{sec:8.14} and \ref{sec:8.15} we list the $\C O_{ij}$ in the
classical cases. We will show in all cases relevant to our goal that
$k=1$. This is done in the next lemma and corollary.
}

Let $\sig$ be conjugation on $G_c$ so that the fixed points are $G,$
and denote by $\sig$ the corresponding conjugation on $\fk g_c$ as
well. 
{
\begin{lemma}
  \label{l:fiber}
Let $E$ be induced from $e\in\fk m,$ with $\fk p=\fk m+\fk n.$
Assume that every connected component of $C_{G_c}(E)$ intersects $G.$
Then $k=1$ in formula \ref{eq:9.1.3}. 
\end{lemma}
\begin{proof}
In the notation of definition \ref{d:8.4}, let $E=e+X$ be a
representative of one of the induced orbits. 
The claim is proved if we show that if\newline $E'=e+X'\in e+\fk n$ is
conjugate to $E$ by an element in $G,$ then it is conjugate to $E$ by
an element in $P.$ 

By (2) of proposition (1) \ref{p:cxind}, there is $p_c\in P_c$ such
that $\Ad p_cE=\Ad g E=E'$ with $g\in G.$ 
Thus $p_c^{-1}\sig(p_c)\in C_{P_c}(E),$
and $p_c^{-1}g=x\in C_{G_c}(E).$ Since $C_{P_c}(E)\backslash C_{G_c}(E)$ is
finite, the connected component of the identity of $C_{G_c}(E)$ is
contained in $C_{P_c}(E).$ The element $p_c$ can be replaced by any
other element in the coset $p_c C_{P_c}(E).$ The assumption of the
lemma implies then that we can choose $p_c$ so that $x\in C_G(E)$. But
then $p_c=gx^{-1}\in P_c\cap G=P$ as claimed.
\end{proof}
\begin{remark}
The same conclusion and proof hold if we only assume that each coset
$xC_P(E)\subset C_G(E)$ which is fixed by $\sig$, contains a point fixed
by $\sig.$
\end{remark}
}
\medskip

{

\begin{corollary}
  \label{c:fiber}
The assumptions of lemma \ref{l:fiber} hold in the classical cases.
\end{corollary}
\begin{proof}
  Let $\{E,H,F\}$ be a Lie triple such that $\sig(E)=E,$ $\sig(H)=H$ and
  $\sig(F)=F.$ Then each component of $C_G(E)$ and $C_{G_c}(E)$
  respectively, intersects the component group of $C_G(E,H,F)$ and
  $C_{G_c}(E,H,F)$ respectively. In the classical cases these are as
  follows. The inclusion $C_G\subset C_{G_c}$ is the usual one. Recall
  the notation   for nilpotent orbits in \ref{2.1}.

\medskip
\noindent\textbf{Type A.}\ $C_{G_c}(E,H,F)=\prod GL(r_i,\bC)$ and
$C_{G}(E,H,F)=\prod GL(r_i,\bR).$   

\noindent\textbf{Type B,D.}\ For $G_c=O(n,\bC)$, 
$C_{G_c}(E,H,F)$ is a product of $Sp(r_i,\bC)$ for $a_i$ even and
$O(r_i,\bC)$ for $a_i$ odd. For $C_G(E,H,F)$, $Sp(r_i,\bC)$ is
replaced by $Sp(r_i,\bR)$ and $O(r_i,\bC)$ by $O(r_i^+,r_i^-)$
where $r_i^+ +r_i^-=r_i.$ In the classification of nilpotent orbits
referred to in section \ref{sec:8.15}, $r_i^+$ is the number of rows
of size $a_i$ starting with $+,$ $r_i^-$ is the number of rows of size
$a_i$ starting with $-.$ There are minor modifications for $G_c=SO(n,\bC).$

\noindent\textbf{Type C.}\   For $G_c=Sp(2n,\bC)$, 
$C_{G_c}(E,H,F)$ is a product of $Sp(r_i,\bC)$ for $a_i$ odd and
$O(r_i,\bC)$ for $a_i$ even. For $C_G(E,H,F)$, $Sp(r_i,\bC)$ is
replaced by $Sp(r_i,\bR)$ and $O(r_i,\bC)$ by $O(r_i^+,r_i^-)$ with
$r_i^++r_i^-=r_i.$  In the classification of nilpotent orbits referred
to in section \ref{sec:8.14}, $r_i^+$ is the number of rows of size
$a_i$ starting with a $+,$ $r_i^-$ is the number of rows of size $a_i$
starting with a $-.$

\medskip
The claim of the corollary follows from this description.
\end{proof}
}

\subsection{$\theta-$stable induction}\label{sec:thetaind}
{
We will also use an analogue of (\ref{eq:9.1.3}) for $\theta$-stable
induction (derived functors construction or cohomological induction) 
and $AV(\pi)$. 

Let $\fk q_c=\fk l_c + \fk u_c$ be a
$\theta$-stable parabolic subgroup, and write $\ovl{\fk q_c}=\fk
l_c+\ovl{\fk u_c}$ for its complex conjugate. Let $e\in\fk l_c\cap\fk
s_c$ be a nilpotent element. 
\begin{proposition}\label{p:8.5}
 There is a unique $K_c-$orbit $\C O_{K_c}(E)$ so that its intersection
 with $\C O_{L_c\cap K_c}(e)+(\fk u_c\cap \fk s_c)$ is open and
 dense. Similarly there is a unique $G_c$-orbit such that its
 intersection with $\C O_{L_c}(e) +\fk u_c$ is open and dense.
\end{proposition}
\begin{proof}
This follows from the fact that $e+(\fk u_c\cap \fk s_c)$ is formed of
nilpotent orbits, there are a finite number of nilpotent orbits, and
being complex, the $K_c-$orbits have even real dimension. Similarly
for the case of $G_c$-orbits.
\end{proof}
\begin{definition}
  \label{d:thind}
The  unique orbit which
intersects $\C O_{L_c\cap K_c}(e) +(\fk u_c\cap s_c)$  
in a dense open set is called
$\ind_{\fk u_c\cap\fk s_c}^{\fk s_c}(e).$ 

 Similarly write $\ind_{\fk q_c}^{\fk g_c}(e)$ for the unique
 nilpotent $G_c-$orbit which intersects $\C O_{L_c}(e)+\fk u_c$ in an open set.
\end{definition}
\subsection*{Remark} The induced orbit is characterized by the
property that it is the (unique) largest dimensional one  which meets
$e +\fk u_c\cap\fk s_c.$ 
It depends on $e$ as well as $\fk q_c,$ not just $e$ and $\fk l_c.$\qed

We need a special case which is well known, \eg
\cite{Tr1} and \cite{Tr2}. 
Fix a regular dominant integral infinitesimal character $\chi,$ and
let $\C D_\chi$ be the sheaf of twisted differential operators on the
flag variety $\mathscr B:=G_c/B_c.$  The
Beilinson-Bernstein localization functor $\Delta_\chi$ provides an
equivalence of categories between admissible $(\fk g_c,K)$ modules
with infinitesimal character $\chi$, and
$K_c-$equivariant holonomic $\C D_\chi-$modules.   
Let $\fk q_c=\fk l_c +\fk u_c$ be a $\theta-$stable parabolic
subalgebra, $\mathscr P$ be the
generalized flag variety of parabolic subalgebras of type $\fk q_c,$ and
$f:\mathscr B\longrightarrow \mathscr P$ the canonical projection
map. Then $\mathscr R:=K_c\cdot\fk q_c\subset\mathscr P$ is closed, and
let $R$ be the $K_c$-orbit which is open in $f^{-1}(\mathscr{R}).$
Assume $\chi$ is the (dominant integral) infinitesimal character of
the trivial representation. Then $\C D_\chi$ is the sheaf of usual
differential operators $\C D_\mathscr B.$ Let $\C D_{\mathscr P}$ be
the sheaf of differential operators on $\mathscr P.$ Via the
Riemann-Hilbert correspondence, regular $K_c-$equivariant holonomic
sheaves of $\C D_\mathscr B-$modules on $\C B$ are equivalent to
$K_c-$equivariant perverse sheaves (theorem 7.9 in \cite{ABV} and
reference therein). In turn 
irreducible equivariant perverse sheaves correspond to local systems on
$K_c$ orbits. If  $\ovl{\mathscr S}$ is the $\C D_\mathscr P-$module
corresponding to the trivial local system on $\mathscr{R},$ let
$\mathscr S$ be the irreducible module corresponding to the trivial
local system on $R$, and $\pi$ be the irreducible $(\fk g_c,K)$ module
corresponding to $\mathscr S$ by localization. Let
$Ch(\pi):=Ch(\mathscr S)$ be the \textit{characteristic cycle}. A review of its
definition and properties can be found for example in \cite{ABV}
chapter 19. Then $Ch(\ovl{\mathscr{S}})=T^*_{\mathscr{R}}\mathscr P,$
the conormal bundle of $\mathscr{R}.$  By theorem 20.1 of \cite{ABV},
$Ch(\mathscr S)=T^*_R\mathscr B.$   

Let 
$$
\mu:T^*\mathscr B\longrightarrow\fk g^*
$$
be the moment map. Then $\mu(T^*_R\mathscr B)=\ovl{K_c\cdot(\fk
  u_c\cap \fk s_c)}$ is the closure of $ind_{\fk u_c\cap s_c}^{\fk s_c}(0).$ 

\begin{proposition}\label{p:ac}
Assume that $\chi$ is the infinitesimal character of the trivial
representation (dominant integral), and $\mu\mid_{T^*_R\mathscr B}$ is
birational. Then  
$$
\Ad
G_c\ind_{\fk u_c\cap\fk s_c}^{\fk s_c}(0)=\ind_{\fk  q_c}^{\fk g_c}(0),
$$  
and
$$
AV(\pi)=\ind_{\fk u_c\cap\fk s_c}^{\fk s_c}(0).
$$
\end{proposition}
\begin{proof}
The first assertion is a dimension count; recall that for $e\in\fk
s_c,$ $\dim\Ad G_ce=2\dim\Ad K_ce.$ The second assertion follow from 
the results in \cite{Ch}, particularly corollary
2.5.6. Under the assumptions of the proposition, and the notation
in \cite{Ch} preceding 2.5.6,  $F_\nu$ is a single point, and 
$c_{F_\nu}=1.$ 
\end{proof}

Let $\xi$ be a unitary character of the Levi component of a
$\theta-$stable parabolic subalgebra $\fk q_c=\fk l_c+\fk u_c.$
Cohomological induction, detailed \eg in \cite{KnV}, associates to $\xi$ 
a family of $(\fk g_c,K)$ modules $\C R^i_{\fk q_c}(\xi).$ Recall
$S=\dim(\fk u_c\cap\fk k_c)$ from \cite{KnV}.  Extending the
notion of $AC$ or $AV$ by linearity to the Grothendieck group, we can
talk about $AC(\sum (-1)^i \C R^i_{\fk q_c}(\xi))$ or 
$AV(\sum (-1)^i \C R^i_{\fk q_c}(\xi))$. 


\begin{corollary}[1] Assume that the moment map
  $\mu\mid_{T^*_Q\mathscr B}$ is birational. Then 
$$
AV(\sum (-1)^i \C R^i_{\fk q_c}(\xi))=(-1)^S.
$$
In particular if $\xi$ is such that $\C R^i_{\fk q_c}(\xi)=0$ for
$i\ne S,$ then 
$$
AV(\C R^S_{\fk q_c}(\xi))=ind_{\fk q_c\cap\fk
  s_c}^{\fk s_c}(0)
$$
(in particular $\C R^S_{\fk q_c}(\xi) $ is also nonzero). 
\end{corollary}
\begin{proof}
Embed $\pi$ into a coherent family $\pi_\nu$ as in \cite{Ch}. The assumptions
imply that
\begin{equation}
  \label{eq:av}
AC(\pi_\nu)=const\cdot Ch(\mathscr S)\prod_{\al\in\Delta^+(\fk
  l_c)}\langle\nu,\cha\rangle/\prod_{\al\in\Delta^+(\fk
  l_c)}\langle\rho_{\fk l_c},\cha\rangle\ ind_{\fk u_c\cap\fk
  s_c}^{\fk s_c}(0). 
\end{equation}
Setting $\nu=\xi+\rho_{\fk l_c},$ proposition \ref{p:ac} implies that
$const\cdot CH(\mathscr S)=1.$  The proof follows from the fact that
$(-1)^S\pi_\nu=\sum (-1)^i \C R^i_{\fk q_c}(\xi)$ for $\nu=\xi
+\rho_{\fk l_c}$.  
\end{proof}
}

\subsection{u(p,q)}\label{sec:8.7} Let $V$ be a complex finite dimensional
vector space of dimension $n.$ There are two inner classes of
real forms of {$GL(V).$} One is such that $\theta$ is an outer
automorphism. It consists of the real form $GL(n,\bR),$ and when $n$ is
even, also $U^*(n).$ The other one is such that $\theta$ is inner, and
consists of the real forms $U(p,q)$ with $p+q=n.$ In sections
\ref{sec:8.7}-\ref{sec:8.13}, we investigate $\rho$ and $\theta$
induction for the forms $u(p,q)$, and then derive the corresponding results
for $so(p,q)$ and $sp(n,\bR)$ from them in sections
\ref{sec:8.14}-\ref{sec:8.15}. The corresponding results for the other
real forms are easier, the case of $GL(n,\bR)$ is well known, and we
will not need $u(n)^*$. The
usual description of $u(p,q)$ is that $V$ is endowed with a  hermitian
form $(\ ,\ )$ of signature $(p,q),$ and $u(p,q)$ is the
Lie algebra of skew hermitian matrices with respect to this
form. Fix a positive definite hermitian form $\langle\ ,\ \rangle.$ We
will identify the complexification of $\fk g:=u(p,q)$ with $\fk
g_c:=gl(V),$ and the complexification of $U(p,q)$ with $GL(V).$ Up to
conjugacy by $GL(V),$    
\begin{equation}
  \label{eq:8.7.1}
  (v,w)=\langle\theta v,w\rangle,\qquad \theta^2=1,
\end{equation}
The eigenspaces of $\theta$ on $V$ will be denoted $V^\pm$.
The Cartan decomposition is $\fk g_c=\fk k_c +\fk s_c,$ where
{$\fk k_c$
is the $+1$ eigenspace, and $\fk s_c$ the $-1$ eigenspace of $\Ad\theta.$}

\medskip
The classification of nilpotent orbits of $u(p,q)$ is by signed
tableaus as in theorem 9.3.3 of \cite{CM}. The same parametrization
applies to $\theta$ stable orbits under $K_c.$ 
\subsection{}\label{sec:8.8} A parabolic subalgebra of $gl(V)$ is the
stabilizer of a generalized flag
\begin{equation}
  \label{eq:8.8.1}
(0)=W_0\varsubsetneq  W_1\varsubsetneq \dots \varsubsetneq W_k=V.
\end{equation}
Fix complementary spaces $V_i,$
\begin{equation}
  \label{eq:8.8.2}
W_{i}=W_{i-1} + V_i,\qquad i>0.
\end{equation}
They determine a Levi component 
\begin{equation}
  \label{eq:8.8.1.2}
\fk l\cong gl(V_1)\times\dots\times gl(V_k).  
\end{equation}

Conjugacy classes under $K_c$ of $\theta$-stable parabolic subalgebras
are \newline parametrized by ordered pairs $(p_1,q_1),\dots ,(p_k,q_k)$ such
that the sum of the $p_i$ is $p,$ and the sum of the $q_i$ is $q$. A
realization in terms of flags is as follows.  Choose the
$W_i$ to be stable under $\theta,$ or equivalently that
the restriction of the hermitian form to each $W_i$ is 
nondegenerate. In this case we may assume that the $V_i$ are
$\theta$-stable as well, and let $\fk q_c=\fk l_c + \fk u_c$ be the
corresponding  parabolic subalgebra of $gl(V).$ The signature 
of the form restricted to $V_i$ is $(p_i,q_i),$ so that 
\begin{equation}
  \label{eq:8.9.1}
\fk l_c\cap \fk g\cong u(p_1,q_1)\times\dots\times u(p_k,q_k).
\end{equation}

The following algorithm for computing the induced
orbit in the case $\fk g\cong u(p,q)$ and a maximal parabolic
subalgebra holds. 

\medskip
{
\textit
{Suppose the signature of $V_1$ is $(a_+,a_-).$ Then add $a_+$
$+$'s to the beginning of largest possible rows of $e$ starting with  $-$
and $a_-$ $-$'s to the largest possible rows of $e$ starting with 
$+.$ If $a_+$ is larger than the number of rows starting with $-,$ add a new
row of size $1$ starting with $+.$ The similar rule applies to $a_-.$ }
}

If $e\in gl(V_1),$ the analogous procedure applies, but the $a_+$
$+$'s are added at the end of the largest possible rows finishing in
$-$ and $a_-$ $-$'s to the end of the largest possible rows finishing in $+.$

Because induction is transitive, the above algorithm can be
generalized to compute the $\theta$-induced of any nilpotent orbit. We
omit the details. See \cite{Tr1} and \cite{Tr2}.

\subsection{}\label{sec:8.10} 

Conjugacy classes under $G$ of real parabolic
subalgebras are given by ordered subsequences $n_1,\dots ,n_k$ and a
pair $(p_0,q_0)$ such that $\sum n_i+p_0=p$ and $\sum n_i+q_0=q.$ The
complexification of the corresponding \textit{real} parabolic
subalgebra is given as follows. Start with a partial flag  
\begin{equation}
  \label{eq:8.10.1}
(0)=W_0\varsubsetneq \dots\varsubsetneq W_k  
\end{equation}
such that the hermitian form is trivial when restricted to $W_k,$ and
complete it to 
\begin{equation}
  \label{eq:8.10.2}
 (0)=W_0\varsubsetneq \dots\varsubsetneq W_k\varsubsetneq
 W_k^*\varsubsetneq\dots\varsubsetneq W_0^*=V  
\end{equation}
Choose transverse spaces
\begin{equation}
  \label{eq:8.10.3}
  W_i=W_{i-1}+ V_i,\qquad W_i^*=W_{i-1}^* + V_i^*,\qquad  W_k^*=W_k + V_0.
\end{equation}
They determine a Levi component 
\begin{equation}
  \label{eq:8.8.10.4}
\fk l_c = gl(V_1)\times\dots\times
gl(V_{k})\times gl(V_0)\times gl(V_{k}^*)\times\dots\times gl(V_1^*),
\end{equation}
so that
\begin{equation}
  \label{eq:8.10.5}
  \fk l_c\cap \fk g = gl(V_1,\bC)\times\dots\times gl(V_{k},\bC)\times
  u(p_0,q_0).
\end{equation}
Then $n_i=\dim V_i,$ and $(p_0,q_0)$ is the signature of $V_0.$

\subsection{}\label{sec:8.13} Suppose $\fk p_c=\fk m_c +\fk n_c$ is
the complexification of a 
real parabolic subalgebra corresponding to the flag $(0)\subset
V_1\subset V_1+V_0\subset V_1 +V_0 + V_1^*,$ and let $e\subset
gl(V_0)$ be a real nilpotent element. The rest of the notation is as
in section \ref{sec:8.4}. 
\begin{theorem}\label{t:8.13}
The tableau of an orbit $\Ad G(E_i)$ which is in the $\rho-$induced
set $ind_{\fk p}^{\fk g}(\fk c),$  is obtained from
the tableau of $e$ as follows. 

Add two boxes  to the end of each of $\dim V_1$ of the largest rows
such that the result is still a signed tableau. 
\end{theorem}
\begin{proof}
We use (\ref{eq:8.4.2}) and (\ref{eq:8.4.3}). Let
$\al\in\Hom[V_1,V_1^*]\oplus\Hom[V_1^*,V_1]$ be nondegenerate such
that $\al^2=Id\oplus Id,$ and extend it to an endomorphism
$\xi\in gl(V)$ so that its restriction to  $V_0$ is zero.  This is an
element such that the centralizer of $\ad\xi$ is $\fk m,$ in particular,
$[\xi,e]=0.$ Let 
\begin{equation}
  \label{eq:8.13.1}
  P(X)=X^m+a_{m-1}X^{m-1}+\dots + a_0
\end{equation}
be any polynomial in $X\in gl(V)$. Suppose $t_i\in\bR$ are such that
$t_i\to 0,$ and assume there are $g_i\in K$ such that $t_ig_i(\xi
+e)g_i^{-1}\to E.$ Then 
\begin{equation}
  \label{eq:8.13.2}
  \ker t_i^mP(g_i(\xi+e)g_i^{-1})\cong\ker P(\xi+e).
\end{equation}
On the other hand,
\begin{equation}
  \label{eq:8.13.3}
  \begin{aligned}
&t_i^mP(g_i(\xi +e)g_i^{-1})=
[t_ig_i(\xi+e)g_i^{-1}]^m + \\
&+a_{m-1}t_i[t_ig_i(\xi+e)g_i^{-1}]^{m-1}+\dots +
t_i^m Id \to E^m,
\end{aligned}
\end{equation}
as $t_i\to 0.$ Thus
\begin{equation}
  \label{eq:8.13.4}
  \dim\ker E^m\mid_{V_\pm}\ge \dim\ker P(\xi+e)\mid_{V_\pm}.
\end{equation}
Choosing $P(X)=(X^2-1)X^n,$ we conclude that $E$ must be
nilpotent. Choosing $P(X)=X^m,\ (X\pm 1)X^{m-1}$ or
$P(X)=(X^2-1)X^{m-2},$ we can bound the dimensions of $\ker
E^m\mid_{V_\pm}$ to conclude that it must be in the closure of one of
the nilpotent orbits given by the algorithm of the theorem. 
The fact that these nilpotent
orbits are in (\ref{eq:8.4.3}) follows by a direct calculation which we omit.
\end{proof}

\subsection{sp(V)}\label{sec:8.14} Suppose $\fk g_c\cong sp(V_0),$ where
$(V_0,\langle\ ,\ \rangle)$ is a real symplectic vector space of
dimension $n.$  The complexification $(V,\langle\ ,\ \rangle)$ admits a complex
conjugation $\ovl{\phantom{x}},$ and we define a nondegenerate hermitian  form 
  \begin{equation}
    \label{eq:8.14.1}
    (v,w):=\langle v,\ovl{w}\rangle
  \end{equation}
which is of signature $(n,n).$ Denote by $u(n,n)$
the corresponding unitary group. Since $sp(V_0)$
stabilizes $(\ ,\ ),$ it embeds in $u(n,n),$ and the Cartan
involutions are compatible. 
The classification of nilpotent orbits is as follows. See chapter 9 of
\cite{CM} for a more detailed explanation. 

\noindent\textit{
Each orbit corresponds to a signed tableau so that every odd part
occurs an even number of times. Odd sized rows occur in pairs, one
starting with $+$ the other with $-.$ Even sized rows start with a $+$
or a $-$.
}

A real parabolic subalgebra of $sp(V)$ is the stabilizer of a flag of
isotropic subspaces
\begin{equation}
  \label{eq:8.14.2}
  (0)=\C W_0\subset\dots\subset \C W_k,
\end{equation}
so that the symplectic form restricts to  0 on $\C W_k.$ 
As before, complete this to a flag
\begin{equation}
  \label{eq:8.14.3}
 (0)=\C W_0\subset\dots\subset \C W_k\subset \C
    W_k^*\subset\dots\subset\C W_0^*=V.  
\end{equation}
We choose transverse spaces
\begin{equation}
  \label{eq:8.14.4}
\C W_{i}=\C W_{i-1} + V_i,\quad \C W_k^*=\C W_k + \C W,\quad 
\C W_{i-1}^*=\C W_{i}^* + V_i^*
\end{equation}
in order to fix a Levi component.  We get
\begin{equation}
  \label{eq:8.14.5}
  \fk l\cong gl(V_1)\times\dots\times gl(V_k)\times sp(\C W).
\end{equation}
If we assume that $V_i,\ \C W$ are $\theta$-stable, then the
corresponding parabolic subalgebra is $\theta$-stable as well, and the
real points of the Levi component are
\begin{equation}
  \label{eq:8.14.6}
\fk l_0\cong u(p_1,q_1)\times\dots\times u(p_k,q_k)\times sp(\C W_0).
\end{equation}
where $(p_i,q_i)$ is the signature of $V_i.$ The Levi component of the 
parabolic subalgebra corresponding to (\ref{eq:8.14.4}) in $gl(V)$
satisfies
\begin{equation}
  \label{eq:8.14.7}
  \fk l'\cong u(p_1,q_1)\times\dots\times u(p_k,q_k)\times
  u(n_0,n_0)\times u(q_k,p_k)\times\dots\times u(q_1,p_1).
\end{equation}

For a maximal $\theta$-stable parabolic subalgebra, 
the Levi component $\fk l$ satisfies $\fk l\cong u(p_1,q_1)\times
sp(\C W_0).$ Let
$e\in sp(W)$ be a $\theta$-stable nilpotent element. The algorithm for
induced nilpotent orbits in section \ref{sec:8.8} implies the
following algorithm for {$ind_{\fk q_c}^{\fk g_c} (e).$}

\begin{enumerate}
\item add $p$ boxes labelled $+$'s to the beginning of the longest rows
  starting with $-$'s, and $q$ $-$'s to the beginning of the longest
  rows starting with $+$'s.
\item add $q$ $+$'s to the ending of the longest possible rows
  starting with $-$'s, and $p$ $-$'s to the beginning of the longest
  possible rows starting with $+$'s.
\end{enumerate}

Unlike in the complex case, the result is automatically the signed
tableau corresponding to a nilpotent element in $sp(V).$ Again see
\cite{Tr1} and \cite{Tr2} for details.

\medskip
For a maximal real parabolic subalgebra, we must assume that
 $\ovl{V_1}=V_1,\ \ovl{\C W}=\C W.$ Let $V_{1,0}$ and $\C W_0$ be
 their real points. The Levi component satisfies
 \begin{equation}
   \label{eq:8.14.8}
   \fk l\cong gl(V_{1,0})\times sp(\C W_0).
 \end{equation}
The results in section \ref{sec:8.13} imply the
 following algorithm for real induction.

 \begin{enumerate}
 \item add two boxes  to the largest $\dim V_1$  rows of $e$ so that
   the result is still a signed tableau for a nilpotent orbit.
\item Suppose $\dim V_1$ is odd and the last row that would be increased by 2
  is odd size as well. In this case  there is a pair of rows of this
  size, one starting with $+$ the other with $-.$ In this case
  increase these two rows by one each.
 \end{enumerate}
\subsection{so(p,q)}\label{sec:8.15} Suppose $\fk g_c\cong so(V_0),$ where
 $(V_0,\langle\ ,\ \rangle)$  is a real nondegenerate quadratic space
 of signature $(p,q).$ The complexification admits a hermitian form
 $\langle\ , \ \rangle$ with signature $(p,q)$ as well as a complex
 nondegenerate quadratic form $(\ ,\ ),$ which restrict to $\langle\
 ,\ \rangle$  on $V_0.$  The form $\langle\ , \
 \rangle$ gives an embedding of $o(p,q)$ into $u(p,q)$ with compatible
 Cartan involutions. The classification of nilpotent
orbits of $so(V_0)$ or equivalently $\theta$-stable nilpotent
orbits is as follows. See chapter 9 of \cite{CM} for more details.

\noindent\textit{
Orbits correspond to signed tableaus so that every even part
occurs an even number of times. 
Even sized rows occur in pairs, one starting with $+$ the other with
$-.$ An odd sized row starts with a $+$ or a $-.$ 
When all the rows have even sizes, there are two nilpotent
orbits denoted I and II.
}

A parabolic subalgebra of $so(V)$ is the stabilizer of a flag of
isotropic subspaces
\begin{equation}
  \label{eq:8.15.1}
  (0)=\C W_0\subset\dots\subset \C W_k,
\end{equation}
so that the quadratic form restricts to  0 on $\C W_k.$ 
As before, complete this to a flag
\begin{equation}
  \label{eq:8.15.2}
 (0)=\C W_0\subset\dots\subset \C W_k\subset \C
    W_k^*\subset\dots\subset\C W_0^*=V.  
\end{equation}
We choose transverse spaces
\begin{equation}
  \label{eq:8.15.3}
\C W_{i}=\C W_{i-1} + V_i,\quad \C W_k^*=\C W_k + \C W,\quad 
\C W_{i-1}^*=\C W_{i}^* + V_i^*
\end{equation}
in order to fix a Levi component,
\begin{equation}
  \label{eq:8.15.4}
 \fk l\cong gl(V_1)\times\dots\times gl(V_k)\times so(W).
\end{equation}
To get a $\theta$-stable parabolic subalgebra we must assume $V_i,\ W$
are $\theta$-stable and so $\ovl{V_i}=V_i^*,\  \ovl{W}=W.$ If the
signature of $V_i$ with respect to $\langle\ , \ \rangle$ is
$(p_i,q_i),$ and that of $W$ is $(p_0,q_0),$ then
\begin{equation}
  \label{eq:8.15.5}
\fk l_0\cong u(p_1,q_1)\times\dots\times u(p_k,q_k)\times so(p_0,q_0).
\end{equation}
The  parabolic subalgebra corresponding to (\ref{eq:8.15.2}) in $gl(V)$
satisfies
\begin{equation}
  \label{eq:8.15.6}
  \fk l'\cong u(p_1,q_1)\times\dots\times u(p_k,q_k)\times
  u(p_0,q_0)\times u(p_k,q_k)\times\dots\times u(p_1,q_1).
\end{equation}
For a maximal $\theta$-stable parabolic subalgebra, 
the Levi component $\fk l$ satisfies $\fk l\cong u(p_1,q_1)\times
so(\C W_0).$ Let
$e\in so(W)$ be a $\theta$-stable nilpotent element. The algorithm for
induced nilpotent orbits in section \ref{sec:8.8} implies the
following algorithm for $ind_{\fk q_c}^{\fk g_c} (e).$
\begin{enumerate}
\item add $p_1$ $+$'s to the beginning of the longest possible rows
  starting with $-$'s, and $q_1$ $-$'s to the beginning of the longest
  possible rows starting with $+$'s.
\item add $p_1$ $+$'s to the ending of the longest possible rows
  starting with $-$'s, and $q_1$ $-$'s to the beginning of the longest
  possible rows starting with $+$'s.
\end{enumerate}
Unlike in the complex case, the result is automatically a signed tableau
for a nilpotent element in $so(V).$

\medskip
For a maximal real parabolic subalgebra, we must assume that
 $\ovl{V_1}=V_1,\ \ovl{\C W}=\C W.$ Let $V_{1,0}$ and $\C W_0$ be
 their real points. The Levi component satisfies
 \begin{equation}
   \label{eq:8.15.7}
   \fk l\cong gl(V_{1,0})\times so(W_0).
 \end{equation}
The results in section \ref{sec:8.13} imply the
 following algorithm for real induction.

\medskip
 \noindent\textit{
 Add two boxes to $\dim V_1$ of the largest possible rows so that the
 result is still a signed tableau for a nilpotent orbit.
Suppose $\dim V_1$ is even and the last row that would be increased by 2
is even size as well. In this case  there is a pair of rows of this
size, one starting with $+$ the other with $-.$ Increase
these two rows by one each so that the result is still a signed tableau.
When there are only even sized rows and $\dim V_1$ is even as
well, type I goes to type I and type II goes to type II.
}

\section{Unitarity}\label{sec:9}

\medskip
In this section we prove the unitarity of the representations of
the form $L(\chi)$ where $\chi=\vh/2.$
As already mentioned, 
in the $p-$adic case this is done in \cite{BM1}. It amounts to the observation
that the Iwahori-Matsumoto involution preserves unitarity, and takes
such an $L(\chi)$ into a tempered representation. 

\medskip
The idea of the proof in the real case is described in  \cite{B2}.  
We will do an  induction on rank.
We rely heavily on the properties of the
wave front set, asymptotic support and associated variety, and their
relations to primitive ideal cells and Harish-Chandra cells. 
We review the needed facts in sections \ref{sec:9.1}-\ref{sec:9.4}.
Details are in \cite{BV1}, \cite{BV2}, and \cite{B3}. In
sections \ref{sec:9.5}-\ref{sec:9.7} we give details of the proof of
the unitarity in the case of $SO(2n+1).$ The proof is simpler than in \cite{B2}.

\bigskip
\subsection{}\label{sec:9.2}
We will make heavy use of the results in \cite{V2}. 
Fix a regular integral  infinitesimal character $\chi_{reg}.$ Let 
$\C G^G(\chi_{reg})$ be the set of parameters of irreducible admissible 
$(\fk g_c,K)$ modules with infinitesimal character
$\chi_{reg}.$ We will suppress the superscript $G$ whenever there is
no danger of confusion. The explicit description of $\C
G^G(\chi_{reg})$ is well known, called the Langlands or Vogan
classification. We will use the description in \cite{V2}.   
Denote by $\bZ\C G(\chi_{reg})$ the corresponding
Grothendieck group of characters. 
Recall from \cite{V2} (and references therein)  
that there is an  action of the Weyl group
on $\bZ\C G(\chi_{reg}),$ called the {\it coherent continuation action.} 
As a set, $\C G(\chi_{reg})$ decomposes into a disjoint union
of blocks  $\C B$, and so $\bZ\C G$ decomposes into a direct sum  
\begin{equation}
  \label{eq:9.2.1}
\bZ\C G(\chi_{reg})=\bigoplus \bZ\C G_{\C B}(\chi_{reg}).  
\end{equation}
Each $\bZ\C G_\C B(\chi_{reg})$ is preserved by the coherent continuation
action. We give the explicit description of the $\bZ\C G_\C B$  in all
classical cases. Many of the details, and complementary explanations
can be found in \cite{McG}.

\noindent\textbf{Type B:\ } In order to conform
to the duality between type B and type C in \cite{V2}, we only consider
the real forms with $p>q.$ The Weyl group acts on the linear space
$\sum\bZ\C G^{SO(p,q)}(\chi_{reg})$, and the representation is
\begin{equation}
\label{eq:9.2.2}
\begin{aligned}
&\sum_{p>q}\bZ\C G^{SO(p,q)}(\chi_{reg})=\\
&=\sum_{a,b,\tau}Ind_{W_a\times
  W_b\times W_{2s}\times S_t}^{W_n} 
[sgn\otimes sgn\otimes\sig[\tau,\tau]\otimes triv],
\end{aligned}
\end{equation}
where $\tau$ is a partition of $s,$ and $a+b+2s+t=n.$ 
The multiplicity of a $\sig[\tau_L,\tau_R]$ in one of the induced
modules in (\ref{eq:9.2.2}) is as follows. Choose a $\tau$ that fits
inside both $\tau_L$ and $\tau_R,$ and label it by $\bullet$'s. 
{Add ``a'' $r$ and ``b''
$r'$ to $\tau_R,$ at most one to each row for inducing from the sign
representations of $W_a$ and $W_b$. Add ``t''
$c$, at most one to each column, to  $\tau_L$ or
$\tau_R$ for inducing from the trivial representation on $S_t.$} The
multiplicity of $\sig$ in the induced module for a given $(\tau,a,b)$
is then the number of ways that $\tau_L,\ \tau_R$ can be filled in 
this way. This procedure uses induction in stages, and the well known formula
\begin{equation}
Ind_{S_n}^{W_n}(triv)=\sum_{k+l=n} \sig[(k),(l)].
\label{eq:9.2.3}\end{equation}

\subsection*{Example} Let $\fk g_c=so(5).$ The real
forms are $so(3,2),so(4,1), so(5).$ The choices of $(\tau,a,b,t)$ are

  \begin{align}  \label{eq:9.2.4}
&  (1,1,0,0),&& (1,0,1,0), && (1,0,0,1), && && &&\\
& (0,2,0,0), && (0,1,1,0)),&& (0,1,0,1), && (0,0,2,0), && (0,0,1,1), 
&&(0,0,0,2). \notag    
  \end{align}
Let $\sig=\sig[(1),(1)].$ Then its multiplicity is given by the number
of labelings
\begin{equation}
  \label{eq:9.2.5}
  \begin{aligned}
& (\bullet,\bullet) &&\emptyset, &&\emptyset, && && &&\\
& \emptyset, && \emptyset, &&(c,r), &&\emptyset, &&(c,r'), &&
(c,c).    
  \end{aligned}
 \end{equation}
For $\sig=\sig[(0),(2)]$ we get
\begin{equation}
  \label{eq:9.2.6}
  \begin{aligned}
&\emptyset, &&\emptyset, &&\emptyset &&         &&       &&\\
&\emptyset,   && (0,rr'),   &&(0,rc),  &&\emptyset, &&(0,r'c), &&(0,cc).    
  \end{aligned}
 \end{equation}
\qed

The following formula sorts the representations according to the
various real forms $SO(p,q)$ with $p+q=2n+1.$  Each real form gives
a single block. A representation
occuring in $\C G,$ labelled as above, occurs in $\bZ\C G^{SO(p,q)}$ with
\begin{equation}
  \label{eq:9.2.7}
 p=n+1 + \mid \# r' -\#r \mid - \ep, \text{ where } \ep=\begin{cases} 
 0 &\text{ if } \# r' \ge  \# r,\\
            1 &\text{ otherwise.}
\end{cases} 
\end{equation}
In the above example, $(\bullet,\bullet)$, $(c,c)$, $(0,rr'),$
$(0,rc)$ and $(0,cc)$ belong to $so(3,2)$ while $(c,r')$ and $(0,r'c)$
belong to $so(4,1).$ 

\medskip
To each pair of partitions parametrizing a representation of $W,$ 
\begin{equation}
  \label{eq:9.2.8}
  \tau_L=(r_0,\dots ,r_{2m}),\qquad\tau_R=(r_1,\dots ,r_{2m-1}),\qquad
  r_i\le r_{i+2},
\end{equation}  
Lusztig attaches a \textit{symbol} 
\begin{equation}
\begin{pmatrix} r_0&   &r_2+1&     &      &\dots       &            &r_{2m}+m\\
            &r_1&     &r_3+1&\dots &            &r_{2m-1}+m-1&       
\end{pmatrix}.
\label{eq:9.2.9}
\end{equation} 
The symbol  is called special if 
\begin{equation}
  \label{eq:9.2.10}
  r_0\le r_1\le r_2+1\le r_3+1\le \dots\le r_{2m}+m. 
\end{equation}
Two representations belong to the same double cell if and only if
their symbols have the same entries. Given a special symbol of the form
(\ref{eq:9.2.9}), the corresponding nilpotent orbit $\CO_c$ has
partition obtained as follows. Form the set 
\begin{equation}
  \label{eq:9.2.11}
\{2r_{2i}+2i+1, 2r_{2j-1}+2j-2{: 0\le i\le m, 1\le j\le m}\},
\end{equation}
and order the numbers in increasing order, $x_0\le \dots \le x_{2m}.$ The
partition of $\CO_c$ is
\begin{equation}
  \label{eq:9.2.12}
  (x_0,x_1-1,\dots ,x_i-i,\dots ,x_{2m}-2m).
\end{equation}

\medskip
\noindent\textbf{Type C:\ } The representation $\C G(\chi_{reg})$ is
obtained from the one in type B by tensoring with sign. Thus
\begin{equation}
\bZ\C G(\chi_{reg})=\sum_{a,b,\tau}Ind_{S_t\times W_{2s}\times W_a\times W_b}
^{W_n}[sgn\otimes\sig[\tau,\tau]\otimes triv\otimes triv],
\label{eq:9.2.13}\end{equation}
where $\tau$ is a partition of $s,$ and $a+b+2s+t=n.$ This takes into
account the duality in \cite{V2} of 
types $B$  and $C$. We write $r$ for the sign representation of
$S_t,$ and $c$ and $c'$ for the trivial representations of $W_a,\ W_b.$
A representation of $W$ is parametrized by a pair of partitions
$(\tau_L,\tau_R),$ with
\begin{equation}
  \label{eq:9.2.14}
  \tau_L=(r_0,\dots ,r_{2m}),\qquad\tau_R=(r_1,\dots ,r_{2m-1}),\qquad
  r_i\le r_{i+2}.
\end{equation}
The associated symbol is 
\begin{equation}
\begin{pmatrix} r_0&   &r_2+1&     &      &\dots       &            &r_{2m}+m\\
            &r_1&     &r_3+1&\dots &            &r_{2m-1}+m-1&       
\end{pmatrix},
\label{eq:9.2.15}\end{equation} 
and it is called special if 
\begin{equation}
  \label{eq:9.2.16}
  r_0\le r_1\le r_2+1\le r_3+1\le \dots\le r_{2m}+m. 
\end{equation}
Two representations belong to the same double cell if their symbols
have the same entries. Given a special symbol as in (\ref{eq:9.2.15}),
the nilpotent orbit $\CO_c$ attached to the double cell has partition
obtained as follows. Order the set
\begin{equation}
  \label{eq:9.2.17}
 \{2r_{2i}+2i, 2r_{2j-1} + 2j-1 {: 0\le i\le m, 1\le j\le m}\}
\end{equation}
in increasing order, $x_0\le \dots \le x_{2m}.$ Then the partition of
$\CO_c$ is
\begin{equation}
  \label{eq:9.2.18}
(x_0,\dots ,x_j-j,\dots ,x_{2m}-2m ).
\end{equation}
The decomposition into blocks is obtained from the one for type B by
tensoring with $sgn.$

\medskip
\noindent\textbf{Type D:\ } Since in this case 
$\sig[\tau_L,\tau_R]$ and $\sig[\tau_R,\tau_L]$ parametrize the same
representation, (except of course when $\tau_L=\tau_R$ which
corresponds to two nonisomorphic representations), we assume that the
size of $\tau_L$ is the larger one. The Cartan subgroups are parametrized by
integers $(t,u,2s,p,q),$  $p+q+2s+t+u=n.$ There are
actually two Cartan subgroups for each $s>0,$ related by the outer
automorphism of order 2. Then 
\begin{align}
\label{eq:9.2.19}
&\sum \bZ\C G^{SO(p,q)}(\chi_{reg})=\\
&=\sum
Ind_{W_a\times W_b\times W'_{2s}\times W_t\times W_u}^{W'_n}
[sgn\otimes sgn \otimes\sig[\tau,\tau]_{I,II}\otimes triv\otimes triv].
\notag  \end{align}
where $a+b+2s+t+u=n.$ 
The sum is also over $\tau$ which is a partition of $s$. 
We label the $\sig$ by $\bullet$'s, trivial representations by $c$
and $c'$ and the $sgn$ representations by $r$ and $r'.$ These are
added to  $\tau_L$  when inducing. In this case the sum on the left is
over the real forms $SO(p,q)$ with all $p+q=2n.$ Each real form gives rise to
a single block. A representation labelled as above belongs to the block
$\bZ\C G^{SO(p,q)}$ with $p=n+\# r' - \# r.$ 

If 
\begin{equation}
  \label{eq:9.2.20}
  \tau_L=(r_0,\dots ,r_{2m-2}),\qquad \tau_R=(r_1,\dots ,r_{2m-1}),
\end{equation}
then the associated symbol is
\begin{equation}\label{eq:9.2.21}
\begin{pmatrix} 
r_0&r_2+1&\dots &r_{2m-2}+m-1\\
r_1&r_3+1&\dots &r_{2m-1}+m-1
\end{pmatrix}. 
\end{equation}
A representation is called special if the symbol satisfies
\begin{equation}
  \label{eq:9.2.22}
  r_0\le r_1\le r_2+1\le r_3+1\le\dots \le r_{2m-1}+m-1.
\end{equation}
Two representations belong to the same double cell if their symbols
have the same entries. The nilpotent orbit $\CO_c$ attached to the
special symbol is given by the same procedure as for type B.


\subsection{}\label{sec:9.3} We follow section 14 of \cite{V2}. 
We say that $\pi'\le\pi$ if $\pi'$ is a factor of
$\pi\otimes F$ with $F$ a finite dimensional representation with
highest weight equal to an integer sum of roots. 
Two irreducible representations $\pi,\ \pi'$ are said to be in the same
Harish-Chandra cell if $\pi'\le\pi$ and $\pi\le\pi'.$ The
Harish-Chandra cell of $\pi$ is denoted $\C C(\pi)$. 

Recall the
relation $\lr$ from definition 14.6 of \cite{V2}. The \textit{cone
  above} $\pi$ is defined to be
\begin{equation}
  \label{eq:9.3.1}
\ovl{\C C}^{LR}(\pi):=\{\pi'\ :\ \pi'\lr\pi \}.  
\end{equation}
The subspace in $\bZ\C G$ generated by the elements in $\ovl{\C
  C}^{LR}(\pi)$ is a representation of $W$ denoted $\ovl{\C
  V}^{LR}(\pi).$ The equivalence $\elr$  is defined by
\begin{equation}
  \label{eq:9.3.2}
  \gamma\elr\phi \Longleftrightarrow \gamma\lr\phi\lr\gamma.
\end{equation}
The Harish-Chandra cell $\C C(\pi)$ is then
\begin{equation}
  \label{eq:9.3.3}
  \C C(\pi)=\{\pi'\ :\ \pi'\elr\pi\}.
\end{equation}
{Define 
\begin{equation}
  \label{eq:9.3.5}
  \C C^{LR}_+(\pi)=\ovl{\C C}^{LR}(\pi)\  \backslash\  \C C(\pi)
\end{equation}
and $\C V(\pi)$ and $\ovl{\C V}^{LR}(\pi)_+$ in analogy with $\ovl{\C
  V}^{LR}(\pi).$ Thus there is a representation of $W$ on $\C V(\pi)$
by the natural isomorphism
\begin{equation}
  \label{eq:9.3.6}
  \C V(\pi)\cong \ovl{V }^{LR}(\pi)/\ovl{\C V}^{LR}(\pi)_+.
\end{equation}
}
Let $\C O_c\subset \fk g_c$ be a nilpotent orbit. We say
that a Harish-Chandra cell {$\C C(\pi)$} is attached to
a complex orbit $\C O_c$ if 
$$
\ovl{\Ad G_c (AS(\pi))}=\ovl{\C O_c}.
$$ 
The sum of the Harish Chandra cells attached to $\C O_c$ is denoted
$\C V(\C O_c).$

{
Let $\fk h_a \subset \fk g_c$ be an abstract Cartan subalgebra and let
$\Pi_a$ be a set of (abstract) simple roots.  Let $\gamma$ be a
Langlands parameter of an irreducible representation $\C L(\gamma)$ with
infinitesimal character the same as the trivial representation (or
some fixed finite dimensional representation). 
Denote by $\tau(\gamma)$ the
$\tau-$invariant as defined in \cite{V2}.} Given a block $\C B$ and  disjoint
orthogonal sets $S_1,\ S_2\subset \Pi_a,$ define 
\begin{equation}
  \label{eq:9.2.23}
 \C B(S_1,S_2)=\{\gamma \in \C B\ |\ S_1\subset \tau(\gamma),\ S_2\cap
\tau(\gamma)=\emptyset \}\enspace. 
\end{equation}
If in addition we are given a nilpotent orbit $\C O_c\subset \fk g_c,$ we
can also define
\begin{equation}
  \label{eq:9.2.24}
 \C B(S_1,S_2,\C O_c)=\{\gamma \in \C B(S_1,S_2) |\ AS(\C L(\gamma))
\subset \ov {\C O_c} \}\enspace .  
\end{equation}
{If an $S_i$ is absent from the notation, assume
  $S_i=\emptyset.$}

\medskip
Let $W_i=W(S_i)$, and define
\begin{equation}
\begin{aligned}
m_S(\sigma)&=[\sigma :Ind^{W}_{W_1\times W_2}(Sgn\otimes Triv)],\\
m_{\C B}(\sigma)&=[\sigma : \C G_{\C B}(\chi_{reg})]\enspace .
\end{aligned}
\end{equation}
In the case of $G_c$ viewed as a real group, the cones defined by
(\ref{eq:9.3.1}) are parametrized by nilpotent orbits in $\fk g_c.$ In
other words, $\pi'\lr\pi$ if and only if $AS(\pi')\subset
\ovl{AS(\pi)}.$ So let $\C C(\C O_c)$ be the cone corresponding to $\C O_c.$
Note that in this case $W_c\cong W\times W,$ and the representations
are of the form $\sig\otimes\sig.$
\begin{theorem} 
 {Let $G$ be the real points of a connected linear
  reductive group $G_c.$} Then 
$$
|\C B(S_1,S_2,\C O_c)|=\sum_{\sigma\otimes \sigma \in \C C(\C O_c)}
m_{\C B}(\sigma)m_S(\sigma)\enspace.
$$
\end{theorem}
\begin{proof}
Consider $\bZ\C B(S_1,S_2,\C O_c)\subset \bZ\C B(\emptyset,\emptyset,\C
O_c).$ Then $\C B(\emptyset,\emptyset,\C O_c)$ is a representation of
$W$ which consists of the representations in $\C V(\pi)$ with
$AS(\pi)\subset \ovl{\C O_c}.$   The fact that the representation $\C
V(\pi)$ is formed  of $\sig$ with $\sig\otimes\sig\in \C C(\C O_c)$
follows from the argument before theorem 1 of \cite{McG}. This
accounts for $m_\C B(\sig)$ in the sum. The
expressions of the action of $W$ given by lemma 14.7 in \cite{V2}
and Frobenius reciprocity imply that the dimension of $\bZ\C
B(S_1,S_2,\C O_c)$ equals the left hand side of the formula in the
theorem, and it equals the cardinality of $\C B(S_1,S_2,\C O_c)$.
\end{proof}
\subsection{}\label{sec:9.4}
Assume that $\vO$ is even.  
Then  $\lambda:=\vh/2$ is integral, and  it defines a set $S_2$ by
\begin{equation}
  \label{eq:9.2.25}
S_2=S(\la)=\{\alpha \in \Pi_a | (\alpha,\lambda)=0 \}.  
\end{equation}
Let $\C O_c$ be the nilpotent orbit attached to $\vO$ by the duality
in \cite{BV3}. Then the \textit{special unipotent representations
$Unip(\cCO$ attached to $\cCO$}, 
are defined to be the representations $\pi$ with
infinitesimal character $\la$ and $AS(\pi)\subset \C O_c.$ 
{The 
translation principle discussed in
greater detail in chapter 16 particularly lemma 16.2 and theorem 16.4
of \cite{ABV} (see also definition 8.6 in \cite{ABV}), 
establish a  bijection between the irreducible admissible representations in
$\bigcup \C B(\emptyset,S(\la))$ at regular integral infinitesimal character
$\la'$, and irreducible admissible representations at infinitesimal
character $\la.$ Because this bijection is realized as tensoring with
a finite dimensional representation and projecting onto the
appropriate infinitesimal character, $AS(\pi)$ is preserved. Thus
there is a bijection 
\begin{equation}
  \label{eq:9.2.26}
Unip(\cCO)\longleftrightarrow\bigcup_{\C B}\C B(\emptyset,S(\la),\C O_c). 
\end{equation}
So we can use theorem \ref{t:9.3} to count the number of unipotent
representations.} 
In the  classical groups case, $m_{\C B}(\sigma)$ is straightforward to
compute.  For the special unipotent case, $m_{S}(\sigma)$
equals 0 except for the representations occuring in a particular left cell
sometimes also called the \textit{Lusztig cell}, which we denote 
{$\ovl{\C C}^L(\C O_c).$ The multiplicities of the
representations occuring in  $\ovl{\C C}^L(\C O_c)$ are all
1. } These
representations are in 1-1 correspondence with the conjugacy classes
in Lusztig's quotient of the component group $\ovl{A}({\cCO}).$ See
\cite{BV2} for details.  
\begin{theorem}[1]\label{t:9.2.2}
$$
|Unip(\cCO)|=\sum_{\C B}\sum_{\sigma\in \ov{\C C}^L(\C O_c)}
m_{\C B}(\sigma)\enspace . 
$$ 
\end{theorem}

\begin{theorem}[2, \cite{McG}]\label{t:9.3} 
In the classical groups $Sp(n),\  SO(p,q),$ each
  Harish-Chandra cell is of the form  $\ovl{\C C}^L(\C O_c).$ 
\end{theorem}

\begin{definition}\label{d:9.2}
We say that a nilpotent orbit $\C O_c$ is {\rm smoothly cuspidal} if it
satisfies
\begin{description}
\item[Type B] all odd sizes except for the largest one occur an even
  number of times, 
\item[Type C] all even sizes occur an even number of times,
\item[Type D] all odd sizes occur an even number of times, 
\end{description}
\end{definition}
For $\CO(\bb R), $ a real form of $\CO_c,$  write $A(\CO(\bb R))$ for its
(real) component group.
\begin{proposition}[1]\label{p:9.2} 
{Let $\C O_c$ be a smoothly cuspidal orbit with dual
  $\vO.$ Then $A(\vO)=\ovl{A}(\vO),$  and in particular, 
$|\C C^L(\C O_c)|=|A(\vO)|.$} Furthermore, 
\begin{equation*}
|Unip(\cCO))|= \sum_{\CO(\bb R)} |A(\CO(\bb R))|
\end{equation*}
where the sum is over all real forms $\C O(\bR)$ of $\CO_c$. 
\end{proposition}
\begin{proof}
This is theorem 5.3 in  \cite{B2}. The proof consists of a direct
calculation of multiplicities in the coherent continuation
representation using the results developed earlier in this section.
\end{proof}
{
\subsection*{Remark} This proposition does not hold for orbits that
are not smoothly cuspidal. An example is when $\C O_c$ is the principal
nilpotent orbit.\qed
} 
{

\bigskip
We omit the proof of the following proposition. It involves results of
Borho, MacPherson, Joseph and Kashiwara which say that the $AS-$set
of an irreducible highest weight module is the closure of a single
nilpotent orbit. See references in \cite{BV3}, and also another proof
in \cite{V4}.

\begin{proposition}[2]
  \label{p:9.4}
Assume $\vO$ is even. Then the $AS-$set of $L(\chi_\cCO)$
satisfies the property that 
$\ovl{\Ad G_c (AS(L(\chi_\cCO)))}$ is the closure of the special orbit $\C O_c$
dual {(in the sense of \cite{BV3})} to $\cCO.$
\end{proposition}
}
\subsection{}\label{sec:9.5}
We now return to type $G=SO(2n+1).$ Consider the spherical irreducible
representation $L(\chi_\vO)$ with $\chi_\vO=\vh/2$ 
corresponding to a nilpotent orbit $\cCO$
in $sp(n).$  If the orbit $\cCO$ meets a proper Levi component
$\cfm,$ then $L(\cCO)$ is a subquotient of a representation which is
unitarily induced from a unipotent representation on $\fk m$ 
{(compare with (\ref{eq:1.3})).}
 By 
induction, $L(\chi_\vO)$ is unitary. 
Thus we only consider the cases when
$\cCO$ does not meet any proper Levi component. This means
\begin{equation}\label{eq:9.5.1}
\cCO=(2x_0,\dots , 2x_{2m}),\qquad 0\le x_0<\dots <x_i<x_{i+1}<\dots
<x_{2m} ,
\end{equation}
so these orbits are even.
{
\begin{proposition}\label{p:9.5} Let $\vO$ be as in (\ref{eq:9.5.1}). 
The Harish-Chandra cell of $L(\chi_\vO)$ corresponds to $\C O_c$
with partition 
\begin{equation}
  \label{eq:9.5.2}
(\unb{r_1}{1,\dots ,1},\unb{r_2}{2,\dots ,2},\dots
,\unb{r_{2m}}{2m,\dots ,2m},\unb{r_{2m+1}}{2m+1,\dots ,2m+1}),
\end{equation}
where 
\begin{equation*}
  \begin{aligned}
 r_{2i+1}&=2(x_{2m-2i}-x_{2m-2i-1}+1),\\
r_{2i}&=2(x_{2m-2i+1}-x_{2m-2i}-1),\\
r_{2m+1}&=2x_0+1.
  \end{aligned}
\end{equation*}
The columns of $\CO_c$ are $(2x_{2m}+1,2x_{2m-1}-1,\dots
,2x_0+1).$ 
\end{proposition}
\begin{proof}
The proof is a calculation involving Lusztig cells and their duals. We
omit the details. 
\end{proof}
 }

{
\subsection*{Remark} In general if $\vO$ is distinguished, then $\C
O_c$ is smoothly cuspidal. But there are many more smoothly cuspidal
orbits. For example the orbit with partition $(2,2)$ in $sp(4)$ is
smoothly cuspidal, but its dual is $\vO$ with partition $(1,1,3)$ which is
not distinguished.
}

\begin{definition} Given an orbit $\CO_c$ with partition
  (\ref{eq:9.5.2}) or more generally a smoothly cuspidal orbit,
we call the {\rm split real form} $\CO_{spl}$ the one which satisfies for each
row size,  
\begin{description}
\item[Type C,D] the number of rows starting with $+$ and $-$ is equal,
\item[Type B] in addition to the condition in types C,D for rows of
  size less than $2m+1$, for size $2m+1,$ the number of
  starting with $+$ is one more than those starting with $-.$ 
\end{description} 
\end{definition}
\begin{theorem}
The $WF-$set of the spherical representation $L(\chi_\vO)$ with $\cCO$
satisfying (\ref{eq:9.5.1}) is the closure of the
split real form $\CO_{spl}$ of the (complex) orbit $\CO_c$  given by
(\ref{eq:9.5.2}). 
\end{theorem}
\begin{proof}
The main idea is outlined in \cite{B2}. We use the fact that if $\pi$
is a factor of $\pi',$ then $WF(\pi)\subset WF(\pi').$ We do an
induction on $m$, {the number of parts in the partition of
  $\cCO$.} The claim amounts to showing that if $E$ occurs in
$WF(L(\chi_\vO)),$ then the signatures of $E,\ E^2,\dots $ are greater than
the pairs 
\begin{equation} \label{eq:9.5.3}
\begin{aligned}
 &(x_{2m}+1,x_{2m}),(x_{2m}+x_{2m-1},x_{2m}+x_{2m-1}),\dots ,\\
&\dots (x_{2m}+\dots + x_1,x_{2m}+\dots + x_1),\\
&(x_{2m}+\dots + x_1+x_0+1,x_{2m}+\dots + x_1+x_0).
\end{aligned}
\end{equation}
{
The first pair is the signature of $E,$ the second the signature of
$E^2$ and so on. The first pair of numbers is obtained by counting the
$+$'s and $-$'s in the first column of the nilpotent orbit, the second
pair by counting the same in the first two columns and so on.
}

The statement is clear when $m=0;$ $L(\chi_\vO)$ 
 is the trivial representation. Let $\cCO_1$ be
the nilpotent orbit corresponding to 
\begin{equation}
  \label{eq:9.5.4}
  (2x_0,\dots, 2x_{2m-2}).
\end{equation}
By induction, $WF(L(\chi_{\cCO_1}))$ is the split real form of the nilpotent
orbit corresponding to the partition
\begin{equation}
  \label{eq:9.5.5}
(\unb{r'_1}{1,\dots ,1},\unb{r'_2}{2,\dots ,2},\dots
,\unb{r'_{2m-2}}{2m-2,\dots ,2m-2},\unb{r'_{2m-1}}{2m-1,\dots ,2m-1}),
\end{equation}
where the columns are $(2x_{2m-2}+1,2x_{2m-3}-1,\dots ,2x_0+1)$.
Let $\fk p$ be the real parabolic subalgebra with Levi component $\fk
g(n-x_{2m}-x_{2m-1})\times gl(x_{2m}+x_{2m-1}).$ There is a character
$\chi$ of $gl(x_{2m}+x_{2m-1})$ such that $\pi:=L(\chi_\vO)$ is a factor
of $\pi':=Ind_{\fk p}^{\fk g} [L(\chi_{\cCO_1})\otimes\chi].$ But by section
\ref{sec:8}, $WF(\pi')$ is in the closure of nilpotent orbits
corresponding to partitions 
{
\begin{equation}
  \label{eq:9.5.6}
(\unb{r_1/2+r_2}{2,\dots ,2},\dots,\unb{r_{2m}}{2m,\dots  ,2m},
\unb{r_{2m+1}}{2m+1,\dots ,2m+1}),{\text{ if}}\ r_1
\text{ even,}
\end{equation}
{or}
\begin{equation}
  \label{eq:9.5.7}
(1,1,\unb{(r_1-1)/2+r_2}{2,\dots,2},\dots
,\unb{r_{2m}}{2m,\dots ,2m},\unb{r_{2m+1}}{2m+1,\dots ,2m+1}),
{\text{ if}}\ r_1 \text{ odd,}
\end{equation}
and the dots refer to sizes as in (\ref{eq:9.5.2}).}
It follows that the signatures for $E^k$ in
$WF(L(\chi_\vO))$ are greater than the pairs
\begin{equation}
  \label{eq:9.5.8}
(a_+,a_-),\ (x_{2m}+x_{2m-1},x_{2m}+x_{2m-1}),\dots\ ,
\end{equation}
for some $a_+ + a_-=x_{2m}+1,$ {and the dots are pairs as in (\ref{eq:9.5.3})}.
Also, each row size greater than two and less than $2m+1$ has an equal
number that start with  $+$ and $-.$ For size $2m+1$ there is one
more row starting with $+$ than $-.$

{
We now do the same argument with $\cCO_2$ corresponding to 
\[
(2x_0,\dots,{2x_{2m-3},}
\widehat{2x_{2m-2}},\widehat{2x_{2m-1}},2x_{2m}).
\]
The $WF-$set correponding to $L(\chi_{\vO_2})$ is the split real form
of the nilpotent orbit with partition
\begin{equation}
  \label{eq:9.5.8.1}
 (\unb{r^{''}_1}{1,\dots ,1},\unb{r^{''}_2}{2,\dots ,2},\dots
,\unb{r^{''}_{2m-2}}{2m-2,\dots ,2m-2},\unb{r^{''}_{2m-1}}{2m-1,\dots ,2m-1}), 
\end{equation}
with columns $(2x_{2m}+1,2x_{2m-3}-1,\dots ,2x_0+1).$ 
Let $\fk p$ be the real parabolic subalgebra with Levi component $\fk
g(n-x_{2m-1}-x_{2m-2})\times gl(x_{2m-1}+x_{2m-2}).$ There is a character
$\chi$ of $gl(x_{2m-1}+x_{2m-2})$ such that $\pi:=L(\chi_\vO)$ is a factor
of $\pi^{''}:=Ind_{\fk p}^{\fk g} [L(\chi_{\cCO_2})\otimes\chi].$ But by section
\ref{sec:8}, $WF(\pi^{''})$ is in the closure of nilpotent orbits
corresponding to partitions  
\begin{equation}
  \label{eq:9.5.8.2}
(\unb{r_1+r_2/2}{1,\dots ,1},\unb{r_2/2+r_3}{3,\dots ,3},\dots,
\unb{r_{2m}}{2m,\dots  ,2m},
\unb{r_{2m+1}}{2m+1,\dots ,2m+1}),{\text{ if}}\ r_2
\text{ even,}
\end{equation}
{or}
\begin{equation}
  \label{eq:9.5.8.3}
(\unb{r_1+(r_2-1)/2}{1,\dots ,1},{2,2},\unb{(r_2-1)/2+r_3}{3,\dots ,3}\dots
,\unb{r_{2m}}{2m,\dots ,2m},\unb{r_{2m+1}}{2m+1,\dots ,2m+1}),
{\text{ if}}\ r_2 \text{ odd.}
\end{equation}
where the dots refer to pairs as in (\ref{eq:9.5.2}).
Thus $WF(L(\chi_\vO))$ is also contained in the closure of the
nilpotent orbits with signatures 
\begin{equation}
  \label{eq:9.5.9}
  \begin{aligned}
 &(x_{2m}+1,x_{2m}),\ (x_{2m}+1+a_+,x_{2m}+a_-),\\
 &(x_{2m}+1+x_{2m-1}+x_{2m-2},x_{2m}+1+x_{2m-1}+x_{2m-2}),\ \dots \ ,
\end{aligned}
\end{equation}
 for some $a_++a_-=x_{2m-1},$ {and the dots refer to pairs as in
   (\ref{eq:9.5.3}).} 

The claim follows.
} 
\end{proof}

\subsection{}\label{sec:9.6}
Consider the special case when
\begin{equation}
  \label{eq:9.6.1}
x_0=x_1-1\le x_2=x_3-1\le \dots \le x_{2m-2}=x_{2m-1}-1\le x_{2m}.
\end{equation}
{This makes $\CO_c$ of the form (\ref{eq:9.5.2}) with
  $r_{2i}=0,$ \ie the partition of $\C O_c$ has only odd terms.}
The (Lusztig) cell {$\ovl{\C C}^L(\C O_c)$} has size $2^m$ 
{(which is true for any $\C O_c$ with $\cCO$ distinguished, not just for the
  particular form (\ref{eq:9.6.1}))}.
We produce $2^m$ distinct irreducible representations with {$AC-$set} equal to
$\CO_{spl}.$  So $\fk g$ is $so(2p+1,2p)$. 
Let $\fk h$ be the compact Cartan subalgebra. We write the coordinates 
\begin{equation}
  \label{eq:9.6.2}
  (a_1,\dots , a_p\ |\ b_1,\dots ,b_p)
\end{equation}
where the first $p$ coordinates before the $|$ are in the Cartan
subalgebra of $so(2p+1)$ the last $p$ coordinates are in $so(2p).$ The
roots $\ep_i\pm \ep_j, \ep_i$ with $1\le i,j\le p$ are all compact and so
are $\ep_{p+k}\pm\ep_{p+l}$ with $1\le k,l\le p.$ The roots $\ep_i\pm
\ep_{p+k},\ \ep_{p+k}$ are noncompact. Let  $\fk q_c=\fk l_c +\fk u_c$ be a
$\theta$-stable parabolic subalgebra with Levi component

\begin{align}
  \label{eq:9.6.3}
  \fk l:=\fk l_c\cap\fk g=&u(x_{2i_1+1},x_{2i_1})\times
  u(x_{2i_2},x_{2i_2+1})\times\dots\times\notag\\ 
   & u(x_{2i_{m-1}},x_{2i_{m-1}+1})\times so(x_{2m}+1,x_{2m}),\notag\\
\text{ or }& \\
\fk l:=\fk l_c\cap\fk g=&u(x_{2i_1+1},x_{2i_1})\times
  u(x_{2i_2},x_{2i_2+1})\times\dots\times\notag\\ 
    &u(x_{2i_{m-1}+1},x_{2i_{m-1}})\times so(x_{2m}+1,x_{2m}),\notag
\end{align}
{(depending on the parity of $m$)}
where the $i_j$ are the numbers $0,\dots ,m-1$ in some order.  
The parabolic subalgebra $\fk q_c$ corresponds to the weight
\begin{equation}
  \label{eq:9.6.4}
  \begin{aligned}
&\xi=(m^{x_{2i_1+1}},\dots ,1^{x_{2i_{m-1}}+1},0^{x_{2m}}\ \mid\
  m^{x_{2i_1}},\dots , 1^{x_{2i_{m-1}}},0^{x_{2m}}),\\
&\text{ or } \\
&\xi=(m^{x_{2i_1+1}},\dots ,1^{x_{2i_{m-1}}},0^{x_{2m}}\ \mid\
  m^{x_{2i_1}},\dots , 1^{x_{2i_{m-1}}+1},0^{x_{2m}}),\\
  \end{aligned}
\end{equation}
depending whether $m$ is odd or even.

The derived functor modules $\C R_{\fk q_c}^i(\xi)$ from characters on
$\fk l_c$ have $AC$-set contained in $\CO_{spl}.$ To get infinitesimal
character $\chi_\cCO,$ these characters can only be
\begin{equation}
  \label{eq:9.6.5}
\xi_{i_j}^\pm:=\pm (1/2,\dots ,1/2),
\end{equation}
on the unitary factors $u(x_{2i_j+1},x_{2i_j})$ or
$u(x_{2i_j},x_{2i_j+1}),$ and trivial on $\fk g(x_{2m}).$ We need to
show that there are choices of parabolic subalgebras $\fk q_c$ as in
(\ref{eq:9.6.3}) and characters as in (\ref{eq:9.6.5}) so that we get
$2^m$  nonzero and distinct representations. For this we have to
specify the Langlands parameters.  

\medskip
For each subset  $A:=\{k_1,\dots ,k_r\}\subset\{0,\dots, m-1\},$ $k_j$ in
decreasing order, label the complement $A^c:=\{\ell_1,\dots ,
\ell_t\},$  {$A^c$ also labeled in decreasing order,} and
consider the $\theta-$stable parabolic subalgebra $\fk 
q_{c,A}$ as in (\ref{eq:9.6.3}) and (\ref{eq:9.6.4}) corresponding to 
\begin{equation}
  \label{eq:9.6.6}
 \{ i_1,\dots ,i_{m-1}\}=\{k_1,\dots ,k_r,\ell_1,\dots,\ell_t\}.
\end{equation}
We will consider the representations $\C R_{\fk q_{c,A}}(\xi_A),$ where
$\xi_A$ is the concatentation of the $\xi^\pm_{i_j}$ with $+$ for the
first $r$, and $-$ for the last $t.$ 
\begin{proposition}
  \label{p:9.6}

\begin{equation*}
  \CR_{\fk q_{c,A}}^i(\xi_A)=
  \begin{cases}
    &0\text{ if } i\ne \dim(\fk u_{c,A}\cap \fk k_c),\\
    &\text{nonzero irreducible if } i=\dim(\fk u_{c,A}\cap \fk k_c).
  \end{cases}
\end{equation*}
\end{proposition}
\begin{proof}
For the vanishing part we check that the conditions in proposition
5.93 in  \cite{KnV} chapter V section 7 are satisfied. 
It is sufficient to show that 
\begin{equation}
  \label{eq:9.6.7}
ind_{\ovl{\fk q}_{c,A},L\cap K}^{\fk g, K}(Z_{\ovl{\fk q}_{c,A}}^\#):=
U(\fk g)\otimes_{\ovl{\fk q}_{c,A}} Z_{\ovl{\fk q}_{c,A}}^\#   
\end{equation}
is irreducible. Here $Z_{\ovl{\fk q}_{c,A}}^\#$ is the $1-$dimensional
module corresponding to $\xi_A -\rho(\fk u_{c,A}),$ with 
\begin{equation*}
  \rho(\fk u_{c,A}):=\frac12\sum_{\al\in\Delta(\fk u_{A,c})} \al.
\end{equation*}
The derived functors are normalized
so that if $W$ has infinitesimal character $\chi,$ then so do the $\C
R^i_{\fk q_c}(W).$ 
{
We prove a slightly more general result. Write the parameter of an
induced module as in (\ref{eq:9.6.7}) as
\begin{equation}
  \label{eq:vanish1}
  (a_1,\dots ,b_1;a_2,\dots b_2,\dots ;a_m,\dots ,b_m;-r+1/2,\dots ,-1/2)
\end{equation}
where each $(a_i,\dots,b_i)$ is a string \ie the coordinates go up by
one. Each string encodes a 1-dimensional character on a
$gl(b_i-a_i+1)$. The last one with $-r+1/2$  denotes the trivial
character on $\fk g(r).$ 
\begin{lemma}
  \label{l:vanish}
Assume that the parameter in equation (\ref{eq:vanish1}) satisfies the
following.
\begin{enumerate}
\item For each $i,j$ either
  \begin{align*}
    &a_i\le a_j\le b_j\le b_i\\
&\text{ or }\\
    &a_i\le a_j\le b_j\le b_i.
  \end{align*}
\item $a_i +b_i\le 1.$
\item $|a_i|,\ |b_i|\le r-1/2.$
\end{enumerate}
Then the module in equation (\ref{eq:9.6.7}) is irreducible.
\end{lemma}
\begin{proof}
Any factor of (\ref{eq:9.6.7}) would have to be a lowest weight module
which is a generalized Verma module for
$\ovl{\fk q}.$ Thus its parameter would be
\begin{equation}
  \label{eq:vanish2}
  (a'_1,\dots ,b'_1;\dots :a'_k,\dots ,b'_k;-R'+1/2,\dots ,-r'+1/2)
\end{equation}
where the coordinates in each $(a'_i,\dots ,b'_i)$ are increasing and
differ by integers, and the same is true for $r',R'\in \bN.$ These
\textit{strings} encode finite dimensional representations on
$gl(b_i-a_i+1)$ and $\fk g(r)$ respectively. To be a factor of (\ref{eq:9.5.7}),
\begin{enumerate}
\item the coordinates of (\ref{eq:vanish2}) have to be in the same
  Weyl orbit as the coordinates (\ref{eq:vanish1}) (so that the
  infinitesimal character are the same).
\item the difference between (\ref{eq:vanish2}) and (\ref{eq:vanish1})
  must be a sum of roots in $\Delta(\fk u).$
\end{enumerate}
Given the assumptions, (\ref{eq:vanish1}) and (\ref{eq:vanish2}) must
coincide except for some subset of the $i$  for which $a_i+b_i=1,$ and 
\[
(a'_i,\dots ,b'_i)=(-b_i,\dots ,-a_i).
\]
We do an induction on the number of such $i.$ The module
(\ref{eq:vanish1}) is irreducible if there are no such $i.$ Suppose
there are such $i.$ Then the possible factors are of the same kind,
but with fewer such $i.$ The theory of asymptotic expansions of
characters outlined in section \ref{sec:9.1} coming from \cite{BV1}
applies. Any two induced
modules have the same leading term. So if (\ref{eq:vanish1}) contains
a factor of the type (\ref{eq:vanish2}), the remaining factors would
have to have strictly smaller Gelfand-Kirillov dimension. But the
infinitesimal character is of the special unipotent type considered in
\cite{BV2}, so there are no such factors.
\end{proof}


}
\medskip
{
The notions and results about associated cycles in chapter
\ref{sec:8} section \ref{sec:thetaind} apply.
By the corollary   \ref{sec:thetaind}, $\CR_{\fk q_{c,A}}^{\dim(\fk
  u_{c,A}\cap\fk     k_c)}(\xi_A)$ is irreducible or zero because the
multiplicity  of the corresponding nilpotent orbit in $\fk s_c$ in
$AV$ is 1.} 
To show that $\CR_{\fk q_{c,A}}^{\dim(\fk u_{c,A}\cap \fk
  k_c)}(\xi_A)\ne 0,$ we use the bottom layer $K-$ types defined in
chapter V section 6 of \cite{KnV}. To simplify the notation slightly,
we write 
\begin{equation}
  \label{eq:9.6.8}
  \begin{aligned}
&a_1=x_{2k_1+1},\ b_1=x_{2k_1},\dots , a_r=x_{2k_r},\ b_r=x_{2k_r+1}\
r \text{ even},\\
&a_1=x_{2k_1+1},\ b_1=x_{2k_1},\dots ,  a_r=x_{2k_r+1},\ b_r=x_{2k_r}    
\ r \text{ odd},
  \end{aligned}
\end{equation}
{and $\fk u_c$ for $\fk u_{c,A}.$} Let also $a:=\sum a_j,\ b:=\sum b_j.$ 
Note that $|a_j-b_j|=1,$  and also $|a-b|=1.$ Then
{
\begin{equation}
  \label{eq:9.6.9}
\mu:=\xi +2\rho(\fk u \cap \fk s) -\rho(\fk u)=(1^{a},0^{p-a}\mid
1^{b},0^{p-b}) 
\end{equation}
}
is dominant, therefore bottom layer. The aforementioned results then
imply the nonvanishing. 
\end{proof}
We now show that there are $2^m$ distinct representations. 
{To do this we will compute Langlands parameters. We rely
on the notation and conventions in chapter 11 of \cite{ABV}. Many
details are in the references to earlier work
given there. In \cite{ABV} a Langlands parameter is a triple
$\La=(\La^{can},\bR^+_{i\bR},\bR^+_{\bR}),$ where $\La^{can}$  is an
irreducible representation of a $\theta-$stable Cartan subgroup
$H_\bR$ (denoted $T_\bR$ in \cite{ABV}). Attached to such data is a
\textit{standard module} X$(\La)$, with a canonical irreducible
subquotient $M(\La).$  We will determine $\wti\la:=d\La^{can}$ and the
lowest $K-$types of the derived functor modules constructed
earlier. The Lie algebra of $H$ decomposes according to the $\pm 1$
eigenvalues of $\theta$ into $\fk h=\fk t +\fk a.$ Then
$\wti\la=(\la,\nu).$ Recall from \cite{V3} chapter 5  that the
centralizer of $\fk a$ is a quasisplit $\theta-$stable Levi component
$\fk l,$ and that $M(\La)$ is determined by the fact that it contains
a \textit{lowest $K-$type} $\mu$ related to $\la$ by the formula
$\la=\mu+2\rho_c +\rho+v$ where $v=\sum c_i\al_i$  is a combination of
strongly orthogonal imaginary roots with coefficients $0<c_i\le 1/2.$
}

To compute $\la$ and $\mu$, we will
need to use the intermediate parabolic subalgebras  
\begin{equation}
  \label{eq:9.6.10}
  \fk q_{c,A}\subset \fk q_{c,A}'\subset \fk q^{''}_{c,A}\subset \fk g_c
\end{equation}
with Levi components 
\begin{equation}
  \label{eq:9.6.11}
  \begin{aligned}
&  \fk l'_{A}= u(a_1,b_1)\times \dots \times u(a_r,b_r)\times\fk
  g(n-a-b),\\
& \fk l^{''}_{A}= u(a,b)\times\fk g(n-a-b),
  \end{aligned}
\end{equation}
Apply induction in stages from $\fk q_{c,A}$ to $\fk q'_{c,A}$ first. 
On the factor $\fk g (n-a-b)$
the $K-$type $\mu$ in (\ref{eq:9.6.9}) is trivial, so the Langlands
parameter is that of the spherical principal series. Similarly on the
$u(a_j,b_j)$ assume the infinitesimal character is $\chi_j:=(\max
(a_j,b_j),\dots ,\min (a_j,b_j))$ with the coordinates going down by
1, and the Langlands parameter is that
of a principal series with the appropriate $1$-dimensional Langlands
subquotient. Let $\fk h_A\subset \fk l'_{c,A}$ 
be the the most split Cartan subalgebra. In particular the real roots
are  
\begin{equation}
  \label{eq:9.6.12}
\al_d:=\ep_d +\ep_{d + p},\quad \sum_{j\le s}a_j< d< \sum_{j\le
  s} a_j +\min(a_j,b_j) , \quad 0\le s\le r-1.   
\end{equation}
For each factor $u(a_j,b_j)$ the Langlands parameter is of the form
$(\la_j,\nu_j)$ where $\la_j\in \fk h_A\cap \fk k_c,$ and $\nu_j\in \fk
h_A\cap \fk s_c.$ Then
\begin{equation}
  \label{eq:9.6.13}
  \la_j=(1/2^{a_j}\ \mid\ 1/2^{b_j}),
\end{equation}
while
\begin{equation}
  \label{eq:9.6.14}
  \langle\nu_j,\al_d\rangle=\max (a_j,b_j)-(d-\sum_{j\le s}a_j)
\end{equation}

\begin{proposition}\label{p:9.6} The  representations 
$\CR_{\fk q_{c,A}}^{\dim(\fk u_{c,A}\cap \fk k_c)}(\xi_A)$ have
Langlands parameters $(\la^G,\nu)$ where $\la^G$ is obtained by
concatenating the $\la_j$ in (\ref{eq:9.6.13}) and $\nu$ satisfies 
(\ref{eq:9.6.14}).
\end{proposition}
\begin{proof}
{We write $\fk l_c'=\fk l_{c,A}'$ and
  $\fk l'=\fk l_{c,A}'\cap\fk g$.} 
There is a nonzero map $X_{\fk l_c'}(\la^G,-\nu)\longrightarrow
L_{\fk l'}(\la^G,-\nu)$ given by the Langlands classification. Thus
there is a map 
\begin{equation}
  \label{eq:9.6.15}
\begin{aligned}
&\C R_{\fk q'_c,L'\cap K}^{\dim \fk k_c\cap\fk u_c'}[X_{\fk  l'}(\la^G,-\nu)]
\longrightarrow\\
&\longrightarrow\C R_{\fk q'_c,L'\cap K}^{\dim \fk k_c\cap\fk
  u_c'}(L_{\fk l'_c}(\la^G,-\nu))=
\C R_\fk {q_c,L\cap K}^{\dim \fk k_c\cap\fk u_c}(\xi_A),
\end{aligned}
\end{equation}
which is nonzero on the bottom layer $K-$type (\ref{eq:9.6.9}). On the
other hand, because these are standard modules,
\begin{equation}
  \label{eq:9.6.16}
  \C R^i_\fk q(X_{\fk l'}(\la^G,\nu))=
  \begin{cases}
    X(\la^G,\nu) &\text{ if } i=\dim\fk k_c\cap\fk u_c,\\
    0            &\text { otherwise.}
  \end{cases}
\end{equation}
The proof follows.
\end{proof}
\subsection{}\label{sec:9.7} 
{
\begin{theorem} The spherical unipotent representations $L(\chi_\vO)$
with $\vO$ even  are unitary.  
\end{theorem}
}
\begin{proof} {We present the details for $G=SO(n+1,n).$} 
Write $\fk g(n)$ for the Lie algebra containing $\CO_{spl}$ 
the split real form of $\C O_c$, same as  the  support 
$AC(L(\chi_{\vO})).$ There is a (real) parabolic subalgebra $\fk p^+$ 
with Levi component $\fk m^+:=gl(n_1)\times \dots
\times gl(n_k)\times\fk g(n)$ in  $\fk g^+$ of rank $n_1+\dots+ n_k+n,$  
such that the split form $\CO^+_{spl}$ of 
{
$$
\CO^+_c:=(1,1,3,3,\dots,2a-1,2a-1, 2a+1)
$$ 
for some $a$ depending on $(2x_0,\dots ,2x_{2m}).$}
is induced from $\CO$ on $\fg(n),$ trivial on the
$gl$'s.  We will
consider the representation
\begin{equation}
  \label{eq:9.7.1}
  I(\pi):=Ind_{\fk m^+}^{\fg^+}[triv\otimes\dots\otimes triv\otimes\pi].
\end{equation}
We show that the form on $I(\pi)$ induced from $\pi$ is positive
definite; this implies that the form on $\pi$ is definite. We do this
by showing that the possible factors of $I(\pi)$ have to be unitary,
and the forms on their lowest $K-$types are positive definite. 

\medskip
{
By the results in chapters \ref{sec:8} and \ref{sec:9}, more specifically
section \ref{sec:9.4}, and particularly  proposition 
(1) in section \ref{p:9.2}, we conclude that there are $3^a\cdot 2^a$ unipotent
representations in the block of the spherical irreducible
representation; all the factors of $I(\pi)$ are in this
block. Here $3^a$ is the number of real forms of $\C O_c^+,$ and $2^a$
is the number of representations in the corresponding Lusztig cell.
}
We describe how to get $3^a\cdot 2^a$ distinct
representations. For each real form 
{$\CO^+_j$ of $\C O^+_c$} 
we produce one representation $\pi$ such that 
$AC(\pi)=\ovl{\CO_j^+}.$ Then {theorem (2) of \ref{t:9.3}}
implies that there is a Harish-Chandra cell with $2^a$ representations with this
property, {since $\overline A(\cCO)=A(\cCO)=(\bZ/2)^a$}. 
Since these cells must be disjoint, this gives the required number. 

From section \ref{sec:9.1}, each such form $\CO_j^+$ is $\theta$-stable
induced from the trivial nilpotent orbit on a parabolic subalgebra
with Levi component a real form of $gl(1)\times gl(3)\times
\dots\times gl(2a-1)\times \fk g_c(m).$ Using the results in \cite{KnV},
for each such parabolic subalgebra, we can find a derived functor
induced module from an appropriate 1-dimensional character, that is
nonzero and has associated variety equal to the closure of the given real form.
Actually it is enough to construct this derived functor module at regular
infinitesimal character where the fact that it is nonzero irreducibile
is considerably easier. 
 
\medskip
So in this block, there is a cell for each real form of $\CO^+_c,$ and each
cell has $2^a$ irreducible representations with infinitesimal
character $\chi_\cCO.$ In particular for $\CO_{spl}$, the Levi
component is $u(1,0)\times u(1,2)\times u(3,2)\times\dots\times
so(a,a+1).$ For this case, section \ref{sec:9.6}  produced exactly
$2^a$ parameters; their lowest  $K-$types are of the form
$\mu_e(n-k,k).$  These are the only possible  constituents of the
induced from $L(\chi_\vO).$ Since the constituents of the restriction
of a $\mu_e(n-k,k)$ to a Levi component are again $\mu_e(m-l,l)$'s, the
only way $L(\chi_\vO)$ can fail to be unitary is if the form is
negative on one of the $K-$types $\mu_e(n-k,k).$  But sections
\ref{sec:5} and \ref{6.2} show that the form is positive on the
$K-$types $\mu_e$ of $L(\chi_\vO).$  
\end{proof}

\section{Irreducibility}\label{sec:10} 
\subsection{}\label{sec:10.1} To complete the classification
of the unitary dual we also need to prove the following irreducibility
theorem. It is needed to show that the regions in theorem
\ref{thm:3.1} are indeed unitary in the real case. 
\begin{theorem}\label{t:10.1}
Assume $\vO$ is even, and such that $x_{i-1}=x_i=x_{i+1}$ for some
$i.$ Let $\fk
m=gl(x_i)\times \fk g(n-x_i),$ and $\vO_1\subset \fk g(n-x_i)$ be the
nilpotent orbit obtained from $\vO$ by removing two rows of size
$x_i.$ Then
\begin{equation*}
  L(\chi_\vO)=Ind_{GL(x_i)\times G(n-x_i)}^{G(n)}[triv\otimes
  L(\chi_{\vO_1})].  
\end{equation*}
\end{theorem}
In the $p-$adic case this follows from the work of Kazhdan-Lusztig
(\cite{BM1}). In the real case, it follows from the following proposition.
\begin{proposition}\label{p:10.1}
The associated variety of a spherical representation $L(\chi_\vO)$ is
given by the sum with multiplicity one of the following 
{real forms of the complex orbit $\CO_c$}.
\begin{description}
\item[Type B, D]  On the odd sized rows, the difference between the
  number of $+$'s and number of $-$'s is 1, 0 or -1. 

\item[Type C]  On the even sized rows, the difference between the
  number of $+$'s and number of $-$'s is 1, 0 or -1. 

\end{description}
\end{proposition}
The proof of the proposition is lengthy, and follows from more general
results which are unpublished (\cite{B5}). We will give a different
proof of theorem \ref{t:10.1} in the next sections.
\subsection*{Remark} When $\vO_1$ is
even, but $\vO$ is not, and just $x_i=x_{i+1},$ the proof follows from
\cite{BM1} in the $p-$adic case, and the Kazhdan-Lusztig conjectures for
nonintegral infinitesimal character in the real case. We have already
used these results in the course of the paper. \qed

\medskip
The outline of the proof is as follows.
In section \ref{sec:10.2}], we prove some auxiliary reducibility results in the case
when $\vO$ is induced from the trivial nilpotent orbit of a maximal
Levi component. In section \ref{sec:10.3}, we combine these results with
intertwining operator techniques to  complete the proof of theorem
\ref{t:10.1}. 

\subsection{}\label{sec:10.2}  We need to study the $\rho-$induced modules
from the trivial module on $\fk m\subset \fg(n)$ where $\fk m\cong
gl(n),$ or in some cases $\fk m\cong gl(a)\times \fg(b)$ with $a+b=n.$

\subsubsection*{Type B} Assume $\cCO$ corresponds to the
  partition $2x_0=2x_1=2a$ in $sp(n,\bC).$ The infinitesimal character
  is $(-a+1/2,\dots ,a-1/2)$ and the nilpotent orbit $\CO_c$ corresponds to
  $(1,1,\unb{2a-2}{2,\dots ,2},3).$ We are interested in the
  composition series of  
\begin{equation}
  \label{eq:10.2.1o}
  Ind_{GL(2a)}^{G(2a)}[triv].
\end{equation}
There are three real forms   of $\CO_c$ in $so(2a+1,2a),$
  \begin{equation}\label{eq:10.2.1}
    \begin{matrix}
      +&-&+\\
      +&-&\\
      -&+&\\
      \vdots&\vdots&\\
      +&-&\\
      -&+&\\
      +&&\\
      +&&
    \end{matrix}
\qquad\qquad
    \begin{matrix}
      +&-&+\\
      +&-&\\
      -&+&\\
      \vdots&\vdots\\     
      +&-&\\
      -&+&\\
      +&&\\
      -&&
    \end{matrix}
\qquad\qquad
    \begin{matrix}
      -&+&-\\
      +&-&\\
      -&+&\\
      \vdots&\vdots\\
      +&-&\\
      -&+&\\
      +&&\\
      +&&
    \end{matrix}
  \end{equation}
The associated cycle of (\ref{eq:10.2.1o}) is the middle nilpotent orbit in
(\ref{eq:10.2.1}) with multiplicity 2. 
Section \ref{sec:6} shows that there are at least two factors
characterized by the fact that they contain the relevant $K-$types which are the
restrictions to $S[O(2a+1)\times O(2a)]$ of 
\begin{equation}
  \label{eq:10.2.3}
  \begin{aligned}
&\mu({0}^{a};+)\otimes \mu(0^a;+)\\
&\mu(1,0^{a-1};-)\otimes \mu(0^a;-).
\end{aligned}
\end{equation}
Thus because of multiplicity 2, there are exactly two factors.
One of the factors is spherical. The nonspherical factor has Langlands parameter
\begin{equation}
  \label{eq:10.2.1p}
  \begin{aligned}
\la^G&=(1/2,0,\dots ,0\mid 0,\dots ,0),\\
\nu&=(0,a-1/2,a-1/2,\dots ,3/2,3/2,1/2).
  \end{aligned}
\end{equation}
The Cartan subalgebra for the nonspherical parameter is such that the
root $\ep_1$ is noncompact imaginary, $\ep_i,\ep_i\pm\ep_j$ with $j>i\ge 2,$ are
real. The standard module $X(\la^G,\nu)$ which has
$\ovl{X}(\la^G,\nu)$ as quotient is the one for which $\nu$ is
dominant. Thus we conjugate the Cartan subalgebra such that $\ep_{2a}$ is
noncompact imaginary, $\ep_i,\ep_i\pm\ep_j$ with $i<j<2a$ are real,
and the usual positive system $\Delta^+=\{ \ep_i,\ep_i\pm\ep_j\}_{i<j}$.   

\subsubsection*{Type C} Consider $\cCO$ which corresponds to the
partition   $2x_0=2x_1=2a+1<2x_2=2b+1$ in $so(n,\bC).$ The infinitesimal
  character is 
  \begin{equation}
    \label{eq:10.2.2}
    (-a,\dots ,a)(-b,\dots ,-1)
  \end{equation}
The nilpotent orbit $\CO_c$ is induced from the trivial one on
$gl(2a+1)\times \fk g(b)$ and corresponds to 
\begin{equation}
  \label{eq:10.2.3a}
  (\unb{2b-2a-2}{1,\dots ,1},2,2,\unb{2a}{3,\dots ,3}).
\end{equation}
We are interested in the composition series of 
\begin{equation}
  \label{eq:10.2.4a}
  Ind_{GL(2a+1)\times G(b)}^{G(2a+b+1)}[triv].
\end{equation}
There are three real forms of (\ref{eq:10.2.3a}), 
 \begin{equation}\label{eq:10.2.4}
    \begin{matrix}
      +&-&+\\
      -&+&-\\
      \vdots&&\vdots\\
      +&-&+\\
      -&+&-\\
      +&-&\\
      +&-&\\
       +&&\\
      -&&\\
     \vdots&&\\
      +&&\\
      -&&
    \end{matrix}
\qquad\qquad
   \begin{matrix}
      +&-&+\\
      -&+&-\\
      \vdots&&\vdots\\
      +&-&+\\
      -&+&-\\
      +&-&\\
      -&+&\\
       +&&\\
      -&&\\
     \vdots&&\\
      +&&\\
      -&&
    \end{matrix}
\qquad\qquad
  \begin{matrix}
      +&-&+\\
      -&+&-\\
      \vdots&&\vdots\\
      +&-&+\\
      -&+&-\\
      -&+&\\
      -&+&\\
      +&&\\
      -&&\\
     \vdots&&\\
      +&&\\
      -&&
    \end{matrix}
\end{equation}
The AC cycle of (\ref{eq:10.2.4a}) consists of the middle nilpotent orbit in
(\ref{eq:10.2.4}) with multiplicity 2. 
By a similar argument as for type B, we conclude that
the composition series consists of two representations containing the
relevant $K-$types  
\begin{equation}
  \label{eq:10.2.4c}
  \begin{aligned}
    &\mu(0^n),\\
    &\mu(1^{a+1},0^{b-1},(-1)^{a+1}).
  \end{aligned}
\end{equation}
These are also the lowest $K-$types of the representations.
The nonspherical representation  has parameter
\begin{equation}
  \label{eq:10.2.4b}
  \begin{aligned}
  &\la ^G=(1/2^a,0^b,(-1/2)^a),\\
  &\nu=(1/2^a,0^b,1/2^a).
  \end{aligned}
\end{equation}

\subsubsection*{Type D} Let $\cCO$ correspond to the
partition 
  $2x_0=2x_1=2a+1$ in $so(n,\bC).$ The infinitesimal character is
  $(-a,\dots ,a).$ The real forms of the nilpotent orbit $\CO$ are
  \begin{equation}
    \label{eq:10.2.5}
    \begin{matrix}
      +&-&\\
      -&+&\\
      \vdots&\vdots\\
      +&-&\\
      -&+&
    \end{matrix}   
  \end{equation}
There are two nilpotent orbits with this partition labelled $I,\ II.$
Each of them is induced from  $\fk m\cong gl(2a),$ there are two such
Levi components. We are interested in the induced modules 
\begin{equation}
  \label{eq:10.2.6}
  Ind_{GL(2a)}^{G(2a)}[triv].
\end{equation}
The multiplicity of the nilpotent orbit (\ref{eq:10.2.5}) in the AC
cycle of (\ref{eq:10.2.6}) is 1, so the representations are irreducible.

\medskip
We summarize these calculations in a proposition.
\begin{proposition}\label{p:10.2} 
The composition factors of the induced module from the trivial
representation on $\fk m$, {where $\fk m=gl(2a)$ for
  $SO(2a+1,2a)$ and $SO(2a,2a)$, or $gl(2a+1)\times \fg(b)$ in type
  $Sp(4a+2b+2)$,} all have relevant lowest 
$K-$types. In particular, the induced module is generated by spherically
relevant $K-$types. Precisely,

\noindent\textbf{Type B:\ } the representation is generated by the
$K-$types of the form $\mu_e,$

\noindent\textbf{Type C:\ } the representation is generated by the
$K-$types of the form $\mu_o,$

\noindent\textbf{Type D:\ } the representation is generated by the
{spherical $K-$-type $\mu_e(0)=\mu_o(0)$.} 
\end{proposition}

\subsection{}\label{sec:10.3} We now prove the irreducibility result
 mentioned at the beginning of the section in the case of $\fk g$ of type B; the
 other cases are similar. {Recall that we denote
 $\vO=(2x_0,2x_1,\dots,2x_{2m})$.} Let ${\vO_1}$ be the nilpotent orbit where we
 have removed {two  entries
   equal to $2a$ from $\vO$.} Let $\fk m:=gl(2a)\times
 \fg(n-2a).$ Then $L(\chi_\vO)$ is the  spherical subquotient of the
 induced representation 
 \begin{equation}
   \label{eq:10.3.1}
I(a,L(\chi_{\vO_1})):= 
Ind_{\fk m}^{\fg} [(-a+1/2,\dots ,a-1/2)\otimes L(\chi_{\vO_1})].
 \end{equation}
It is enough to show that if a parameter is unipotent, and satisfies 
$x_{i-1}=x_{i}=x_{i+1}=a,$ then $I(a,L(\chi_{\vO_1}))$ is generated by
its $K-$types  $\mu_e.$ This is because by  theorem \ref{sec:5.3}, the
$K-$types of type $\mu_e$ in (\ref{eq:10.3.1}) occur with full multiplicity in
the spherical irreducible subquotient, and the module is
unitary. 
{We will use various irreduciblity results stated in section
\ref{sec:2}) without mention.
} 

\medskip  

First, we reduce to the case when there are no $0<x_j< a.$ Let $\nu$
be the dominant parameter of $L(\chi_\vO),$ and assume $i$ is the smallest
index so that $x_{i-1}=a.$ There is an intertwining
operator 
\begin{equation}
  \label{eq:10.3.2}
  X(\nu)\longrightarrow I(1/2,\dots ,x_0-1/2;\dots ;1/2,\dots ,
  x_{i-2}-1/2;{\nu'}) 
\end{equation}
where $I$ is induced from $gl(x_0)\times\dots \times gl(x_{i-2})\times
\fk g(n-\sum_{j<i-1} x_j)$ with characters on the $gl$'s corresponding to
the strings in (\ref{eq:10.3.2}), and the irreducible module
{$L(\nu')$} on $\fg(n-\sum_{j<i-1} x_j)$. 
 The intertwining operator is onto, and
thus 
the induced module is generated by its spherical vector. By the
induction hypothesis, the induced module from
$(-a+1/2,\dots ,a-1/2)\otimes L(\nu'')$  on  $gl(2a)\times \fk
g(n-\sum_{j\le i} x_j)$ is irreducible. {
The parameter $\nu''$ is such
that the induced module has infinitesimal character $\nu'$, and
contains $L(\nu')$ as its spherical subquotient.} But 
\begin{equation}
  \label{eq:10.3.3}
\begin{aligned}
&I(1/2,\dots ,x_0-1/2;\dots ;1/2,\dots , x_{i-2}-1/2;-a+1/2,\dots
,a-1/2;\nu'')\cong\\
&I(-a +1/2,\dots ,a-1/2;1/2,\dots ,x_0-1/2;\dots ;1/2,\dots ,
x_{i-2}-1/2,\nu'') 
\end{aligned}
\end{equation}
This module maps by an intertwining operator onto $I(a,L(\chi_{\vO_1})),$ so
\newline $I(a,L(\chi_{\vO_1}))$ is generated by its spherical vector.

\medskip 
So we have reduced to the case when  
\begin{equation}
\begin{aligned}\label{eq:10.3.4}
x_0=x_1=x_2=a,\qquad  \text{ or }\\
x_0=0<x_1=x_2=x_3=a.
\end{aligned} 
\end{equation}

\medskip
Suppose we are in the first case of (\ref{eq:10.3.4}) and $m=1$. The
infinitesimal character is 
\begin{equation*}
(a-1/2,a-1/2,a-1/2,\dots ,1/2,1/2,1/2),  
\end{equation*}
each coordinate occuring three times. The induced module
\begin{equation}
  \label{eq:10.3.5}
I(-a+1/2,\dots , a-1/2) 
\end{equation}
of $\fg(2a)$ is a direct sum of irreducible factors computed in section
\ref{sec:10.2}; in particular it is generated by  $K-$types of the form
$\mu_e(2a-k,k)$ (with $k=0,1$). Consider the module
\begin{equation}
  \label{eq:10.3.6}
I(a-1/2;\dots ; 1/2;-a +1/2,\dots ,a-1/2),
\end{equation}
induced from characters on $GL(1)\times\dots\times GL(1)\times GL(2a).$
It is a direct sum of induced modules from the two factors of
(\ref{eq:10.3.5}). Each such induced module is a homomorphic image of
the corresponding standard module with dominant parameter. So
(\ref{eq:10.3.6}) is also generated by its $\mu_e$ isotypic components.   
But then
\begin{equation}
  \label{eq:10.3.7}
\begin{aligned}
&I(a-1/2;\dots ; 1/2;-a +1/2,\dots ,a-1/2)\cong \\
&I(-a +1/2,\dots ,a-1/2;a-1/2;\dots ; 1/2) 
\end{aligned}
\end{equation}
so the latter is also generated by its $\mu_e$ isotypic components.
Finally, the intertwining operator
\begin{equation}
  \label{eq:10.3.8}
I(a-1/2;\dots ;1/2)\longrightarrow I(1/2,\dots ,a-1/2)
\end{equation}
is onto, and the image of the intertwining operator
\begin{equation}
  \label{eq:10.3.9}
I(1/2,\dots ,a-1/2)\longrightarrow I(-a+1/2,\dots ,-1/2)
\end{equation}
is onto $L(-a+1/2,\dots ,-1/2).$ Thus the module induced from
$gl(2a)\times \fg(a)$, 
\begin{equation}
  \label{eq:10.3.10}
I(-a+1/2,\dots ,a-1/2;L(-a+1/2,\dots ,-1/2)),  
\end{equation}
is generated by its $\mu_e$ isotypic components. Since the
multiplicity of these $K-$types in (\ref{eq:10.3.10}) is the same as
in the irreducible spherical module, it follows that they must be equal.

\medskip
Now suppose that we are in the first case of (\ref{eq:10.3.4}) and
$m>1$, or in the second case, and $m>2$.  
The parameter has another $x_{2m-1}\le x_{2m}.$ 
We use an argument similar to the one above to show that the module
\begin{equation}
  \label{eq:10.3.11}
I(-x_{2m-1}+1/2,\dots ,x_{2m}-1/2,L(\chi_{\vO_2})),  
\end{equation}
{induced from $gl(x_{2m-1}+x_{2m})\times
  \fg(n-x_{2m-1}-x_{2m})$,} where $\vO_2$ is  the nilpotent orbit with
partition obtained from 
$\vO$ by removing $2x_{2m-1}, 2x_{2m},$ 
is generated by its $\mu_e$ isotypic components. The claim then
follows because the induced module is a homomorphic image of 
(\ref{eq:10.3.11}).  Precisely, $X(\nu)$ maps onto 
\begin{equation}
  \label{eq:10.3.12}
  \begin{aligned}
I(x_{2m-1}+1/2,\dots ,&x_{2m}-1/2;1/2,\dots ,x_0-1/2;\dots ; 1/2,\dots
,x_{2m-2}-1/2;\\
&L(-x_{2m-1}+1/2,-x_{2m-1}+1/2,\dots ,-1/2,-1/2))    
  \end{aligned}
\end{equation}
So this module is generated by its spherical vector. 
Replace $L(-x_{2m-1}+1/2,-x_{2m-1}+1/2,\dots ,-1/2,-1/2)$ by 
$I(-x_{2m-1}+1/2,\dots , x_{2m-1}-1/2).$ The ensuing module is a
direct sum of two induced modules by section \ref{sec:10.2}. 
They are both homomorphic images of standard modules, so generated by
their lowest $K-$types, which are of type $\mu_e$. Next observe that the map
\begin{equation}
  \label{eq:10.3.13}
 \begin{aligned}
I(x_{2m-1}+1/2,\dots ,x_{2m}-1/2;&1/2,\dots ,x_0-1/2;\dots ; 1/2,\dots
,x_{2m-2}-1/2;\\
&-x_{2m-1}+1/2,\dots ,x_{2m-1}-1/2) \longrightarrow \\
I(-x_{2m-1}+1/2,\dots x_{2m}-1/2;&1/2,\dots ,x_0-1/2;\dots ; 1/2,\dots
,x_{2m-2}-1/2)    
\end{aligned}  
\end{equation}
is onto. So the target module is generated by its $\mu_e$ isotypic
components.  The module
\begin{equation}
  \label{eq:10.3.14}
I(1/2,\dots ,x_0-1/2;\dots ; 1/2,\dots,x_{2m-2}-1/2) 
\end{equation}
(the string $(-x_{2m-1}+1/2,\dots ,x_{2m}-1/2)$ has been removed) 
has $L(-x_{2m-2}+1/2,\dots ,1/2)$ as its unique irreducible quotient,
because it is the homomorphic image of {an $X(\wti\nu)$ with $\wti\nu$
dominant.}
Therefore it is generated by its spherical vector.
Combining this with the induction assumption, we conclude that
\begin{equation}
  \label{eq:10.3.15}
I(-x_{2m-1}+1/2,\dots ,x_{2m}-1/2;-a+1/2,\dots ,a-1/2;L(\cCO_3))  
\end{equation}
is generated by its $\mu_e$ isotypic components. {Here
  $\cCO_3$ is obtained from $\cCO_2$ by removing two entries equal to
  $2a$.} It is isomorphic to 
\begin{equation}
  \label{eq:10.3.16}
 I(-a+1/2,\dots ,a-1/2;-x_{2m-1}+1/2,\dots ,x_{2m}-1/2;L(\cCO_3)).
\end{equation}
Finally, the multiplicities of the $\mu_e$ isotypic components of
$I(-x_{2m-1}+1/2,\dots ,x_{2m}-1/2;L(\cCO_3))$ are the
same as for the irreducible subquotient $L(\vO_1)$. This completes
the proof of the claim in this case.

\medskip
Remains to consider the case when $m=2$ and $x_0=0< x_1=x_2=x_3=a\le x_4.$ In
this case, the module
\begin{equation}
  \label{eq:10.3.17}
I(a+1/2,\dots ,x_4-1/2;-a+1/2,\dots ,a-1/2;-a+1/2,\dots ,a-1/2)
\end{equation}
is generated by its $\mu_e$ isotypic components because of proposition
\ref{p:10.2}, and arguments similar to the above. Therefore the same holds for 
\begin{equation}
  \label{eq:10.3.18}
I(-a+1/2,\dots ,x_4-1/2;-a+1/2,\dots ,a-1/2), 
\end{equation}
which is a homomorphic image via the intertwining operator which
interchanges the first two strings. But this is isomorphic to
\begin{equation}
  \label{eq:10.3.19}
I(-a+1/2,\dots , a-1/2,-a+1/2,\dots , x_4-1/2).
\end{equation}
Then $I(-a+1/2,\dots ,a-1/2,L(-x_4+1/2,\dots ,-1/2,-1/2)$ is a
homomorphic image of (\ref{eq:10.3.19}) so it is generated by its
$\mu_e$ isotypic components. By section \ref{sec:5.3}, the
multiplicities of the $\mu_e$ isotypic components are the same in
 $I(-a+1/2,\dots ,a-1/2,L(-x_4+1/2,\dots ,-1/2,-1/2))$ as in
 $L(\chi_\vO).$ 

This completes the proof of theorem \ref{t:10.1}.\qed

\ifx\undefined\bysame
\newcommand{\bysame}{\leavevmode\hbox to3em{\hrulefill}\,}
\fi

\end{document}